\documentclass[10pt]{article}
\usepackage{hyperref}
\usepackage{amssymb}
\usepackage{amsmath}

\usepackage{amsmath}
\usepackage{latexsym}
\usepackage{amssymb}

\usepackage{amsthm}

\newcommand{\mlabel}[1]{\marginpar{#1}\label{#1}}

\newcommand{\g}{{\mathfrak g}}

\newcommand{\h}{{\mathfrak h}}

\newcommand{\fa}{{\mathfrak a}}

\newcommand{\fe}{{\mathfrak e}}
\newcommand{\ff}{{\mathfrak f}}
\newcommand{\fg}{{\mathfrak g}}
\newcommand{\fh}{{\mathfrak h}}

\newcommand{\fk}{{\mathfrak k}}
\newcommand{\fl}{{\mathfrak l}}

\newcommand{\fn}{{\mathfrak n}}

\newcommand{\fq}{{\mathfrak q}}
\newcommand{\fp}{{\mathfrak p}}

\newcommand{\fs}{{\mathfrak s}}
\newcommand{\ft}{{\mathfrak t}}
\newcommand{\fu}{{\mathfrak u}}

\newcommand{\fz}{{\mathfrak z}}

\renewcommand\sp{\mathfrak {sp}}

\renewcommand{\:}{\colon}
\newcommand{\1}{\mathbf{1}}

\newcommand{\cB}{\mathcal{B}}

\newcommand{\cD}{\mathcal{D}}
\newcommand{\cE}{\mathcal{E}}

\newcommand{\cH}{\mathcal{H}}

\newcommand{\cM}{\mathcal{M}}

\newcommand{\cO}{\mathcal{O}}

\newcommand{\cS}{\mathcal{S}}
\newcommand{\cT}{\mathcal{T}}

\newcommand{\cV}{\mathcal{V}}
\newcommand{\cW}{\mathcal{W}}

\newcommand\bx{{\bf{x}}}
\newcommand\by{{\bf{y}}}

\newcommand{\bO}{\mathbf{O}}

\newcommand{\eset}{\emptyset}

\newcommand{\subeq}{\subseteq}
\newcommand{\supeq}{\supseteq}

\newcommand{\into}{\hookrightarrow}

\newcommand{\N}{{\mathbb N}}
\newcommand{\Z}{{\mathbb Z}}
\newcommand{\R}{{\mathbb R}}
\newcommand{\C}{{\mathbb C}}

\newcommand{\bP}{{\mathbb P}}

\renewcommand{\H}{{\mathbb H}}

\newcommand{\bH}{{\mathbb H}}
\newcommand{\bS}{{\mathbb S}}

\renewcommand{\tilde}{\widetilde}

\renewcommand{\L}{\mathop{\bf L{}}\nolimits}

\newcommand{\GL}{\mathop{{\rm GL}}\nolimits}
\newcommand{\SL}{\mathop{{\rm SL}}\nolimits}

\newcommand{\PGL}{\mathop{{\rm PGL}}\nolimits}

\newcommand{\PSL}{\mathop{{\rm PSL}}\nolimits}
\newcommand{\SO}{\mathop{{\rm SO}}\nolimits}

\newcommand{\U}{\mathop{\rm U{}}\nolimits}

\newcommand{\Sym}{\mathop{{\rm Sym}}\nolimits}

\newcommand{\Skew}{\mathop{{\rm Skew}}\nolimits}

\newcommand{\gl}  {\mathop{{\mathfrak{gl} }}\nolimits}

\newcommand{\fsl} {\mathop{{\mathfrak{sl} }}\nolimits}

\newcommand{\su}  {\mathop{{\mathfrak{su} }}\nolimits}
\newcommand{\so}  {\mathop{{\mathfrak{so} }}\nolimits}

\newcommand{\Exp}{\mathop{{\rm Exp}}\nolimits}

\newcommand{\ad}{\mathop{{\rm ad}}\nolimits}
\newcommand{\Ad}{\mathop{{\rm Ad}}\nolimits}

\newcommand{\tr}{\mathop{{\rm tr}}\nolimits}

\newcommand{\Herm}{\mathop{{\rm Herm}}\nolimits}
\newcommand{\Aherm}{\mathop{{\rm Aherm}}\nolimits}

\renewcommand{\>}{>\!\!>}

\newcommand{\Aut}{\mathop{{\rm Aut}}\nolimits}

\newcommand{\End}{\mathop{{\rm End}}\nolimits}
\newcommand{\id}{\mathop{{\rm id}}\nolimits}
\newcommand{\rk}{\mathop{{\rm rank}}\nolimits}

\renewcommand{\dim}{\mathop{{\rm dim}}\nolimits}

\newcommand{\Inn}{\mathop{{\rm Inn}}\nolimits}
\newcommand{\Int}{\mathop{{\rm int}}\nolimits}

\newcommand{\cone}{\mathop{{\rm cone}}\nolimits}

\newcommand{\sgn}{\mathop{{\rm sgn}}\nolimits}

\newcommand{\Spann}{\mathop{{\rm span}}\nolimits}

\newcommand{\dS}{\mathop{{\rm dS}}\nolimits}

\renewcommand{\phi}{\varphi}

\newcommand{\Rarrow}{\Rightarrow}
\newcommand{\nin}{\noindent} 
\newcommand{\oline}{\overline}

\newcommand{\la}{\langle}
\newcommand{\ra}{\rangle}

\newcommand{\up}{\mathop{\uparrow}}
\newcommand{\down}{\mathop{\downarrow}}
\newcommand{\res}{\vert}

\newcommand{\Spec}{{\rm Spec}}

\newcommand{\ssssarr}{\hbox to 15pt{\rightarrowfill}}
\newcommand{\sssarr}{\hbox to 20pt{\rightarrowfill}}
\newcommand{\ssarr}{\hbox to 30pt{\rightarrowfill}}
\newcommand{\sarr}{\hbox to 40pt{\rightarrowfill}}
\newcommand{\arr}{\hbox to 60pt{\rightarrowfill}}
\newcommand{\larr}{\hbox to 60pt{\leftarrowfill}}
\newcommand{\Arr}{\hbox to 80pt{\rightarrowfill}}

\def\theoremname{Theorem}
\def\propositionname{Proposition}
\def\corollaryname{Corollary}
\def\lemmaname{Lemma}
\def\remarkname{Remark}
\def\conjecturename{Conjecture} 

\def\definitionname{Definition}
\def\exercisename{Exercise}
\def\examplename{Example}
\def\examplesname{Examples}
\def\problemname{Problem}
\def\problemsname{Problems}

\def\satzname{Satz} 
\def\koroname{Korollar}
\def\folgname{Folgerung}
\def\bemerkname{Bemerkung}
\def\aufgname{Aufgabe}

\def\beisname{Beispiel}
\def\beissname{Beispiele}
\def\bewname{Beweis}

\def\@thmcounter#1{\noexpand\arabic{#1}}
\def\@thmcountersep{}
\def\@begintheorem#1#2{\it \trivlist \item[\hskip 
\labelsep{\bf #1\ #2.\quad}]}
\def\@opargbegintheorem#1#2#3{\it \trivlist
      \item[\hskip \labelsep{\bf #1\ #2.\quad{\rm #3}}]}
\makeatother
\newtheorem{theor}{\theoremname}[section]
\newtheorem{propo}[theor]{\propositionname}
\newtheorem{coro}[theor]{\corollaryname}
\newtheorem{lemm}[theor]{\lemmaname}

\newenvironment{thm}{\begin{theor}\it}{\end{theor}}

\newenvironment{prop}{\begin{propo}\it}{\end{propo}}

\newenvironment{cor}{\begin{coro}\it}{\end{coro}}

\newenvironment{lem}{\begin{lemm}\it}{\end{lemm}}

\newtheorem{rema}[theor]{\remarkname}

\newenvironment{rem}{\begin{rema}\rm}{\end{rema}}

\newtheorem{stepnow}[theor]{}

\newtheorem{defin}[theor]{\definitionname} 

\newenvironment{defn}{\begin{defin}\rm}{\end{defin}}

\newtheorem{exerc}{\exercisename}[section]

\newtheorem{exa}[theor]{\examplename}

\newenvironment{ex}{\begin{exa}\rm}{\end{exa}}

\newtheorem{exas}[theor]{\examplesname}

\newtheorem{conj}[theor]{\conjecturename}

\newtheorem{pro}[theor]{\problemname}

\newtheorem{prs}[theor]{\problemsname}

\newtheorem{aufg}{\aufgname}[section]

\newenvironment{prf}{\begin{proof}}{\end{proof}}
 
\newcommand{\pmat}[1]{\begin{pmatrix} #1 \end{pmatrix}}

{\hfill\qed\end{trivlist}}

\newenvironment{beweis*}{\begin{trivlist}\item[\hskip%
\labelsep{\bf\bewname.\quad}]}%
{\end{trivlist}}

\newtheorem{satzn}[theor]{\satzname}

\newtheorem{koro}[theor]{\koroname}

\newtheorem{folg}[theor]{\folgname}

\newtheorem{bem}[theor]{\bemerkname}

\newtheorem{aufgn}[theor]{\aufgname}

\newtheorem{beis}[theor]{\beisname}

\newtheorem{beiss}[theor]{\beissname}

\numberwithin{equation}{section}

\usepackage{color}

\renewcommand{\down}{{\mathop{\downarrow}}}

\newcommand\be{{\bf{e}}}

\renewcommand\up{{\uparrow}}

\renewcommand{\rk}{\mathop{{\rm rk}}\nolimits}
\newcommand{\oconv}{\mathop{{\rm oconv}}\nolimits}

\renewcommand{\bO}{\mathbb O}

\newcommand{\comp}{\mathop{{\rm comp}}\nolimits}

\renewcommand{\phi}{\varphi} 

\newcommand{\AdS}{\mathop{{\rm AdS}}\nolimits}

\renewcommand\mlabel{\label}

\begin{document}
\title{Modular geodesics and wedge domains \\ 
in non-compactly causal symmetric spaces \\ 
{\tt modgeo.tex}} 

\author{Vincenzo Morinelli, Karl-Hermann Neeb, Gestur \'Olafsson\footnote{
 The research of V. Morinelli VM was 
partially supported by a  Humboldt Research Fellowship for Experienced Researchers; the University of Rome through the MIUR Excellence Department Project 2023-2027, the  ``Tor Vergata'' CUP E83C23000330006 and ``Tor Vergata''   ``Beyond Borders'' CUP E84I19002200005,  Fondi di Ricerca Scientifica d'Ateneo 2021, OAQM, CUP E83C22001800005, and the European Research Council Advanced Grant 669240 QUEST. The research of K.-H. Neeb was partially supported by DFG-grants NE 413/10-1 and NE 413/10-2. The research of G. \'Olafsson was partially
supported by Simons grant 586106.}}

\maketitle

\bigskip

\begin{abstract}
  We continue our investigation of
  the interplay between causal structures
  on symmetric spaces and geometric aspects
  of Algebraic Quantum Field Theory. 
  We adopt the perspective that the geometric implementation
  of the modular group is given by the flow generated by
  an Euler element of the Lie algebra
  (an element defining a $3$-grading).
  Since   any Euler element of a semisimple Lie algebra specifies a canonical
  non-compactly causal symmetric space $M = G/H$, we turn in this paper
  to the geometry of this flow. Our main results
  concern the positivity region $W$ of the flow
  (the corresponding wedge region): If $G$ has trivial center,
  then $W$ is connected,  it coincides with
  the so-called observer domain, specified by a trajectory
  of the modular flow which at the same time is a causal geodesic,
  it can also be characterized in terms
  of a geometric KMS condition, and it has a natural structure
  of an equivariant fiber bundle over a Riemannian symmetric space
  that exhibits it as a real form of the crown domain of $G/K$.
  Among the tools that we need for these results are two observations
  of independent interest: a polar decomposition of the positivity domain
  and a convexity theorem for $G$-translates of open $H$-orbits
  in the minimal flag manifold specified by the $3$-grading.
\end{abstract}

\tableofcontents 
 
\section{Introduction} 
\mlabel{sec:1}

A new Lie theoretical approach to localization on  spacetimes involved in Algebraic Quantum Field Theory (AQFT) has been introduced in the recent years
  by the authors and collaborators in a series of works,
  see  \cite{MN21,NO17,NO21,NO23a,NO23b,MNO23, Oeh22}.  In the current paper we continue the investigation of the structure of wedge regions in non-compactly causal symmetric spaces, started in 
  \cite{NO23b}. First we briefly recall the motivation form AQFT,
  then we introduce tools and details to formulate our results.

Symmetric spaces are quotients $M = G/H$, where
$G$ is a Lie group, $\tau$ is an involutive automorphism of $G$
and $H \subeq G^\tau$ is an open subgroup
(cf.\ \cite{Lo69}). A causal symmetric space carries
a $G$-invariant field 
of pointed generating closed convex cones $C_m \subeq T_m(M)$
in their tangent spaces.
Typical examples are de Sitter space
$\dS^d \cong \SO_{1,d}(\R)_e/\SO_{1,d-1}(\R)_e$ and
anti-de Sitter space $\AdS^d \cong \SO_{2,d-1}(\R)_e/\SO_{1,d-1}(\R)_e$. 
These are Lorentzian, but we do not require
our causal structure to come from a Lorentzian metric,
which creates much more flexibility and a richer variety
of geometries.
Causal symmetric spaces permit to
study causality aspects of spacetimes in a highly symmetric environment.
Here we shall always assume that $M$ is
{\it non-compactly causal}
in the sense that the causal curves define a global order structure
with compact order intervals (they are called globally hyperbolic)
and in this context one can also prove the existence of a global
``time function'' with group theoretic methods (see \cite{Ne91}). 
We refer to the monograph \cite{HO97}
for more details and a complete exposition
of the classification of irreducible causal symmetric spaces.

Recent interest in causal symmetric spaces in relation with 
representation theory arose from their role as analogs of
spacetime manifolds in the context of
Algebraic Quantum Field Theory 
in the sense of Haag--Kastler.
A model in AQFT is specified by a {\it net} of von Neumann algebras
$\cM(\cO)$ acting on a fixed Hilbert space indexed by open subsets $\cO$
of the chosen spacetime $M$ (\cite{Ha96}). 
The hermitian elements of the algebra $\cM(\cO)$ represent
observables  that can be measured in the ``laboratory'' $\cO$.
These nets are supposed to satisfy fundamental
quantum and relativistic assumptions:
\begin{itemize}
\item[\rm(I)] Isotony: 
  $\cO_1 \subeq \cO_2$ implies $\cM(\cO_1) \subeq \cM(\cO_2).$ 
\item[\rm(L)] Locality: 
  $\cO_1 \subeq \cO_2'$ implies $\cM(\cO_1) \subeq \cM(\cO_2)'$,
  where $\cO'$ is the ``causal complement'' of $\cO$, i.e., the maximal
  open subset that cannot be connected to $\cO$ by causal curves.
\item[\rm(RS)] Reeh--Schlieder property: There exists a unit vector
  $\Omega\in \cH$ that is cyclic for $\cM(\cO)$ if $\cO \not=\eset$.
\item[\rm(Cov)] Covariance: 
  There is a Lie group $G$ acting on $M$ and a unitary representation
  $U \: G \to \U(\cH)$ such that 
  $U_g \cM(\cO) U_g^{-1} = \cM(g\cO)$ for $g \in G$. 
\item[\rm(BW)] Bisognano--Wichmann property: 
  $\Omega$ is separating for some ``wedge region'' $W \subeq M$
and there exists an element
$h \in \g$ with $\Delta^{-it/2\pi} = U(\exp th)$ for $t\in \R$,
where $\Delta$ is the modular operator
corresponding to $(\cM(W),\Omega)$
in the sense of the Tomita--Takesaki Theorem (\cite[Thm.~2.5.14]{BR87}). 
\item[\rm(Vac)] Invariance of the vacuum: $U(g)\Omega = \Omega$   for every $g \in G$. 
\end{itemize}
The (BW) property gives a geometrical meaning to the dynamics
provided by the  modular group $(\Delta^{it})_{t \in \R}$
of the von Neumann algebra $\cM(W)$ 
associated to wedge regions with respect to the vacuum state
specified by $\Omega$. On Minkowski/de Sitter
spacetime it provides an identification of the one-parameter group
$(\Lambda_W(t))_{t \in \R}$ of
 boosts in the Poincar\'e/Lorentz group with the Tomita--Takesaki modular operator:
\[ U(\Lambda_W(t))=\Delta^{-it/2\pi}. \] 
Here $\Lambda_{W}=g\Lambda_{W_1}g^{-1}$ 
is a one-parameter group of boosts associated with $W = g.W_1$, where
$W_1=\{x\in M :|x_0|<x_1\}$ is the standard right wedge and
\[ \Lambda_{W_1}(t) =(\cosh(t)x_0+\sinh(t)x_1,\cosh(t)x_1+\sinh(t)x_0,x_2,\ldots,x_d)\]
describes the boosts associated to $W_1$.

The homogeneous spacetimes occurring naturally in AQFT
  are causal symmetric spaces associated to their symmetry groups 
    (Minkowski spacetime for the Poincar\'e group,
  de Sitter space for the Lorentz group and
  anti-de Sitter space for $\SO_{2,d}(\R)$), and the localization
  in wedge regions is ruled by the acting group. 
  The rich interplay between the geometric and algebraic objects in AQFT allowed a 
  generalization of fundamental localization properties and the subsequent definition of 
  fundamental models (second quantization fields), having as initial data a general Lie group with distinguished elements
  (Euler elements) in the Lie algebra.
  Given an AQFT on Minkowski spacetime $M = \R^{1,d}$ (or  de Sitter spacetime $\dS^d
  \subeq \R^{1,d}$),
 the Bisognano--Wichmann (BW) property  allows an identification of geometric and
 algebraic objects in both free and interacting theories in all dimensions
 (\cite{DM20, Mu01, BW76}). This plays a central role in many results in AQFT and
 is a building block of our discussion.

Nets of von Neumann algebras can be
  constructed on causal symmetric spaces
  with representation theoretic methods.  We refer
to \cite{NO21} for left invariant nets on reductive Lie groups,
to \cite{Oeh21} for left invariant nets on non-reductive Lie groups,
to \cite{NO23a} for nets on compactly causal symmetric spaces,
and to \cite{FNO23} for nets on non-compactly causal symmetric spaces.
These papers construct so-called
  one-particle nets on symmetric spaces from which nets of von Neumann
  algebras can be obtained by second quantization functors (cf.~\cite{CSL22}). 

We know from \cite{MN23} 
that, in the general context, the potential generators $h \in \g$ of the modular groups
in (BW) are {\it Euler elements}, i.e.,
$\ad h$ defines a $3$-grading
\[ \g = \g_1(h) \oplus \g_0(h) \oplus \g_{-1}(h), 
\quad \mbox{ where } \quad \g_\lambda(h)
= \ker(\ad h - \lambda \1).\]
This leads to  the question how the existence and
the choice of the Euler element affects the geometry of
the associated symmetric space.
The (BW) property establishes a one-to-one correspondence between
    ``wedge regions'' $W \subeq M$ and the associated Euler elements. 
    So these {\it fundamental localization regions} can be determined
      in terms of Euler elements.
    This allowed the following generalization of nets of von Neumann 
    algebras on Minkowski/de Sitter spacetime: 
\begin{itemize}
\item     Given a Lie group $G$ with Lie algebra $\fg$, then 
the couples $(h,\tau_h)$, where $h \in \g$ is an Euler element and
  $\tau_h$ an involutive automorphism of $G$, inducing on $\g$ the involution
  $\tau_h = e^{\pi i \ad h}$, allow the definition of an
 ordered, $G$-covariant set of ``abstract wedge regions''
    carrying also some locality information \cite{MN21}.
In particular, they encode the commutation relation property of the Tomita operators 
(modular operator and modular conjugation). 
\item In the setting we are going to consider here, given an (anti-)unitary representation 
$(U,\cH)$ of $G \rtimes \{\1,\tau_h\}$,
  the bosonic second quantization net of von Neumann algebras associated to
  wedge regions can be constructed on the associated
  Fock space and, in this sense, it defines a generalized AQFT.   
  Here the net is constructed from the local one-particle spaces in $\cH$ 
  both abstract-algebraic, as in \cite{MN21,BGL02}, and through distribution vectors, as in 
  \cite{NO21,NO23a, FNO23}. In general, in particular when $Z(G)\neq \{e\}$, one may 
  consider other second quantization schemes as in \cite{GL95, CSL22}. This has not yet been done, see also \cite{MN21}. 
\item Causal symmetric spaces provide manifolds and a causal structure
  supporting nets of algebras. Here the wedge regions can be defined as open
  subsets in several ways.
  The equivalence of various characterizations
  has been shown  in \cite{NO23a,NO23b};
  see also the discussion below. 
\end{itemize}
The whole picture complies with Minkowski, de Sitter and anti-de Sitter spacetimes
and the associated free fields. 
A generalization of wedge regions of the Minkowski or de Sitter spacetime on general curved spacetimes have been proposed by many authors, see for instance \cite{DLMM11} and references therein. In our framework, on non-compactly causal symmetric spaces,
the rich geometric symmetries allow different characterizations of wedge regions, in particular in terms of \textit{positivity of the modular flow},
or  \textit{geometric KMS conditions} and in terms of \textit{polar decompositions}
as described in \cite{NO23b}. Some of them directly accord with the
literature, for instance for positivity of the modular flow, see \cite[Defin.~3.1]{CLRR22} 
and in particular \cite{NPT96} for the connection to thermodynamics
on de Sitter space. 
To see how these definitions apply to wedges in de Sitter space,
see \cite[App.~D.3]{NO23b} and \cite{BM96}.
On the other hand, the definition of wedge regions investigated in
\cite{NO23b,MNO23} coincides, up to specification of a  connected component
(cf.~\cite[Thm.~7.1]{NO23b}).
In Theorem~\ref{thm:connwedgedom} we prove that the
identification is actually complete for the adjoint groups 
since the wedge region defined in term of positivity of the modular flow
is connected. 
This contrasts the situation for compactly causal symmetric spaces,
where wedges are in general not connected,
cf. anti-de Sitter spacetime \cite[Lemma~11.2]{NO23b}.

Before we turn to the detailed formulation of our results,
we recall some basic terminology concerning symmetric Lie algebras 
(see \cite{NO23b} for more details). 

\begin{itemize}
\item A  {\it symmetric Lie algebra}
is a pair $(\g,\tau)$, where $\g$ is a finite-dimensional real Lie algebra 
and $\tau$ an involutive automorphism of~$\g$.
We write $\fh = \g^\tau = \ker(\1 - \tau)$ and 
$\fq = \g^{-\tau} = \ker(\1 + \tau)$ for the $\tau$-eigenspaces. 
\item A {\it causal symmetric Lie algebra}
is a triple $(\g,\tau,C)$, where $(\g,\tau)$ is a symmetric Lie algebra 
and $C \subeq \fq$ is a pointed generating 
closed convex cone invariant under the group
$\Inn_\g(\fh) = \la e^{\ad \fh}\ra$ acting in $\fq$. 
We call $(\g,\tau,C)$ {\it compactly causal~(cc)} if 
$C$ is {\it elliptic} in the sense that, for $x \in C^\circ$ (the interior of $C$
in $\fq$), the 
operator $\ad x$ is semisimple with purely imaginary spectrum. 
We call $(\g,\tau,C)$ {\it non-compactly causal (ncc)} if 
$C$ is {\it hyperbolic} in the sense that, for $x \in C^\circ$, the 
operator $\ad x$ is diagonalizable. 
\end{itemize}

As explained in detail in \cite{MNO23}, Euler elements
in reductive Lie algebras $\g$ lead naturally
to ncc symmetric Lie algebras:
For an Euler element $h \in \g$, choose a Cartan involution
$\theta$ of $\g$ with $\fz(\g) \subeq \g^{-\theta}$ such that
$\theta(h) = -h$. Then $\tau_h := e^{\pi i \ad h}$ is an involutive automorphism
of $\g$ commuting with $\theta$, so that
$\tau := \tau_h \theta$ defines a symmetric Lie algebra
$(\g,\tau)$ and there exists a pointed generating $\Inn_\g(\fh)$-invariant 
hyperbolic cone $C$ with $h \in C^\circ$. Under the assumption that
$\fh = \g^\tau$ contains no non-zero ideal of $\g$, there is a unique
minimal cone $C_\fq^{\rm min}(h)$ with this property. It is generated by
the orbit $\Inn_\g(\fh)h \subeq \fq$.

  Let $(\g,\tau,C)$ be an ncc symmetric Lie algebra
  and $(G, \tau^G, H)$ a corresponding {\it symmetric Lie group,} i.e.,
  $G$ is a connected Lie group, $\tau^G$ an involutive automorphism of $G$
  integrating $\tau$, and $H \subeq G^{\tau^G}$ an open subgroup.
  If, in addition, $\Ad(H)C = C$, then we
  call the quadruple $(G,\tau^G, H,C)$ a {\it causal symmetric Lie group}.
On $M=G/H$ we then obtain the structure of a  {\it causal symmetric space},
specified by the $G$-invariant field of open convex cones 
\begin{equation}
  \label{eq:v+}
  V_+(gH) := g.C^\circ \subeq T_{gH}(M).
  \begin{footnote}
{Note that the cones $V_+(m)$ are open, whereas $C$ is closed.}   
  \end{footnote}
\end{equation}
We further assume that 
\[ S := H \exp(C) = \exp(C)H \subeq G \] 
is a closed subsemigroup for which the polar map $H \times C \to S, 
(h,x) \mapsto h \exp x$ is a homeomorphism.
Then 
\begin{equation}
  \label{eq:ordgh1}
 g_1 H \leq g_2 H \quad \mbox{ if } \quad g_2^{-1} g_1 \in S 
\end{equation}
defines on $M$ a partial order, called the {\it causal order on $M$}. 
According to Lawson's Theorem (\cite{La94} and Theorem~\ref{thm:lawson}),
this is always the case if $\fz(\g) \subeq \fq$ and 
$\exp\res_{\fz(\g)}$ is injective. 
The second condition is always satisfied if $G$ is simply connected.

For an Euler element $h \in \g$ we consider the associated 
{\it modular flow} on $M = G/H$, defined by 
\begin{equation}
  \label{eq:modflow1}
  \alpha_t(gH) = \exp(th)gH.
\end{equation}
We study orbits of this flow which are geodesics 
$\gamma \: \R \to M$ with respect to the symmetric space structure 
and causal in the sense that $\gamma'(t) \in V_+(\gamma(t))$ for $t \in \R$.
We call them {\it $h$-modular geodesics}. All these are contained
in the  {\it positivity domain} 
\begin{equation}
  \label{eq:1}
 W_M^+(h) := \{ m \in M \: X^M_h(m) \in V_+(m)\} 
\end{equation}
of the vector field $X^M_h$ generating the modular flow. 
We refer to \cite{NO23b} for a detailed analysis of the latter domain in the 
special situations where the modular flow on $M$ has fixed points,
which is equivalent to the adjoint orbit $\cO_h = \Inn(\g)h$ intersecting~$\fh$.

We show for ncc symmetric Lie algebras, which are direct sums
of irreducible ones, that: 
\begin{itemize}
\item Causal modular geodesics exist if and only if the adjoint orbit 
  $\cO_h = \Ad(G)h\subeq \g$ intersects the interior of the cone
  $C \subeq \fq$, and then the centralizer $G^h
  = \{ g \in G \: \Ad(g)h = h\}$ of $h$ 
  acts transitively on the union of the corresponding curves
  (Proposition~\ref{prop:3.1}(c)).
\item Suppose that the cone is maximal, i.e., $C_\fq = C_\fq^{\rm max}$
  (see \eqref{eq:cmax-gen} and \cite[\S 3.5.2]{MNO23} for details).
  Let $\fq_\fk= \fq \cap \fk$ for
  a Cartan decomposition $\g = \fk \oplus \fp$ with~$h \in \fq_\fp$ and
  consider the domain
  \[   \Omega_{\fq_\fk}
    = \Big\{ x \in \fq_\fk \: \rho(\ad x) < \frac{\pi}{2}\Big\},\]
where $\rho(\ad x)$ is the spectral radius of $\ad x$.  
Then the connected component 
$W := W_M^+(h)_{eH}$ of the base point $eH$ in the positivity domain is
\[ W = G^h_e.\Exp_{eH}(\Omega_{\fq_\fk}) \] 
(Theorem~\ref{thm:pos-polar}).
\item We associate to any modular geodesic a connected open subset 
$W(\gamma) \subeq M$; the corresponding {\it observer domain}. 
For de Sitter space $\dS^d$, we thus obtain 
the familiar wedge domain obtained by intersecting 
$\dS^d$ with a Rindler wedge in Minkowski space~$\R^{1,d}$ 
(Example~\ref{ex:8.15}). 
In Theorem~\ref{thm:equality} we show that it coincides with~$W$,
provided that 
$H = K^h \exp(\fh_\fp)$ and $C_\fq = C_\fq^{\rm max}$. 
\item A key step in the proof of Theorem~\ref{thm:equality} is the 
following Convexity Theorem. 
Let 
\[ P^- := \exp(\g_{-1}(h)) G^h \subeq G \] 
 be the ``negative'' parabolic 
subgroup of $G$ specified by $h$ and identity 
$\g_1(h)$ with the open subset $\cB := \exp(\g_1(h)).eP^- \subeq G/P^-$. 
Then $\cD := H.0 \subeq \cB$ is an open convex subset, and for any 
$g \in G$ with $g.\cD \subeq \cB$, the subset $g.\cD \subeq \cB$ 
is convex (Theorem~\ref{thm:contheo}).
\item In Section~\ref{sec:7} we further show that, for
$C= C_\fq^{\rm  max}$ and $\g$ simple, that the
  real tube domain $\fh + C^\circ$ intersects the set
  $\cE(\g)$ of Euler elements in a connected subset (Theorem~\ref{thm:x.1}). 
  As a consequence, we derive that
  $W_M^+(h') \not=\eset$ if and only if $h' \in \cO_h$
  (Corollary~\ref{cor:7.3}). In particular, only one conjugacy class of
  Euler elements possesses non-empty positivity regions.   
\item In Theorem~\ref{thm:connwedgedom} we show 
  that the positivity domain $W_M^+(h)$ is connected for
  $G = \Inn(\g)$ and $\g$ simple,
  and this implies that
  \[ W(\gamma) = W = W_M^+(h).\] 
  From this in turn we derive that  the stabilizer group
    $G_W = \{ g \in G \: g.W = W\}$ coincides with $G^h$ 
  (Proposition~\ref{prop:GW}), so that the {\it wedge space}
  $\cW(M) := \{g.W \: g \in G\}$ of wedge regions in $M$ can be identified, 
  as a homogeneous $G$-space,  with the adjoint orbit
  $\cO_h= \Ad(G)h \cong G/G^h$.  In particular
    $\cW(M)$ also is a symmetric space. 
\item Finally, we show in Theorem~\ref{thm:kms} that
    $W$ coincides with the KMS wedge domain
\[   W^{\rm KMS} = \{  m \in M \:  \alpha_{it}(m) \in \Xi \mbox{ for }
  0 < t < \pi\}.\]
\end{itemize}

We conclude this introduction with some more motivation
  from AQFT. 
The analysis of the properties of the modular flow on symmetric spaces
   is also motivated by the investigation of energy inequalities in quantum and relativistic theories.
In General Relativity, there exist many solutions to the Einstein equation that, for various reasons, may not be physical.  Energy conditions such as the pointwise non-negativity of the energy density, which  ensures that the gravity force is attractive, can be required to discard non-physical models (\cite{Wal84,Few12}). In quantum and relativistic theories, the energy conditions need to be rewritten. For instance, it is well known that the energy density at individual spacetime points is unbounded from below, even if the energy density integrated over a Cauchy surface is non-negative (see \cite{FH05, Few12} and references therein).

Families of inequalities have been discussed in several models, employing different mathematical and physical approaches (see for instance \cite{Few12,KS20, FLPW20, Ve00, KLLSM18, MPV23}). In the recent years operator algebraic techniques have been very fruitful for the study of the energy inequalities because of the central role played by the modular hamiltonian in some of these energy conditions. This object corresponds to the logarithm of the modular group of a local algebra of a specific region, which in some cases can be identified with the generator of a one-parameter group of spacetime symmetries by the Bisognano--Wichmann property.  In this regard, we cite---without entering into details---the ANEC (Averaged Null Energy Condition) and the QNEC (Quantum Null Energy Condition) and their relation with the Araki relative entropy, an important quantum-information quantity, defined in terms of relative modular operators (see for instance \cite{MTW22,Lo20, CLRR22,Lo19, CLR20, Ara76,LX18,LM23}). We stress that, in this analysis, the study of the modular flow on the manifold can be particularly relevant. Moreover, in order to find regions where to prove energy inequalities, one may also need to deform the modular flow (\cite{MTW22, CF20}). 
In our abstract context the Euler element specifies the flow that can be
implemented by the modular operator, hence the modular Hamiltonian, when
the Bisognano--Wichmann property holds. In particular, the identification of specific flows on symmetric spaces (modular flows), the characterization in terms of modular operators of covariant local subspaces attached to specific regions (wedges) deeply
motivate an analysis of modular flows on non-compactly causal symmetric spaces
pursued in our project. 

In this respect, the wedge regions are the first fundamental open
subsets of spacetime to be studied in detail. 
Following General Relativity (see for instance \cite{CLRR22,DLMM11}
and references therein), one can define them as an open connected,
causally convex subregion $W$ of a spacetime $M$,
associated to a Killing flow $\Lambda$ preserving $W$, which is 
timelike and time-oriented on $W$. 
On Minkowski spacetime the flow $\Lambda$, a one-parameter group of boosts, corresponds to 
the time-evolution of a uniformly accelerated observer moving within $W$. Then $W$ is a  horizon for this observer:
he cannot send a signal outside $W$ and receive it back.  Then the vacuum state becomes a 
thermal state for the algebra of observables inside the wedge region
$W$ by the Bisognano--Wichmann property (\cite{GL03,Ha76,Lo97}). 
In our general context, we recover the definition (and equivalent ones)
of  wedge regions. Then, by the Bisognano--Wichmann property, 
the thermal property of the vacuum state holds
when nets of algebras/standard subspaces are considered (\cite{MNO23,NO23b,MN21}). In this 
paper we focus on the related properties of the wedge regions
in non-compactly causal symmetric spaces.\\

\nin {\bf Notation:} 
\begin{itemize}
\item If $M$ is a topological space and $m\in M$, then $M_m$ denotes the connected component of $M$ containing
$m$. In particular we write $e \in G$ for the identity element in the Lie group~$G$ 
 and $G_e$ for its identity component.

\item Involutive automorphisms of $G$ are typically denoted $\tau^G$,  
  and $\tau$ is the corresponding automorphism of the Lie algebra $\g = \L(G)$. 
 We write $\fg^\tau = \ker (\1-\tau)$ and $\fg^{-\tau} = \ker (\1 + \tau)$.
\item For $x \in \g$, we write $G^x := \{ g \in G \: \Ad(g)x = x \}$ 
for the stabilizer of $x$ in the adjoint representation 
and $G^x_e = (G^x)_e$ for its identity component. 
\item For $h \in \g$ and $\lambda \in \R$, we write 
$\g_\lambda(h) := \ker(\ad h - \lambda \1)$ for the corresponding eigenspace 
in the adjoint representation.
\item If $\g$ is a Lie algebra, we write $\cE(\g)$ for the set of 
{\it Euler elements} $h \in \g$, i.e., $\ad h$ is non-zero and diagonalizable 
with $\Spec(\ad h) \subeq \{-1,0,1\}$. The corresponding involution
is denoted $\tau_h = e^{\pi i \ad h}$.  
\item For a Lie subalgebra $\fs \subeq \g$, we write 
$\Inn_\g(\fs)= \la e^{\ad \fs} \ra \subeq \Aut(\g)$ for the subgroup 
generated by $e^{\ad \fs}$.
\item For a convex cone $C$ in a vector space $V$, we write
  $C^\circ := \Int_{C-C}(C)$ for the relative interior of $C$ in its span. 
\item We use the notation 
\begin{equation}
  \label{eq:specrad}
\rho(A) := \sup \{ |\lambda| \: \lambda \in \Spec(A) \}
\end{equation}
for the spectral radius of a linear operator~$A$. 
\end{itemize}

\nin{\bf Acknowledgment:}
We thank J.~Wolf for discussions concerning
  the Convexity Theorem~\ref{thm:contheo}.

\section{Causal Euler elements and ncc symmetric spaces} 
\mlabel{sec:2}

In this section we recall some basic results on
  Euler elements and their relation with non-compactly causal
  symmetric spaces. Most of these statements are discussed in
  detail in \cite{MNO23}.

 Recall from above that an Euler element in a Lie algebra $\fg$ is an element $h$ defining a $3$-grading of
 $\g$ by $\g = \g_{-1}\oplus \fg_0\oplus \g_{+1}$
 with $\g_j = \ker (\ad h - j\1)$, $j =-1,0,1$. We write
$\cE (\g)$ for the set of Euler elements in $\g$.
In this section we recall some results on 
from  \cite{MNO23} on Euler elements  that are crucially used in the following.

\begin{defn} Let $\g$ be a reductive Lie algebra.

  \nin (a) A {\it Cartan involution} of $\g$ is an involutive automorphism
  $\theta$ for which $\fz(\g) \subeq \g^{-\theta}$ and $\g^\theta$ is maximal
  compactly embedded in the commutator algebra $[\g,\g]$. We then write,  using
  the notation from the introduction,
  \[  \g = \fk \oplus \fp \quad \mbox{ with } \quad
    \fk = \g^\theta 
    \quad\mbox{ and }\quad \fp = \g^{-\theta} 
    \]

  \nin (b) If $\tau$ is another involution on $\g$ commuting with
  $\theta$, $\fh := \g^\tau$ and $\fq := \g^{-\tau}$, then we   have
  \[ \fh = \fh_\fk \oplus \fh_\fp, \quad
  \fq = \fq_\fk \oplus \fq_\fp \quad \mbox{ with } \quad
  \fh_\fk = \fh \cap \fk, \quad\fh_\fp = \fh \cap \fp,\quad 
  \fq_\fk = \fq \cap \fk, \quad\fq_\fp = \fq \cap \fp.\]

\nin (c) The {\it Cartan dual} of the symmetric Lie algebra $(\g,\tau)$
is  the symmetric Lie algebra $(\g^c, \tau^c)$ with 
\[ \g^c = \fh + i \fq \quad \mbox{ and } \quad 
\tau^c(x+ iy) = x - iy \quad \mbox{ for } \quad x \in \fh, y \in \fq.\] 
 Note that  $\g^c =(\fg_\C)^{\overline{\tau}}$
where $\oline\tau$ is the conjugate-linear extension of $\tau$ to $\g_\C$; 
in particular $\g^c$ is a real form of $\g_\C$.
\end{defn}

\begin{defn} Let $(\g,\tau)$ be a symmetric Lie algebra and
  $h \in \cE(\g) \cap \fq$. 
  We say that $h$ is {\it causal} if there exists an
  $\Inn_\g(\fh)$-invariant 
  closed pointed generating convex cone $C$ in~$\fq$
  with $h\in C^\circ$.
We write $\cE_c(\fq) \subeq \cE(\g) \cap \fq$ for the 
{\it set of causal Euler elements in $\fq$}.
Recall that the triple $(\g,\tau ,C)$ is ncc if $C$ is hyperbolic. 
\end{defn}

\begin{lem}  \mlabel{lem:2.5}
  Let $(\g,\tau,C)$ be a {\bf simple} ncc symmetric Lie algebra
  and $h \in \fq$ be a causal Euler element.
  Then the following assertions hold:
  \begin{itemize}
  \item[\rm(a)] There exist closed
  convex pointed generating $\Inn_\g(\fh)$-invariant cones
  \[ C_\fq^{\rm min}(h) \subeq C_\fq^{\rm max}(h) \]
  such that $h \in  C_\fq^{\rm min}(h)^\circ$ and either 
  \[ C_\fq^{\rm min}(h) \subeq C \subeq C_\fq^{\rm max}(h) \quad \mbox{ or } \quad 
     C_\fq^{\rm min}(h) \subeq -C \subeq C_\fq^{\rm max}(h).\]
 \item[\rm(b)] If $(G,\tau^G, H)$ is a connected
   symmetric Lie group with 
   symmetric Lie algebra $(\g,\tau)$, then
   two mutually exclusive cases occur:
\begin{itemize}
\item[$\bullet$] $\Ad(H)h = \Ad(H_e)h$ and $G/H$ is causal.
\item[$\bullet$] $-h \in \Ad(H)h$ and $G/H$ is not causal.
\end{itemize}
\end{itemize}
 \end{lem}

\begin{prf} (a) follows from \cite[Sect.~3.5.2]{MNO23}
  and (b) from \cite[Prop.~4.18]{MNO23}. 
\end{prf}

If $h$ is an  Euler element in the {\bf reductive} Lie algebra $\g$
and $\theta$ a Cartan involution with $\theta(h) = -h$,
$\tau := \theta \tau_h$ 
and $\fz(\g) \subeq \g^{-\theta}$, then \cite[Thm.~4.2]{MNO23}
implies  that there exists an $\Inn_\g(\fh)$-invariant pointed closed convex
cone $C \subeq \fq$ with $h\in C^\circ$,
so that $(\fg, \tau ,C)$ is ncc.
Further, all ideals of $\g$ contained in $\g^\tau = \fh$ are compact.
We have a decomposition
\begin{equation}
  \label{eq:gdeco}
  \g = \g_k \oplus \g_r \oplus \g_s,
\end{equation}
where $\g_s$ is the sum of all simple ideals not commuting with $h$
(the strictly ncc part), $\g_r$ is the sum of the center
$\fz(\g)$ and all non-compact simple
ideals commuting with $h$ on which $\tau = \theta$
(the non-compact Riemannian part),
and $\g_k$ is the sum of all simple compact ideals (they commute with~$h$).
All these ideals are invariant under $\theta$ and $\tau = \tau_h \theta$,
so that we obtain decompositions
\begin{equation}
  \label{eq:qs}
 \g_s = \fh_s\oplus \fq_s, \quad
  \g_r = \fh_r\oplus \fq_r \quad \mbox{ and } \quad
  \g_k = \fh_k,
\end{equation}
where $\fh_r \oplus \fh_k$ is a compact ideal of~$\fh$,
$\g_r = \fh_r\oplus \fq_r$ is a Cartan decomposition 
and $\fq_\fp=\fq_{\fp,s}\oplus \fq_r$.
In particular $\fq = \fq_s \oplus \fq_r$. 
Let $p_s: \fq \to \fq_s$ be the projection onto $\fq_s$
  with kernel $\fq_r$. Then
\cite[Prop. B.4]{MNO23} implies that
every $\Inn_\g (\fh)$-invariant closed convex cone $C$ satisfies 
\[C_s := p_s(C) = C\cap \fq_s\quad\text{and}\quad
  C_s^\circ = C^\circ \cap \fq_s .\]

By Lemma~\ref{lem:2.5}(a) we obtain a pointed $\Inn_\g(\fh)$-invariant cone
$C_{\fq_s}^{\rm min}(h) \subeq \fq_s$, adapted to the decomposition into irreducible
summands, whose dual cone 
$C_{\fq_s}^{\rm max}(h)$ with respect to the Cartan--Killing form
$\kappa(x,y) = \tr(\ad x \ad y)$ satisfies 
$C_{\fq_s}^{\rm min}(h) \subeq C_{\fq_s}^{\rm max}(h)$.
Put 
\begin{equation}
  \label{eq:cmax-gen}
 C_\fq^{\rm min}(h) := C_{\fq_s}^{\rm min}(h)  
  \subeq 
  C_\fq^{\rm max}(h) := C_{\fq_s}^{\rm max}(h)  + \fq_r.
\end{equation}
Both cones are adapted to the decomposition of $(\g,\tau)$ 
into irreducible summands.  
Further, each pointed generating $\Inn_\g(\fh)$-invariant cone $C$ containing
$h$ satisfies
\begin{equation}
  \label{eq:min-max-q2}
  C_\fq^{\rm min}(h) \subeq C \subeq C_\fq^{\rm max}(h).
\end{equation}
Here the first inclusion is obvious, and the second one follows from the 
fact that $h$ is also contained in the dual cone
\[ C^{\star}
= \{ y \in \fq \: (\forall x \in C) \ \kappa(x,y) \geq 0\}.\]
This leads to
$C_\fq^{\rm min}(h) \subeq C^\star$, and thus to
$C \subeq  C_\fq^{\rm min}(h)^\star = C_\fq^{\rm max}(h)$
(cf.\ \cite[\S 3.5]{MNO23} for more details). 

\begin{lem}\mlabel{lem:2.x} If $x \in (C_\fq^{\rm max})^\circ$, then
  the centralizer   $\fz_\fh(x) = \fh \cap  \ker(\ad x)$
  is compactly embedded in~$\g$, i.e., consists of elliptic elements.   
\end{lem}

\begin{prf} First we observe that
  the cone $C_\fq^{\rm max}$ is adapted to the decomposition
  $\g = (\g_k + \g_r) + \g_s$ and so is the centralizer of $x = x_r + x_s$
  in $\fh = (\g_\fk + \fk_r) + \fh_s$. Hence
  the assertion follows from the fact that
  $\g_k + \fk_r$ is compactly embedded and
  $\fz_{\fh_s}(x) = \fz_{\fh_s}(x_s)$ is compactly embedded because the cone
  $C_{\fq_s}^{\rm max}$ is pointed (\cite[Prop.~V.5.11]{Ne00}). 
\end{prf}

\begin{thm} {\rm(Uniqueness of the causal involution)}
  {\rm(\cite[Thm.~4.5]{MNO23})}
  \mlabel{thm:uniqueinv} 
Let $(\g,\tau,C)$ be a semisimple ncc symmetric Lie algebra 
for which all ideals of $\g$ contained in $\fh$ are compact,
$\g_s$ the sum of all non-Riemannian ideals, $\fq_s := \g_s \cap \fq$,
$C_s := C \cap \fq_s$, and $\theta$ a Cartan  involution commuting with~$\tau$. 
Then the following assertions hold: 
\begin{itemize}
\item[\rm(a)] $C_s^\circ \cap \fq_\fp$ contains a unique Euler element~$h$,
  and this Euler element satisfies $\tau = \tau_h \theta$. 
\item[\rm(b)] $\Inn_\g(\fh)$ acts transitively on $C_s^\circ \cap \cE(\g)$. 
\item[\rm(c)] For every Euler element $h \in C_s^\circ$, the involution 
$\tau\tau_h$ is Cartan. 
\end{itemize}
\end{thm}

  \begin{prop} \mlabel{prop:hhcompact} 
    Let $(G,\tau^G, H,C)$ be a connected semisimple ncc symmetric 
    Lie group for which $\fh = \g^\tau$ contains no non-compact ideal of $\g$ 
    {\rm(}$\g = \g_r + \g_s${\rm)} and let $h \in C_s^\circ$
    {\rm(cf.\ Theorem~\ref{thm:uniqueinv})} be a causal
    Euler element. Then
the following assertions hold: 
\begin{itemize}
\item[\rm(a)] $H = H_e H^h$,
i.e., every connected component of $H$ meets $H^h$.     
\item[\rm(b)] 
  $\Ad(H^h) = \Ad(H)^h$ is a maximal compact subgroup  of $\Ad(H)$. 
\item[\rm(c)] $\Ad(H)^h = \Ad(H)^{\tau_h}$ and  $\tau_h:= e^{\pi i \ad h}$
  induces a Cartan involution on $\Ad(H)$. 
\item[\rm(d)] $\tau$ induces a Cartan involution on $\Ad(H)^h$ for which 
$\Ad(H^h_e)^\tau = e^{\ad \fh_\fk}$ is connected. 
\end{itemize}
  \end{prop}

  \begin{prf} The statements on the adjoint group $\Ad(G) = \Inn(\g)$
    follows from \cite[Cor.~4.6]{MNO23}
    because $\Ad(H) \subeq \Inn(\g)^\tau$ preserves $C$.
    Further,   $\Ad(H)^h = \Ad(H^h)$ and $H^h = \Ad^{-1}(\Ad(H)^h)$
    imply with (a) (for $\Ad(G)$)
    that $H = H_e H^h$.
  \end{prf}

\begin{defn} If $\g$ is a simple hermitian Lie algebra,
  $\theta$ a Cartan involution of $\g$ and   $\fa \subeq \fp$ maximal
  abelian, then the restricted root system
  $\Sigma(\g,\fa)$ is either of type $C_r$ or $BC_r$.
  In the first case, we say that $\g$ is {\it of tube type}.
\end{defn}

Recall that if $(\g,\tau)$ is simple ncc, then either $\g^c$ is simple
hermitian or $\g^c \cong \fh_\C$, where $\fh = \g^\tau$ is simple hermitian
(\cite[Rem.~4.24]{MNO23}).

\begin{prop} \mlabel{prop:testing} {\rm(\cite[Lemma~5.1, Prop.~5.2]{MNO23})} 
Let $(\g,\tau,C)$ be a simple ncc symmetric Lie algebra. 
Pick a causal Euler element $h \in C^\circ$ and 
$\ft_\fq \subeq \fq_\fk$ maximal abelian and set
$s := \dim \ft_\fq$. 
Then the following assertions hold: 
\begin{itemize}
\item[\rm(a)] The Lie algebra $\fl$ generated by $h$ and $\ft_\fq$ 
is reductive. 
\item[\rm(b)] The  commutator algebra $[\fl,\fl]$ is
 isomorphic to $\fsl_2(\R)^s$ 
\item[\rm(c)] $\fz(\fl) = \R h_0$ for some 
hyperbolic element $h_0$ satisfying $\tau(h_0) =- h_0 = \theta(h_0)$ 
which is zero if and only if $\g^c$ is of tube type. 
\item[\rm(d)] The Lie algebra $\fl$ is $\tau$-invariant 
and $\fl^\tau \cong \so_{1,1}(\R)^s$. 
\item[\rm(e)] 
 For $x \in \ft_\fq$, we have 
$\rho(\ad x) = \rho(\ad x\res_{\fs}),$ 
where $\rho$ denotes the spectral radius. 
With the basis 
\[ z^j = \Big(0,\ldots, 0, 
\frac{1}{2}\pmat{0 & -1 \\ 1 & 0},0,\cdots, 0\Big), \quad j =1,\ldots, s, \]
in $\so_2(\R)^s$ we have for $x = \sum_{j = 1}^s x_j z^j$ 
\begin{equation}
  \label{eq:rhosl2}
 \rho(\ad x) = \max \{ |x_j| \:  j =  1,\ldots, s \}.
\end{equation}
\end{itemize}
\end{prop}

Note that (c) implies that $\fl$ is semisimple, i.e., $h \in [\fl,\fl]$,  
if and only if $\g^c$ is of tube type. 

  \begin{prop} \mlabel{prop:hk=kh} {\rm(\cite[Prop.~7.10]{MNO23})}
Let $(\g,\tau,C)$ 
be a semisimple modular non-compactly causal semisimple 
symmetric Lie algebra, 
where $\tau = \tau_h \theta$, $h \in \fq_\fp \cap C$ a causal Euler element,  
\[ G := \Inn_{\g_\C}(\g)\cong \Inn(\g), \quad K := G^\theta = \Inn_\g(\fk),\quad \mbox{ 
and }  \quad G^c := \Inn_{\g_\C}(\g^c).\]
 Then $H := G \cap G^c$ satisfies 
\[ H = K^h \exp(\fh_\fp) \quad \mbox{ and } \quad 
H \cap K = K^h.\] 
In particular $K^h \subeq K^{\tau_h} = K^\tau$ implies $H \subeq G^\tau$. 
\end{prop}

\section{The positivity domain and modular geodesics} 
\mlabel{sec:3}

Let $(G, \tau^G, H, C)$ be a connected {\bf semisimple}
causal symmetric Lie group 
with ncc symmetric Lie algebra $(\g,\tau, C)$.
We fix a causal Euler element $h \in C^\circ$
(Theorem~\ref{thm:uniqueinv}) and
write $M = G/H$ for the associated symmetric 
space.

One of  our goals in this paper
is to describe the structure of the  {\it positivity domain} 
\[  W_M^+(h) := \{ m \in M \: X^M_h(m) \in V_+(m)\} \] 
of the vector field $X^M_h$ generating the modular flow.
Our first major result is the identification of the connected component $W$ 
of the base point $eH$ in the positivity domain $W_M^+(h)$ as 
\begin{equation}
  \label{eq:W1}
  W := W_M^+(h)_{eH} =  G^h_e \Exp_{eH}(\Omega_{\fq_\fk})
\end{equation}
(Theorem~\ref{thm:pos-polar}).

Some of the results in this section 
had been obtained in \cite{NO23b} for the special case of
ncc symmetric Lie algebras for which $\fh$ contains an Euler 
element, whereas here 
we are dealing with general non-compactly causal symmetric 
Lie algebras. 

\subsection{Modular geodesics}

In this subsection we introduce the concept of an $h$-modular
geodesic in a non-compactly causal symmetric space
$M$ and discuss some of its immediate properties.
We also show that, in compactly causal spaces, non-trivial causal
modular geodesics do not exist.

\begin{defn} (Geodesics and causality) 
  Let  $M = G/H$ as above.
\begin{itemize}
\item We call a geodesic $\gamma \: \R \to M$ {\it causal} if 
  $\gamma'(t) \in V_+(\gamma(t))$ for every $t \in \R$
  (see \eqref{eq:v+}). 
\item Let $h \in \g$ be an Euler element. The flow on $M$ defined by 
\begin{equation}\label{eq:modflow} 
\alpha_t(gH) = \exp(th) gH  =g\exp (\Ad (g^{-1})h)H
\end{equation}
is called the {\it modular flow} (associated to $h$).
Its infinitesimal generator is denoted $X^M_h \in \cV(M)$. 
\item A geodesic $\gamma \: \R \to M$ is called {\it $h$-modular} if 
  $\gamma(t) = \alpha_t(\gamma(0))$ holds for all $t \in \R$, i.e.,
  $\gamma$ is an integral curve of $X_h^M$. 
\end{itemize}
\end{defn}

\begin{prop} \mlabel{prop:3.1}
Suppose that $(\g,\tau)$ is a direct sum of irreducible ncc symmetric Lie algebras
{\rm(}$\g = \g_s${\rm)}.   The following assertions hold
  for any Euler element ${h} \in \cE(\g)$ and the corresponding
  modular flow $\alpha_t(m) = \exp(t{h}).m$ on $M = G/H$:
  \begin{itemize}
\item[\rm(a)] The orbit under the modular flow 
is a causal geodesic if and only if $m$ is contained in 
\begin{equation}
  \label{eq:mhc2}
  M^{{h}}_C = \{ g H\in G/H \: \Ad(g)^{-1} {h} \in C^\circ\}.
\end{equation} 
\item[\rm(b)] All connected components of $M^{{h}}_C$
  are Riemannian symmetric space of non-compact type:
  For every $m \in M^{{h}}_C$, the exponential map 
\[ \Exp_m \: T_m(M_C^{{h}}) \to (M^{{h}})_m \] 
is a diffeomorphism.
\item[\rm(c)] ${h}$-modular causal geodesics exist 
  if and only if $\cO_{{h}} = \Ad(G){h}$ intersects~$C^\circ$.
In this case   $G^{{h}}$ acts transitively on $M_C^{{h}}$.
  \end{itemize}
\end{prop}

\begin{prf} (a) Assume first that $\Ad (g)^{-1} h\in \fq$.
Then \eqref{eq:modflow} implies that  
the orbit of $m = gH$ under the modular flow is 
a geodesic.
The  causality is by definition equivalent to
$\Ad(g)^{-1}{h}\in C^\circ$.

Suppose, conversely, that $t \mapsto \alpha_t(gH)$ is a 
causal geodesic. Lemma~\ref{lem:geodesic} implies 
that $\Ad(g)^{-1}{h} = x_\fh + x_\fq$, where $[x_\fh, x_\fq] = 0$ 
and $x_\fq \in C^\circ$. By Lemma~\ref{lem:2.x}, 
$x_\fh$ is elliptic and $x_\fq$ is hyperbolic because it is contained
in $C^\circ$. Therefore $\ad x_\fh + \ad x_\fq$ is the unique
Jordan decomposition
of $\ad x$ into elliptic and hyperbolic summand.
As $\Ad(g)^{-1}{h}$ is an Euler element, the elliptic summand vanishes,
and thus $\ad x_\fh = 0$, i.e., $x_\fh \in \fz(\g) \cap \fh = \{0\}$
(recall that $\fz(\g) \subeq \fq$). This shows that
$\Ad(g)^{-1} {h} \in \fq$, so that $gH \in M^{{h}}_C$.

\nin (b) Choosing $m$ as a base point, we
  may assume that $m = eH$, so that 
(a)  implies that $h \in C^\circ \subeq \fq$ is a causal Euler element. 
Pick a Cartan involution $\theta$ commuting 
$\tau$ which satisfies $\theta(h)  = -h$ 
(cf.\ \cite{KN96}), i.e., $h \in \fq_\fp$. 
Then  $\tau = \tau_{h} \theta$ 
follows from Theorem~\ref{thm:uniqueinv}(a).
As $(M^c_C)_m = G^h_e.m$ by Lemma~\ref{lem:mx}, the assertion
now follows from $\g^h =\fh_\fk \oplus \fq_\fp$.

\nin (c) The first assertion follows immediately from (a).
For the second assertion, suppose that $m_0 = g_0H \in M^{{h}}_C$. 
As $M^{{h}}_C$ is $G^{{h}}$-invariant, $G^{{h}}.m_0 \subeq M^{{h}}_C$. 
Let $h_c := \Ad(g_0)^{-1}{h}$, so that 
$\cE(\g) \cap C^\circ = \Inn_\g(\fh) h_c$
by Theorem~\ref{thm:uniqueinv}(b) (recall that $C = C_s$).
If $gH \in M^h_C$, i.e., 
\[ \Ad(g)^{-1} {{h}} \in C^\circ \cap \cE(\g)
\ {\buildrel\ref{thm:uniqueinv} \over =}\ \Inn_\g(\fh) h_c =
  \Inn_\g(\fh)\Ad(g_0)^{-1}h,\] 
then there exists an element $g_1 \in \Inn_\g(\fh)$ with 
$g g_1 g_0^{-1} \in G^{{h}}$, so that $g \in G^{{h}} g_0 g_1^{-1} \in G^{{h}} g_0 H$, 
and therefore $gH \in G^{{h}}.m_0$. 
\end{prf}

We record the following consequence of~\eqref{eq:modflow}:
\begin{lem} \mlabel{lem:wm+}
For any causal  Euler
  element $h \in C^\circ$, we have
  \[ W_M^+(h) = \{ gH \in G/H \: \Ad(g)^{-1}h \in \cT_C\}, \quad \mbox{ where }
    \quad
    \cT_C := \fh + C^\circ.\]
\end{lem}

Due to the hyperbolicity of Euler elements, modular causal geodesics
do not exist for compactly causal symmetric spaces: 

\begin{prop} If $M = G/H$ is a compactly causal symmetric space, 
then non-trivial causal modular geodesics do not exist. 
\end{prop}

\begin{prf} If there exists a modular causal geodesic 
  and $(\g,\tau,C)$ is the infinitesimal data of $M$, then
  there exists a $g \in G$ such that the Euler element $h$ satisfies 
  $\Ad(g)^{-1}h = x_\fh + x_\fq$ with $x_\fq \in C^\circ$ and
  $[x_\fh, x_\fq] =0$ (Lemma~\ref{lem:geodesic}).
  As $C$ is elliptic, $x_\fq$ is elliptic.
  Further the pointedness of $C$ implies that
  $x_\fh \in \ker(\ad x_\fq)$ is elliptic. This implies that
  the Euler element  $\Ad(g)^{-1}h$ is elliptic, a contradiction.
\end{prf}

\subsection{The fiber bundle structure of the positivity domain}
\mlabel{subsec:posdom}

The main result of this section is Theorem~\ref{thm:pos-polar}
  in which we exhibit a natural bundle structure on the wedge domain
  $W \subeq M$ that is equivariant with respect to the connected
  group $G^h_e$, the base is the Riemannian symmetric space
  of this group, and the fiber is a bounded convex subset of
  $\fq_\fk$.

\begin{defn} \mlabel{def:3.1}
  Let $h \in \fq_\fp \cap C^\circ$ be a causal Euler element,
  so that $\tau = \tau_h \theta$. 
Then $\fz_\fh(h) = \fh^{\tau_h} = \fh_\fk$ implies that 
\[\cO_h^\fq := \Inn_\g(\fh)h = e^{\ad \fh_\fp}h\] 
is the non-compact Riemannian symmetric space 
associated to the symmetric Lie algebra~$(\fh,\theta)$. 
\end{defn}

\begin{thm} \mlabel{thm:pos-polar}{\rm(Positivity Domain Theorem)} 
  Suppose that $(G,\tau^G,C,H)$ is a connected
  semisimple non-compactly causal Lie group
  for which   $(\g,\tau)$ contains no $\tau$-invariant Riemannian ideals ($\g = \g_s$)
  and that $h$ is a causal Euler element.
  Suppose that $C := C^{\rm max}_\fq(h)$ is the 
  maximal $\Inn_\g(\fh)$-invariant cone with $h \in C^\circ$.
  Then the following assertions hold:
  \begin{itemize}
  \item[\rm(a)] The connected component  $W= W_M^+(h)_{eH}$
    of $eH$ in the positivity domain $W_M^+(h)$ is given by
      \begin{equation}
    \label{eq:weq}
 W = G^h_e.\Exp_{eH}(\Omega_{\fq_\fk}), 
\quad \mbox{ where } \quad 
\Omega_{\fq_\fk} 
= \Big\{ x \in \fq_\fk \: \rho(\ad x) < \frac{\pi}{2}\Big\}.
  \end{equation}
\item[\rm(b)] The polar map 
$\Psi \: G^h_e \times_{G^h_e \cap H} \Omega_{\fq_\fk} \to W, 
[g,x] \mapsto g.\Exp_{eH}(x)$ is a diffeomorphism 
\item[\rm(c)] $W$ is contractible, hence in particular simply connected. 
\item[\rm(d)] $G^h_e \cap H = K^h_e$. 
\end{itemize}
\end{thm}

\begin{prf} (a)   Recall from \cite[Thm.~6.7]{MNO23}
  that the connected component of $h$ in the open subset 
  $\cO_{h} \cap \cT_{C}$ of $\cO_{h}$ is 
\begin{equation}
  \label{eq:crownchar-gen}
  \Inn_\g(\fh)e^{\ad \Omega_{\fq_\fk}} h
  = \Ad(H_e)e^{\ad \Omega_{\fq_\fk}} h \subeq \cT_{C}.
\end{equation}
 If $x \in \Omega_{\fq_\fk}$, then 
$\rho(\ad x) < \pi/2$, so that \eqref{eq:crownchar-gen} implies that 
$g = \exp x$ satisfies  
 \begin{equation}
   \label{eq:proq}
   \Ad(g)^{-1} h = e^{-\ad x} h \in \cT_{C} = \fh + C^\circ.
 \end{equation}
By Lemma~\ref{lem:wm+} 
\begin{equation}
  \label{eq:inc-a}
\Exp_{eH}(\Omega_{\fq_\fk})\subeq W_M^+(h), 
\quad \mbox{ and thus  } \quad 
G^h_e.\Exp_{eH}(\Omega_{\fq_\fk})\subeq W 
\end{equation}
by $G^h$-invariance of $W_M^+(h)$. 

Conversely, for $gH \in W$,
the element $\Ad(g)^{-1}h \in \cO_h \cap \cT_{C}$ is contained 
in the connected component of $h$, so that
\eqref{eq:crownchar-gen} implies that it is contained in
$\Ad(H_e) e^{\ad \Omega_{\fq_\fk}} h$. Therefore
\[ g H_e \exp(\Omega_{\fq_\fk}) \cap G^h \not=\eset.\] 
This is equivalent to 
$ g H_e \cap G^h \exp(\Omega_{\fq_\fk}) \not=\eset,$
 which implies 
\[ gH \in G^h \exp(\Omega_{\fq_\fk}).eH 
=  G^h \Exp_{eH}(\Omega_{\fq_\fk})\] 
and thus 
\begin{equation}
  \label{eq:18}
 W \subeq G^h.\Exp_{eH}(\Omega_{\fq_\fk}).
\end{equation}

If $g \in G^h$ satisfies 
$g \Exp_{eH}(\Omega_{\fq_\fk}) \cap W \not = \eset,$ 
then $g.W = W$ follows from \eqref{eq:inc-a} and 
the fact that $g$ permutes the connected 
components of $W_M^+(h)$. Therefore \eqref{eq:18}, combined with 
\eqref{eq:inc-a}, leads with $G^h_W := \{ g \in G^h \: g.W = W\}$ to 
\[ W 
\subeq G^h_W.\Exp_{eH}(\Omega_{\fq_\fk})
\subeq G^h_W.W = W,\] 
and this entails 
\begin{equation}
  \label{eq:gwform}
 W = G^h_W.\Exp_{eH}(\Omega_{\fq_\fk}).
\end{equation}

Next we observe that 
the exponential map 
$\Exp_{eH} \: \fq_\fk \to M$ 
is regular in every $x\in \Omega_{\fq_\fk}$ because $\rho(\ad x) < \pi/2 < \pi$ 
(\cite[Lemma~C.3(b)]{NO23b}). 
Thus \cite[loc.cit.]{NO23b} further implies that the map 
\[ \Phi \: G^h \times \Omega_{\fq_\fk} \to M, \quad 
(g,x) \mapsto g.\Exp_{eH}(x) \] 
is regular in $(g,x)$ because $\Spec(\ad x) \subeq (-\pi/2,\pi/2)i$ 
does not intersect $\big(\frac{\pi}{2} + \Z \pi\big)i$ for 
$x \in \Omega_{\fq_\fk}$. 
This implies that the differential of $\Phi$ is surjective in each point 
of $G^h \times \Omega_{\fq_\fk}$; hence the image of every 
connected component is open. Now the connectedness of $W$ 
implies that $W \subeq G^h_e.\Exp_{eH}(\Omega_{\fq_\fk}),$ 
and this completes the proof. 

\nin (b)--(d): 
The surjectivity of $\Psi$ follows from 
Theorem~\ref{thm:pos-polar}. 
As $\g^h = \fh_\fk \oplus \fq_\fp$ is a Cartan decomposition of $\g^h$, 
the polar map 
$K^h_e \times \fq_\fp \to G^h_e, (k,x) \mapsto k \exp x$ 
is a diffeomorphism. In particular, 
\[ G^h_e \cap H \subeq G^h_e \cap G^{\tau^G} = K^h_e \subeq H \] 
implies $G^h_e \cap H = K^h_e$ and thus (b). 

The space $G^h_e \times_{G^h_e \cap H} \Omega_{\fq_\fk}$ is a fiber bundle over 
$G^h_e/K^h_e$ whose fiber is the convex set $\Omega_{\fq_\fk}$. 
Therefore it is homotopy equivalent to the base 
$G^h_e/K^h_e$, which is also contractible 
because the exponential map $\Exp_{eH} \: \fq_\fp \to G^h_e/K^h_e$ 
is a diffeomorphism. 

It therefore suffices to show that $\Psi$ is a diffeomorphism. 
The proof of (a) shows already that its differential 
is everywhere surjective, hence invertible by equality of the 
dimensions of both spaces. So it suffices to check injectivity, 
i.e., that $\Exp := \Exp_{eH} \: \fq \to M$ satisfies 
\begin{equation}
  \label{eq:injrel1}
g_1.\Exp(x_1) = g_2.\Exp(x_2)  \quad \Rarrow \quad 
g_2^{-1}g_1 \in K^h_e, \quad 
x_2 = \Ad(g_2^{-1}g_1)x_1.
\end{equation}

\nin {\bf Step 1:} $\Exp\res_{\Omega_{\fq_\fk}}$ 
is injective. 
If $\Exp(x_1) = \Exp(x_2)$, then applying the quadratic representation 
implies $\exp(2 x_1) = \exp(2 x_2)$ in $G$. 
As $x_1$ and $x_2$ are both $\exp$-regular, 
\cite[Lemma~9.2.31]{HN12} implies that 
\[ [x_1, x_2]= 0 \quad \mbox{ and } \quad \exp(2x_1 - 2 x_2)= e.\]
We conclude that $e^{2 \ad(x_1-x_2)} = \id_\g$, and since 
the spectral radius of $2 \ad(x_1 - x_2)$ is less than 
$2 \pi$, it follows that $\ad(x_1 - x_2) = 0$, so that $x_1 = x_2$. 

\nin {\bf Step 2:} $g.\Exp(x_1) = \Exp(x_2)$ with 
$g \in G^h_e$ and $x_1, x_2 \in \Omega_{\fq_\fk}$ implies $g \in K^h_e$. 
Applying the involution $\theta^M$, we see that 
$g.\Exp(x_1)$ is a fixed point, 
so that 
\[ g.\Exp(x_1) = \theta(g).\Exp(x_1) \] 
entails that $\theta(g)^{-1}g$ fixes $m_1 :=\Exp(x_1)$. 
We now write $g = k \exp z$ in terms of the polar decomposition of $G^h_e$ and 
obtain 
\[ \theta(g)^{-1} g = \exp(2z) \in G^{m_1}.\] 
Applying the quadratic representation, we get 
\begin{equation}
  \label{eq:commrel1}
 \exp(2z) \exp(2x_1) \exp(2z) = \exp(2x_1), 
\end{equation}
which can be rewritten as 
\[ \exp(e^{2\ad x_1} 2z) = \exp(-2z).\] 
Since $\ad z$ has real spectrum, so has $e^{2 \ad x_1} z$. 
Therefore the same arguments as in Step 1 above imply that 
\[ [z,e^{2 \ad x_1} z] = 0 \quad \mbox{ and }  \quad 
\exp(2 e^{2 \ad x_1} z + 2z) = e,\] 
and $e^{2 \ad x_1} z = - z$. 
The vanishing $\fh$-component of this element is 
$\sinh(2 \ad x_1) z$, and since 
${\rho(2\ad x_1) < \pi}$, it follows that $[x_1, z] =0$. 
Now \eqref{eq:commrel1} leads to $\exp(4z) = e$, and further to 
$z = 0$, because the exponential function on $\fq_\fp$ is injective. 
This proves that $g = k \in K^h_e$. 

\nin {\bf Step 3:} From \eqref{eq:injrel1} we derive 
\[ g_2^{-1}g_1.\Exp(x_1) = \Exp(x_2),\] 
so that Step 2 shows that 
$k := g_2^{-1}g_1 \in K^h_e$. We thus obtain 
\[ \Exp(x_2) = k.\Exp(x_1) = \Exp(\Ad(k) x_1),\] 
and since $\Ad(k) x_1 \in \Omega_{\fq_\fk}$, we infer from Step 1 that 
$\Ad(k) x_1 = x_2$. This completes the proof. 
\end{prf}

The following corollary identifies 
the connected component of $M^h_C$ containing $eH$  
as a submanifold (cf.\ Lemma~\ref{lem:mx})
of the wedge domain~$W$.

\begin{cor} 
Assume that $\tau_h^G$ exists and leaves $H$ invariant, so that 
$\tau_h^M$ exists and leaves the base point $eH \in M$ invariant.  
Then $\tau^M_h(W) = W$ and
the fixed point set of $\tau_h^M$ in $W$ is the Riemannian 
symmetric space 
\[  W^{\tau_h^M} = M^h_{eH} = G^h_e.eH = \Exp_{eH}(\fq_\fp). \]   
\end{cor} 

\begin{prf}  For $g \in G^h$ and $x\in \fq_\fk$: 
\[ \tau_h^M(g \Exp_{eH}(x)) 
= g \Exp_{eH}(\tau_h(x)) 
= g \Exp_{eH}(\tau(x)) 
= g \Exp_{eH}(-x).\] 
So  $g \Exp_{eH}(x)$ is a fixed point if and only if 
$\Exp_{eH}(-x) = \Exp_{eH}(x),$ 
which is equivalent to $\exp(2x) \in H$.  Now $\tau(x) = -x$ implies 
$\exp(2x) =\exp(-2x)$. As $\rho(2\ad x) < \pi$, 
\cite[Lemma~C.3]{NO23b} further shows that 
$x-(-x) = 2x \in \fz(\g) = \{0\}$. Therefore 
$g \Exp_{eH}(x)$ is a fixed point if and only if $x = 0$. 

From $W = G^h_e.\Exp_{eH}(\Omega_{\fq_\fk})$ 
and the polar decomposition 
$G^h_e = K^h_e \exp(\fq_\fp) = \exp(\fq_\fp) (H_K)_e$ 
(Theorem\ref{thm:pos-polar}(b)) we derive that 
the fixed point set is 
\[ W^{\tau_h^M} = G^h_e.eH = \Exp_{eH}(\fq_\fp) = M^h_{eH}.\qedhere\] 
\end{prf}

The preceding corollary shows that the 
wedge domain $W \subeq M = G/H$ contains the symmetric subspace 
$M^h_{eH} = \Exp_{eH}(\fq_\fp)$ as the fixed point set of an involution. 
Hence the description of $W$ from Theorem~\ref{thm:pos-polar} as 
\[ W = G^h_e.\Exp_{eH}(\Omega_{\fq_\fk}) \] 
suggest to consider $W$ as a real ``crown domain'' of the Riemannian 
symmetric space $M^h_{eH} \cong G^h/H^h$. 

\begin{rem} Theorem~\ref{thm:pos-polar} has a trivial generalization
  to semisimple non-compactly causal Lie algebras
  of the form $\g = \g_k \oplus \g_r \oplus \g_s$
  because then
  \[  C_\fq^{\rm max} = \fq_r \oplus C_{\fq_s}^{\rm max}
    \quad \mbox{ and }  \quad
    \cT_{C_\fq^{\rm max}} = \g_k + \g_r + \cT_{C_{\fq_s}^{\rm max}}.\]
  For $h = h_r + h_s$ with $h_s \in C_{\fq_s}^\circ$ the relation
  $\Ad(g)^{-1}h \in \cT_{C_\fq^{\rm max}}$ is therefore equivalent to 
  $\Ad(g)^{-1}h_s \in \cT_{C_{\fq_s}^{\rm max}}$.
  If $M = \Inn(\g)/\Inn(\g)^\tau \cong M_r \times M_s$
  is the corresponding product decomposition, we obtain
  \[  W_M^+(h) = M_r \times W_{M_s}^+(h_s)
    \quad \mbox{ for } \quad C = C_{\fq}^{\rm max}.\]
  However, if $\g_r \not=\{0\}$, then
  $C_\fq^{\rm max}$ is not pointed, and there are many
  pointed invariant cones $C$, which are not maximal, for which
  the domain $W_M^+(h)$ may have a more complicated structure.
\end{rem}

\begin{ex} \mlabel{ex:gl2} We consider the reductive Lie algebra
    \[ \g = \gl_2(\R) = \R \1 \oplus \fsl_2(\R).\]
    Any Euler element in $\g$ is conjugate to some
    \[ h = \pmat{\lambda & 0 \\ 0 & \mu}\quad \mbox{ with } \quad
    \lambda - \mu = 1.\]
    The Cartan involution $\theta(x) = -x^\top$ on $\g$
    then satisfies $\theta(h) = -h$ and
    $\tau := \theta \tau_h$ acts by 
    \[ \tau\pmat{a & b \\ c & d} = \pmat{-a & c \\ b & -d}.\]
    With the Euler element
    \[ h_0 := \frac{1}{2}\pmat{ 0 & 1 \\ 1 & 0}, \] we then have 
    \[ \fh = \R h_0  = \so_{1,1}(\R) \quad\mbox{
      and } \quad
    \fq = \R \1 + \underbrace{\R h + \R z}_{\fq_s} \quad \mbox{ with }\quad
    z := \frac{1}{2} \pmat{0 & -1 \\ 1 & 0}.\]

    The group $G := \GL_2(\R)_e$ acts by
    $g.A := g A g^\top$ on the $3$-dimensional space $\Sym_2(\R)$
    of symmetric matrices and the stabilizer of
    $I_{1,-1} := \pmat{1 & 0 \\ 0 & -1}$ is the subgroup
    $H := \SO_{1,1}(\R)$ with Lie algebra~$\fh$. Therefore
    $M := G.I_{1,1} \cong G/H$ 
    can be identified with the subspace $\Sym_{1,1}(\R)$ of
    indefinite symmetric matrices. Note that $\R^\times_e\1 = Z(G)_e$ 
    acts by multiplication with $\lambda^2$ and that
    $\R^\times_+ \times M_1 \to M, (\lambda,A) \mapsto \lambda A$
    is a diffeomorphism, where
    \[ M_1 := \{ A \in M \: \det(A) =-1\} \cong \SL_2(\R)/\SO_{1,1}(\R) \cong
    \dS^2\]
    is a realization of $2$-dimensional de Sitter space.
    Note that the determinant defines a quadratic form of signature $(1,2)$
    on $\Sym_2(\R)$ which is invariant under the action of the subgroup 
    \[ \{ g \in \GL_2(\R) \: |\det(g)| = 1\} \supeq \SL_2(\R)\]
    which acts as $\SO_{1,2}(\R)^\up$. 

For the Euler element $h_s := \frac{1}{2}\pmat{1 & 0 \\ 0 & -1}$, we have
    \[ [h_0, h_s] = z, \quad [z,h_s] = h_0 \quad \mbox{ and } \quad [h_0,z] = h_s.\] 
 
According to \cite[Ex.~3.1(c)]{NOO21}, 
all $\Ad(H)$-invariant cones in $\fq$ are Lorentzian of the form 
\[ C_m = \{ x_0 \1 + x_1 (h_s + z) + x_{-1} (h_s - z) \: 
x_1 x_{-1} - m x_0^2 \geq 0, x_{\pm 1} \geq 0\}
\quad \mbox{ for some } \quad m > 0. \]
Actually $C_0 = C_\fq^{\rm max}$ contains $\R \1$ and is not pointed.

\nin (a) We write
\[ h = \frac{\lambda + \mu}{2} \1 +  h_s\]
to see that $h \in C_m$ is equivalent to
\[ m (\lambda+\mu)^2 \leq 1.\]
We also note that the ``semisimple part'' of $C_m$ is 
\[ C_{m,s} = C_m \cap \fq_s = C_m \cap (\R h + \R z)
= \cone(h_s\pm z)\]
coincides with the projection of $C_m$ to $\fq_s$, so that
$C_{m,s}^\circ = C_m^\circ \cap \fq_s$.

Write $W(C_m,h)$ for the positivity domain of the Euler element
$h$ with respect to the causal structure specified by the cone~$C_m$.
Then Theorem~\ref{thm:pos-polar} implies that
\[ W(C_0, h_s) = G^{h_s}_e.\Exp_{eH}(\Omega_{\fq_\fk}), \quad \mbox{ where } \quad
\Omega_{\fq_\fk} = \Big(-\frac{\pi}{2},\frac{\pi}{2}\Big) z.\]
For $x \in \Omega_{\fq_\fk}$ we have
\[ e^{-\ad x}.h_s \in C_s + \fh \subeq C_m + \fh\]
(see \eqref{eq:proq}) 
and $G^h = G^{h_s}$, 
so that we have 
\[ W(C_0,h_s) = G^h_e.\Exp_{eH}(\Omega_{\fq_\fk}) \subeq W(C_m,h_s)
\subeq W(C_0, h_s)\]
implies the equality
\[ W(C_m,h_s) = G^h_e.\Exp_{eH}(\Omega_{\fq_\fk}) \quad \mbox{ for all } \quad
m > 0.\]

We also note that
\[ W(C_m,h) \subeq W(C_0,h) = W(C_0, h_s) = G^h_e.\Exp_{eH}(\Omega_{\fq_\fk}) \]
because $\Ad(g)^{-1}h = h_z + \Ad(g)^{-1}h_s \in C_0^\circ$ if and only if
$\Ad(g)^{-1}h_s \in C_s^\circ$.

To determine the domain $W(C_m,h)$ in general, we write
\[ h = h_z + h_s = \frac{\lambda+\mu}{2}\1 + \frac{1}{2}\pmat{1 & 0 \\
    0 & -1}.\]
By $G^h_e$-invariance, we have to determine when
$\Exp_{eH}(t z)$, $|t| < \frac{\pi}{2}$, is contained in $W(C_m,h)$.
For $g = \exp(tz)$ we have
\[ \Ad(g)^{-1} h
  = h_z + e^{-t \ad z} h_s = h_z + \cos(t) h_s - \sin(t) h_0, \]
and
\[ p_\fq(\Ad(g)^{-1}h) = h_z + \cos(t) h_s
  = \frac{\lambda+\mu}{2} \1 + \frac{\cos(t)}{2}(h_s + z) +
  \frac{\cos(t)}{2}(h_s - z).\] 
We then have
\[ x_0 = \frac{\lambda + \mu}{2} \quad\mbox{ and } \quad
  x_{\pm 1} = \frac{\cos(t)}{2}.\] 
We conclude that, for $|t| < \frac{\pi}{2}$,
the inclusion $h_z + \cos(t) h_s \in (C_m)^\circ$ is equivalent to
\[ x_1 x_{-1} - m x_0^2 = \frac{1}{4}(\cos^2(t) - m(\lambda + \mu)^2) >  0.\]
We thus obtain the condition
\[ |t| < \arccos(\sqrt{m} |\lambda + \mu|).\] 
For $m> 0$ and $h \not= h_s$, this is specifies a proper subinterval
of $(-\frac{\pi}{2}, \frac{\pi}{2})$.

\nin (b) To determine which cone $C_m$ corresponds to the canonical order 
  on the space $\Sym_{1,1}(\R)$, induced from the natural
  order of $\Sym_2(\R)$ (which is also Lorentzian), we evaluate
  the tangent map $\fq \to  \Sym_2(\R), x \mapsto x I_{11} + I_{11} x^\top$ to
  \[ \1.I_{1,1} = 2 I_{1,1}, \quad
    h_s.I_{1,1} = \1, \quad 
    z.I_{1,1} = \pmat{0 & 1 \\ 1 & 0}.\]
  We thus obtain for
  $x = x_0 \1 + x_1 (h_s + z) + x_{-1} (h_s - z)$ that
  \[ x.I_{1,1} = x I_{11} + I_{11} x^\top
    = \pmat{
      2 x_0 + x_1 + x_{-1} & x_1 - x_{-1} \\ 
      x_1 - x_{-1} & -2 x_0 + x_1 + x_{-1}}.\]
By the Hurwitz criterion, this matrix is positive semidefinite if and only if
  \[  x_1 + x_{-1} \geq |2x_0| \]
  and
\[ (x_1  + x_{-1})^2 - 4 x_0^2 - (x_1 - x_{-1})^2
  = 4( x_1 x_{-1}  - x_0^2) \geq 0.\]
Is $x_1 + x_{-1} \geq 0$, then these two inequalities are equivalent to
$x_1 x_{-1} - x_0^2 \geq 0$. As these two conditions imply
that $x_{\pm 1} \geq 0$, we see that the canonical order on $M$ corresponds to
the cone $C^1$, i.e., to $m = 1$.

\nin (c)   For the modular vector field $X_h$ we have
  \[ X_h\pmat{a & b \\ b & d}
  = \pmat{
    2\lambda a & (\lambda + \mu) b \\ 
    (\lambda + \mu)b & 2\mu d}.\]
The positivity domain of $X_h$ depends on $\lambda$,
and with this formula one can also determine the positivity domain
quite directly for $m = 1$, where $C^1$ corresponds to the canonical order.
\end{ex}

\begin{ex} (cf.\ \cite[Exs.~2.11, 2.25]{NO23b}) 
  Let $G := \GL_n(\R)_+$ and $K := \SO_n(\R)$.
  We consider the Riemannian symmetric space
  \[ M_r := \Sym_n(\R)_+ \cong \GL_n(\R)/\SO_n(\R) \]
  and the corresponding irreducible subspace \
  \[ M_{r,s} := \{ A \in M \: \det(A) = 1\}\cong \SL_n(\R)/\SO_n(\R)\]
  (here the index $s$ refers to ``semisimple''). 
  On $\g = \gl_n(\R)$, we consider the
Cartan involution given by $\theta(x) = - x^\top$  
and write $n = p + q$ with $p,q> 0$. Then 
\begin{equation}
  \label{eq:euler-pq}
  h^p_s := \frac{1}{n}\pmat{ q \1_p & 0 \\ 0 & -p\1_q}  \in \fsl_n(\R)
  \quad \mbox{ and }\quad
    h^p := h^p_s - \frac{q}{n}\1 = \pmat{ 0 & 0 \\ 0 & \1_q} \in \g  
\end{equation}
are Euler elements and 
$\tau := \tau_{h^p}\theta$ leads to a non-compactly causal symmetric 
Lie algebra $(\g,\tau,C)$, where 
\[ \fh = \so_{p,q}(\R) \quad \mbox{ and } \quad 
\fq = \Big\{\pmat{a & b \\ -b^\top & d} \: a^\top = a, d^\top = d \Big\}\] 
To identify $G/H$ in the boundary of the crown domain
in $G_\C/K_\C \cong G_\C.\1 \cong \Sym_n(\C)^\times$,
where $G_\C$ acts on $\Sym_n(\C)$ by $g.A := gAg^\top$
(\cite[Thm.~5.4]{NO23b}), 
we observe that
\[ \exp(it h^p).\1 = e^{2ith^p}
  = \cos(2t h^p) + i \sin(2t h^p)
  = \pmat{ \1_p & 0 \\ 0 &  (\cos(2t) + i \sin(2t)) \1_q},\]
so that we obtain for $t  =\frac{\pi}{2}$ the matrix
\[ \exp\big(\frac{\pi i}{2} h^p\big).\1 =  I_{p,q}.\]
The $G$-orbit of this matrix is the open subset 
\[ M := G.I_{p,q} = \{ gI_{p,q} g^\top \: g \in \GL_n(\R)_+\}
  = \Sym_{p,q}(\R) \]
of symmetric matrices of signature $(p,q)$.
We have
\[ X_{h^p}\pmat{ a & b \\ b^\top & d}
  = h^p \pmat{ a & b \\ b^\top & d} + \pmat{ a & b \\ b^\top & d}h^p
  = \pmat{ 0 & b \\ b^\top & 2 d}.\] 
These matrices are {\bf never} positive definite.
So we have to take $h_s$ instead to find non-trivial positivity
domains.

For the case $p = q = 1$ and $n = 2$ this has been carried out
in Example~\ref{ex:gl2}. We also write
\[ h = \pmat{ \lambda \1_p & 0 \\ 0 & \mu \1_q}
  \quad \mbox{ with } \quad \lambda - \mu = 1.\]
Then 
\[ X_h\pmat{ a & b \\ b^\top & d}
  = \pmat{ 2\lambda a & (\lambda + \mu) b \\
    (\lambda + \mu) b^\top & 2\mu d},\]
so that
\[ X_h(I_{p,q}) = 2 \pmat{ \lambda \1_p & 0 \\  0  & -\mu \1_q} \> 0
\quad \mbox{ if }  \quad 0 < \lambda < 1, \]
which is equivalent to $\lambda\mu < 0$.
\end{ex}

\subsection{The connected components of $M_C^h$} 

The main result in this section is Proposition \ref{prop:39} below
on the subgroup $H_K$ of $K^h$.
We then discuss several examples to clarify the situation.

\begin{prop}\mlabel{prop:39} {\rm(Connected components of $M^h_C$)}
  If $G = \Inn(\g)$ and $(\g,\tau)$ is irreducible ncc with causal
  Euler element $h$, then
  $\pi_0(M^h_C) \cong K^h/H_K$ contains at most two elements.
\end{prop}

\begin{prf}   We recall from Proposition~\ref{prop:3.1}(c)
  that   $M^h_C = G^h.eH$.
 With \cite[Thm.~IV.3.5]{Lo69} we see that 
the symmetric space $G^h.eH \cong G^h/H^h$ is a vector bundle
over $K^h/H_K^h$, hence in particular 
homotopy equivalent to $K^h.eH \cong K^h/H_K^h$. 
In view of Proposition~\ref{prop:hhcompact}(c),
we have for $G = \Inn(\g)$
that $H^h= H^{\tau_h} = H_K \subeq G^\tau$ is a maximal compact subgroup of~$H$.
It follows in particular that
$H^h = H_K^h \subeq K^h$.
We conclude that
$\pi_0(M^h_C) \cong \pi_0(K^h/H_K)$. From \cite[\S 7]{MNO23}, we know that
$\pi_0(G^h)\cong \pi_0(K^h)$ has at most two elements. 
\end{prf}

\begin{ex} (The inclusion $H_K \subeq K^h$ may be proper)
We have $G^h = K^h \exp(\fq_\fp)$ and 
$K^{\tau^G} = K^{\tau_h^G}$ because $K = G^\theta$.
Further $H_K \subeq K^h$ by Proposition~\ref{prop:hhcompact}(a),
  so that the equality $H_K = K^h$ is equivalent to 
  $K^h \subeq H_K$. This may fail for two reasons.
  One is failure in the adjoint group $\Inn(\g)$
  (Proposition~\ref{prop:39}), and the other reason is that
  $Z(G)$ may be non-trivial.

Assume that $\fg$ is semisimple and
$(\fg,\tau, C)$ ncc. Let $G$ be a corresponding connected Lie group
on which $\tau^G$ exists (for $\tau = \tau_h \theta$)
and $H := G^{\tau^G}_e$. For the connected group $K := G^\theta$,
the intersection $H_K := H \cap K =  \la \exp\fh_\fk \ra$ 
is connected but $K^h \supeq Z(G) H_K$ is in general not connected
because $Z(G)$ need not be contained in $H_K$.

This can be seen easily for $\g = \fsl_2(\R)$.
For 
\begin{equation}\label{eq:hsl2}
  h = \frac{1}{2}\begin{pmatrix} 1 & 0 \\ 0 & -1\end{pmatrix}, \quad
  \theta(x) = - x^\top, 
  \quad \mbox{ we have } \quad
  \tau = \theta\tau_h,  
  \quad \begin{pmatrix} a & b\\ c & -a\end{pmatrix}
    \mapsto \begin{pmatrix} -a & c\\ b & a\end{pmatrix}.
\end{equation}
For any connected Lie group $G$ with Lie algebra $\g$, 
the group $K = G^\theta$ is connected $1$-dimensional and $\tau(k) = k^{-1}$
for $k \in K$. Moreover, $K^h = Z(G)$ is a discrete subgroup
which intersects $H = \exp \fh \cong \R$ trivially.
Even the inclusion $K^h \subeq G^{\tau^G}$ fails if $|Z(G)| \geq 3$,
i.e., if $\tau$ acts non-trivially on $Z(G)$.
Note that $Z(G)$ is infinite if $G$ is simply connected.
\end{ex}

\begin{ex} (a)  For  $\g = \fsl_2(\R)$ we consider again the Euler element
$h$ from \eqref{eq:hsl2} 
and the Cartan involution $\theta(x) = -x^\top$. 
By Lemma~\ref{lem:geodesic}, the $\alpha$-orbit of $gH$ 
is a geodesic if and only if $\Ad(g)^{-1}h$ commutes with 
$\tau(\Ad(g)^{-1}h) = - \Ad(\tau(g))^{-1}h$, i.e., if 
\[ \Ad(g\tau(g)^{-1})h \in \fz_\g(h) = \R h.\] 
As $\cO_h \cap \R h = \{  \pm h\}$, this leaves two possibilities: 
\begin{itemize}
\item[(1)] If $\Ad(g\tau(g)^{-1})h= h$, then 
$\Ad(\tau(g))^{-1}h = \Ad(g)^{-1}h$ implies $\Ad(g)^{-1}h \in \fq$. 
\item[(2)] If $\Ad(g\tau(g)^{-1})h= -h$, then 
$-\Ad(\tau(g))^{-1}h = \Ad(g)^{-1}h$ implies  $\Ad(g)^{-1}h \in \fh$. 
In this case $gH$ is a fixed point of the modular flow. 
\end{itemize}

\nin (b) For $\g = \fsl_{2k}(\R)$ with the Cartan involution $\theta (x) = - x^\top$ and  the causal Euler element 
\[ h = \frac{1}{2} \pmat{  \1_k & 0 \\ 0 & -\1_k},\] 
we obtain $\fh = \so_{k,k}(\R)$ for $\tau = \theta \tau_h$.  There exists a subalgebra 
$\fs \cong \fsl_2(\R)^k$, where the 
$\fsl_2$-factors correspond to the coordinates 
$x_j$ and $x_{j + k}$ for $1 \leq j \leq k$. 
Accordingly, $h = \sum_{j = 1}^k h_j$, where the Euler elements 
$h_j$ in the $\fsl_2$-factors are conjugate to Euler elements $h_j'$ in $\fh$. 
Therefore the ``geodesic condition'' is satisfied by all elements 
$\sum_{j = 1}^k \tilde h_j \in \cO_h$, where $\tilde h_j$ is either $h_j$ or 
$h_j'$. 
\end{ex}

The following example shows that modular geodesics also exist 
in symmetric spaces without causal structure. They can be 
``space-like'' 
rather than ``time-like'', resp., causal. 
\begin{ex} \mlabel{ex:hypspace}
The $d$-dimensional hyperbolic space 
\[ \bH^d := \{ x = (x_0, \bx) \in \R^{1,d} \: x_0^2 - \bx^2 = 1, x_0 > 0 \} 
\cong \SO_{1,d}^\up(\R)/\SO_d(\R) \] 
carries a modular flow specified by 
any Euler element $h \in \fq \subeq \so_{1,d}(\R)$ (corresponding to a tangent 
vector of length $1$). Every geodesic of $\bH^d$ is an 
orbit of the flow generated by an Euler element of $\so_{1,d}(\R)$. 
\end{ex}

\begin{rem} \mlabel{rem:4.13} 
Let $(\g,\tau,C)$ be a simple ncc symmetric 
Lie algebra.  In general, we have for a causal Euler element  
$h \in \cE(\g) \cap C^\circ$ a proper inclusion 
\[ \Inn_\g(\fh)h  =  \cO_h \cap C^\circ\subeq \cO_h \cap \fq.\] 
By Lemma~\ref{lem:3.11}, this implies that $M^h$ 
is not connected and $M^h_C \not= M^h$. 

For instance if $\g= \fh_\C$ and $\fh$ is simple hermitian of tube type, 
then we obtain for any pointed generating invariant cone 
$C_\fh \subeq \fh$ a hyperbolic cone $C := i C_\fh \subeq \fq = i \fh$. 
If $h \in \cE(\g) \cap C^\circ$ is a causal Euler element, then 
$-h \in \Ad(G)h$ follows from \cite[Thm.~3.10]{MN21} and the subsequent discussion, 
but $-h \not\in C^\circ$;
see also \cite[Thm.~5.4]{MNO23}. 
\end{rem}

\begin{ex} (a) 
  For de Sitter space $M = \dS^d$ (cf.\ Example~\ref{ex:deSitter1}
  and Appendix~\ref{subsec:dS1}),
  the subspace 
$M^h_{eH} = \Exp_{eH}(\R h)$ is a single geodesic, 
hence in particular $1$-dimensional. Note that $\dim \fq_\fp= 1$ in this case. 
The modular flow on $M$ has the fixed point set 
$M^\alpha \cong \bS^{d-2}$. 

\nin (b) For $M = G_\C/G$, $\g$ hermitian, we have $M^h_{eG} = \Exp_e(i\fk)$ 
with dual symmetric space the group~$K$, considered as a symmetric space.
\end{ex}

\section{Open $H$-orbits in flag manifolds and a convexity theorem} 
\mlabel{sec:convtheo1} 

In this section we prove a convexity theorem that is 
vital to derive the equality $W = W(\gamma)$ in the next section. Here,
as above, $W= W_M^+(h)_{eH}$.

Let 
$P^- := \exp(\g_{-1}(h)) G^h \subeq G$  be the ``negative'' parabolic 
subgroup of $G$ specified by $h$ and identity 
$\g_1(h)$ with the open subset $\cB := \exp(\g_1(h)).eP^- \subeq G/P^-$. 
Then $\cD := H.0 \subeq \cB$ is an open convex subset, and our
convexity theorem  (Theorem~\ref{thm:contheo}) asserts that, for any 
$g \in G$ with $g.\cD \subeq \cB$, the subset $g.\cD \subeq \cB$ 
is convex. 

We consider a {\bf connected semisimple} Lie group $G$ with Lie algebra 
$\g$ and an Euler element $h \in \g$. 
We put 
\[ \fn^\pm := \g_{\pm 1}(h) \quad \mbox{ and } \quad 
N^\pm := \exp(\fn^\pm), \] 
and write 
\[ P^\pm := \{ g \in G \: \Ad(g) \g_{\pm 1}(h) = \g_{\pm 1}(h)\} 
= N^\pm G^h \cong N^\pm  \rtimes G^h\] 
(see \cite[Thm.~1.12]{BN04} for the equality) 
for the corresponding  maximal parabolic subgroups. We write 
\[ \cM_\pm := G/P^\mp \] 
for the corresponding flag manifold. 
The abelian subgroup $N^+$ has an open 
orbit $\cB := N^+.eP^- \subeq \cM_+$, which we call the {\it open Bruhat cell}. 
It carries a natural affine structure because the map 
\[ \phi \: \fn^+ \to \cB := N^+.eP^-, \quad x \mapsto \exp(x) P^- \] 
defines an open embedding. Below we shall always use these coordinates on~$\cB$. 

Choose a Cartan involution $\theta$ with $\theta(h) = -h$ 
and consider the involution $\tau := \theta e^{\pi i \ad h}$. We write 
\[ H := K^h \exp(\fh_\fp) \quad \mbox{ for } \quad  
\fh_\fp = \g^\tau \cap \fp\quad \mbox{ and } \quad 
H_K := K^h.\] 
Then 
\[ P^{\pm}\cap H = G^h \cap H = K^h, \] 
so that  
\begin{equation}\label{eq:boundedReal} \cD_+ := H.eP^- \cong H/H_K 
\end{equation}
is an open $H$-orbit in $\cB \subeq G/P^-$.  It is a real
bounded symmetric domain (\cite[Thm.~5.1.8]{HO97})
and  coincides with the unit 
ball in the positive real Jordan triple 
\begin{equation}
  \label{eq:jts}
V := \fn \quad \mbox{ with }\quad 
\{x,y,z\} := x \square y.z = - \frac{1}{2}[[x,\theta(y)],z]
\end{equation}
(cf.\ \cite[(4.6)]{BN04}) 

\subsection{The open $H$-orbits in $G/P^\pm$}

\begin{lem} {\rm(\cite[Cor.~1.10]{BN04})} 
For $y \in \g_{-1}(h)$ and $x \in \g_1(h)$, 
we have $\exp(y).\exp(x) P^- \in \cB$ if and only if the 
 {\it Bergman operators}
\[ B_+(x,y) := \1 + \ad(x)\ad(y) + \frac{1}{4} (\ad x)^2 (\ad y)^2  
\in \End(\g_1(h)) \]
and 
\[ B_-(y,x) := \1 + \ad(y)\ad(x) + \frac{1}{4} (\ad y)^2 (\ad x)^2 
\in \End(\g_{-1}(h))\] 
are both invertible. 
\end{lem}

\begin{rem} Note that 
\[ \theta B_-(y,x) \theta 
= \1 + \ad(\theta(y))\ad(\theta(x)) + \frac{1}{4} (\ad \theta(y))^2 (\ad \theta(x))^2 
= B_+(\theta(y),\theta(x)).\] 
\end{rem}

\begin{ex} \mlabel{ex:sl2}
We consider the group $G = \SL_2(\R)$ with Lie algebra 
$\g = \fsl_2(\R)$ and the linear basis 
\begin{equation}
  \label{eq:sl2raa}
h := \frac{1}{2} \pmat{1 & 0 \\ 0 & -1},\quad 
e =  \pmat{0 & 1 \\ 0 & 0}, \quad
f =  \pmat{0 & 0 \\ 1 & 0}, 
\end{equation}
satisfying 
\[ [h,e] = e, \quad [h,f] = -f, \quad [e,f] = 2h.\] 

Then 
\[ N^+ = \pmat{1 & \R \\ 0 & 1}, \quad 
N^- = \pmat{1 & 0 \\ \R & 1}, \quad \mbox{ and } \quad 
G^h = \Big\{ \pmat{a & 0 \\ 0 & a^{-1}} \: a \in \R^\times \Big\},\] 
so that 
\[ P^- 
= \Big\{ \pmat{a & 0 \\ c & a^{-1}} \: a \in \R^\times, c \in \R \Big\}
= \{ g \in G \: g e_2 \in \R e_2\}.\] 
For $K = \SO_2(\R)$, we have $K^h = \{ \pm \1\}$. 
Identifying $G/P^-$ with the projective space $\bP(\R^2) = G.[e_2]$, 
the Bruhat cell is 
\[ \cB = \{ [x:1] \: x \in \R \} \cong \R,\] 
and $G$ acts by 
\[ g.x = \frac{ax + b}{cx+d}\quad \mbox{ for } \quad ax + b \not=0.\] 
In particular, we have  
\begin{equation}
  \label{eq:yact}
\exp(yf).x = \frac{x}{1 + xy}.
\end{equation}

We consider the Cartan involution $\theta(x) = -x^\top$, 
so that $\tau := \theta e^{\pi i \ad h}$ acts by 
\[ \tau\pmat{a & b \\ c & -a} 
  = \pmat{-a & c \\ b & a} \quad \mbox{ and } \quad
  \tau^G\pmat{a & b \\ c & d} 
= \pmat{d & c \\ b & a}. \] 
Then 
\[ H = \SL_2(\R)^{\tau^G} = 
\Big\{ \pmat{a & b \\ b  & a} \: a^2 - b^2 = 1\Big\}
= \SO_{1,1}(\R),\] 
so that 
\begin{equation}
  \label{eq:cD-for-sl2}
\cD_+ = H.0 = \{ ba^{-1} \: a^2 - b^2 = 1\} = (-1,1).
\end{equation}
Note that $\Ad(H) \cong H/\{ \pm \1\}$ is connected. 

The Jordan triple product satisfies 
\[ \{e,e,e\} = - \frac{1}{2}[[e,\theta(e)],e] = -\frac{1}{2}[[e,-f],e] 
= \frac{1}{2}[2h,e] = e,\] so that 
\[ \{xe,ye,ze\} =  xyz \cdot e\quad \mbox{ for } \quad x,y,z \in \R.\]
Further 
\[ (\ad f)e = -2h, \quad (\ad f)^2e = -2[f,h] = -2f, \quad 
 (\ad e)f = 2h, \quad (\ad e)^2f  = -2e\] 
implies 
\begin{align*}
 B_+(x e, y f) e
&= e + xy \ad(e) \ad(f)e + \frac{x^2y^2}{4} (\ad e)^2 (\ad f)^2 e\\
&= e + xy \ad(e)(-h) + \frac{x^2y^2}{4} (\ad e)^2 (-2f) \\
&= e + xy 2 e  + \frac{x^2y^2}{4} 4 e 
= (1 + 2 xy + x^2 y^2) e 
= (1 + xy)^2 e.
\end{align*}
Moreover, 
\[  B_+(\theta(yf), \theta(xe))e
= B_+(-ye, -xf)e  = (1 + (-y)(-x))^2 e.
= B_+(xe,yf)e.\] 
As $1 +xy$ is invertible for all $x$ with $|x| < 1$ 
if and only if $|y| \leq 1$, it follows that 
\[ \exp(y).\cD_+ \subeq \cB \quad \Leftrightarrow \quad
|y| \leq 1.\] 
\end{ex}

 Now back to the general case. In the following we write $\|\cdot\|$ for the spectral norm on the 
Jordan triple system $\g_1(h) = \fn^+$. If 
$x = \sum_{j = 1}^k x_j c_j$ with pairwise orthogonal tripotents $c_j$, then 
\begin{equation}
  \label{eq:normg1}
  \|x\| = \max \{ |x_1|, \ldots, |x_k|\}.
\end{equation}
If
\begin{equation}\label{def:Dg}
 \cD_\g := \{ x \in \g_1(h) \: \|x\| < 1 \},
 \end{equation}
then we have
\begin{equation}
  \label{eq:cDg}
  \cD_+ = \exp(\cD_\g) P^- \subeq G/P^-
\end{equation}
(\cite[Thm.~5.1.8]{HO97}). 

\begin{prop} \mlabel{prop:1.1}
The following assertions hold: 
  \begin{itemize}
  \item[\rm(a)] $g.\cD_+ \subeq \cB$ is equivalent to 
$g \in P^+ \exp(y)$ for $y \in \fn^-$ with $\|y\| \leq 1$. 
  \item[\rm(b)] $g.\cD_+ \subeq \cB$ is relatively compact if and only if 
$g \in P^+ \exp(y)$ for $y \in \fn^-$ with $\|y\| <1$. 
  \end{itemize}
\end{prop}

\begin{prf} The condition $g.eP^- \in \cB$ is equivalent to 
$g \in N^+ P^- = N^+ G^h N^- = P^+ N^-.$ 
Let $y \in \fn^-$ with $g \in P^+ \exp(y)$. Then 
the invariance of $\cB$ under $P^+$ 
implies that $g.\cD_+ \subeq \cB$ is equivalent to $\exp(y).\cD_+ \subeq \cB$. 

\nin (a) Suppose first that $\exp(y).\cD_+ \subeq \cB$. 
By the Spectral Theorem for positive Jordan triples 
(\cite[Thm.~VI.2.3]{Ro00}\begin{footnote}{This theorem 
is stated for complex hermitian Jordan triple 
systems, but $V = \g_1(h)$ is a real form of the complex 
JTS $V_\C = \g_{\C,1}(h)$ on which we have an antilinear isomorphism 
$\sigma$ with $V = V_\C^\sigma$. Therefore the uniqueness in 
the spectral decomposition shows that, for $x \in V$, the corresponding 
spectral tripotents are contained in $V$.}\end{footnote}), 
there exist 
pairwise orthogonal tripotents 
$c_1, \ldots, c_k$ and $\beta_1, \ldots, \beta_k \in \R$ with 
\[ y = \sum_{j = 1}^k \beta_j \theta(c_j)\] 
(\cite[Thm.~VI.2.3]{Ro00}). For 
$x = \sum_j \alpha_j c_j$ and 
$z = \sum_j \gamma_j c_j$, we then have 
\[ \{x,\theta(y),z\} = \sum_{j = 1}^k \alpha_j \beta_j \gamma_j  \cdot c_j \] 
(\cite[Prop.~V.3.1]{Ro00}). 
As $x \in \cD_\g$ is equivalent to  
\[ \|x\|= \max \{ |\alpha_j| \: j =1,\ldots, k\} < 1,\] 
the calculations in Example~\ref{ex:sl2} show that 
$\exp(y).\cD_+ \subeq \cB$ implies 
$\|y\| \leq 1$.\begin{footnote}
{As $\cD_+$ is invariant under the group $(H_K)_e$ which acts linearly,  
and this group acts transitively on the set of all maximal 
flat subtriples of $V$ 
(\cite[Lemma~VI.3.1]{Ro00}), it suffices to shows that an element 
with a spectral resolution $x = \sum_{j = 1}^r x_j c_j$  is contained 
in $\cD_\g$ if and only if $|x_j| < 1$ for every~$j$. 
This follows easily from \eqref{eq:cD-for-sl2}.}\end{footnote}

To prove the converse, suppose first that $\|y\| < 1$. 
Then 
\[ \exp(-y) = \theta(\exp(-\theta(y)) \in \theta(H P^-) 
= H P^+ \] 
implies $\exp(y) \in P^+ H,$  so that  
\[ \exp(y).\cD_+ \in P^+H.\cD_+ = P^+.\cD_+ \subeq \cB.\] 

Now we assume that $\|y\| = 1$. 
We observe that 
\[ \exp(y).x = \exp(-th) \exp(th) \exp(y).x 
= \exp(-th) \exp(e^{-t} y).(e^t x),\] 
so that, for $r > 0$, $\exp(y).x \in \cB$ 
is equivalent to $\exp(r^{-1}y).(rx) \in  \cB$. 
For $x \in \cD_\g$, we pick 
$r > 1$ with $rx \in \cD_\g$. Then $\|r^{-1}y\| < 1$ implies 
$\exp(r^{-1}.y). \exp(rx) \in \cB,$ 
and thus $\exp(y).\exp(x) \in \cB$. This shows that 
$\exp(y).\cD_+ \subeq \cB$. 
  
\nin (b) If $\|y\| < 1$, then the argument under (a) shows that 
$\exp(y).\cD_+ \subeq P^+.\cD_+$ is relatively compact. 

Now we assume that $\|y\| = 1$. 
We show that this implies that $\exp(y).\cD_+$ is unbounded. 
As above, we use the Spectral Theorem to write 
\[ y = \sum_{j = 1}^k \beta_j \theta(c_j)\] 
and observe that there exists an $\ell \in \{1,\ldots, k\}$ with $|b_{\ell}| = 1$. 
For $x = \sum_j \alpha_j c_j \in \cD_\g$ we then obtain with 
\eqref{eq:yact}
\[ \exp(y).x = \sum_{j = 1}^k \frac{\alpha_j}{1 + \alpha_j \beta_j} c_j.\] 
For $x = \alpha c_{\ell}$ we get in particular 
\[ \exp(y).x 
= \frac{\alpha}{1 + \alpha \beta_{\ell}} c_{\ell}.\] 
For $\alpha \to -\sgn(\beta_\ell)$ these element leave every compact 
subset of $\cB$. Therefore $\exp(y).\cD_+$ is unbounded. 
\end{prf}

\begin{thm}{\rm(Convexity Theorem for conformal balls)} \mlabel{thm:contheo}
If $g \in G$ is such that 
$g.\cD_+ \subeq \cB$, then $g\cD_+$ is convex. 
If $g.\cD_+$ is relatively compact in $\cB$, then 
there exists an element $p \in P^+$ with $g.\cD_+ = p.\cD_+$, so that 
$g.\cD_+$ is an affine image of $\cD_+$.
\end{thm}

\begin{prf} If $g.\cD_+ \subeq \cB$ is relatively compact, then 
Proposition~\ref{prop:1.1}(b) and its proof imply 
the existence of $p \in P^+$ with $g.\cD_+ = p.\cD_+$. 
In particular $g.\cD_+$ is an affine image of $\cD_+$ and 
therefore convex. 

If $g.\cD_+\subeq \cB$ is not relatively compact, then we put 
$r_n := 1-\frac{1}{n}$. Now  
\[ \exp y.\cD_+ = \bigcup_{n \in \N} \exp y\exp(r_n\cD_\g)P^- \] 
is an increasing union. Therefore it suffices to show that the 
subsets $\exp y\exp(r_n\cD_\g)P^-$ are convex. 
For $r_n = e^t$ we have 
\begin{align*}
 e^{-t}\cdot (\exp y\exp(r_n \cD_\g)P^-) 
&= \exp(-th).(\exp(y)\exp(r_n \cD_\g)P^-) \\
&= \exp(e^t y).\exp(e^{-t}r_n \cD_\g)P^- 
  = \exp(r_n y).\cD_+,
\end{align*}
and these sets are convex by the preceding argument. 
\end{prf}

\begin{ex} \mlabel{ex:deSitter1}
  We consider $G = \SO_{1,d}(\R)_e$ 
as the identity component of the conformal group of the euclidean space 
$\R^{d-1}$, $H = G_{\be_1} = \SO_{1,d-1}(\R)$, and the 
Euler element $h \in \so_{1,1}(\R)$ with $h.\be_0 = \be_1$ and $h.\be_1 = \be_0$. 
As $Z(\SO_{1,d-1}(\R)) \subeq \{ \pm \1\}$ and $G$ preserves the positive 
light cone, the center of $G$ is trivial. 

The symmetric space $M = G/H \cong G.\be_1 \cong \dS^d$ 
is $d$-dimensional de Sitter space, 
$P = N^- \rtimes G^h$ is the stabilizer of the positive light ray 
$\R_+(\be_0 - \be_1)$, and $G/P \cong \bS^{d-1}$ is the sphere 
of positive light rays. 
On the sphere $\bS^{d-1}$, the subgroup $H$ has two open orbits 
which are positive half-spheres separated by the sphere $\bS^{d-2}$ of 
positive light rays in the subspace $\be_1^\bot$. 

In the sphere the Bruhat cells are the point complements and 
if $g.\cD_+ \subeq \cB \cong \R^{d-1}$, then the convexity of $g.\cD_+$ is 
well-known from conformal geometry because conformal images of balls 
are balls or half spaces. 
\end{ex}

\subsection{The subset realization of the ordered
space $M = G/H$} 

As before $G$ is assumed to be a connected semisimple Lie group. 
To simplify the notation we write $\cM$ for  $\cM_+ = G/P^-$.
Recall the following fact about the compression semigroup
of the $H$-orbit $\cD_+ = H.eP^- \subeq \cM_+$, 
which is the Riemannian symmetric space $H/H\cap K$.

\begin{lem} \mlabel{lem:Compress} The compression semigroup of
the open $H$-orbit $\cD_+ = H.eP^-\subeq G/P^-$ is
\begin{equation}
 \label{eq:compcD}
  \comp({\cD_+}) = \{ g \in G \: g.{\cD_+} \subeq {\cD_+} \}
  = H \exp(-C_\fq^{\rm max}(h)),
\end{equation}
\end{lem}

\begin{prf} This result was announced in \cite{Ol81, Ol82} and 
a detailed proof was given in   \cite[Thm.~VI.11]{HN95} 
for the case where $G \subeq G_\C$, $G_\C$ is simply connected
  and $H = G^\tau$. In this case $\Ad(G^\tau)$ preserves
  $C_\fq^{\rm max}(h)$, so that $G^\tau \subeq K^h \exp(\fh_\fp)$.
  Conversely, $K^h$ leaves $\cD_+$ invariant, so that we obtain
  $G^\tau = H = K^h \exp(\fh_\fp)$ in this particular case.

  To see that the lemma also holds in the general case,
  note that the center of $G$ acts trivially on $G/P^-$ and that 
  $Z(G) \subeq K^h\subeq H$.
  Therefore the general assertion follows if 
  the equality \eqref{eq:compcD} holds at least for one
  connected Lie group $G$ with Lie algebra $\g$. 
Hence it follows from the special case discussed above. 
\end{prf}

 We now use this to realize $G/H$ as an ordered symmetric space
as a set of subsets of $\cM$ and describe the ordering in that realization.

 \begin{prop} \mlabel{prop:setreal} 
{\rm(The subset realization of  ncc symmetric spaces)} Let $G$ be 
a connected semisimple Lie group, $h \in \g$ an
Euler element, $\theta$ a Cartan involution with 
$\theta(h) = -h$ and $\tau :=\theta\tau_h$, so  that
$(\g,\tau)$ is a ncc symmetric Lie algebras with $\g = \g_s$. 
Let $\cD_+\subeq G/P^-$ be the open orbit of the base point under 
$H := K^h \exp(\fh_\fp)$. 
We endow the homogeneous space 
\[ M_{\cD_+} :=  \{g.\cD_+ \:  g \in G \},  \] 
consisting of subsets of $\cM$, with the inclusion order. 
Then the stabilizer subgroup $G^{\cD_+}$ of the base point is 
$H$.
The map $gH \mapsto g.\cD_+$ induces an isomorphism 
\[ (M_{\cD_+}, \subeq) \cong (G/H, C_\fq^{\rm max}(h)) \] 
of ncc symmetric spaces, where $C_\fq^{\rm max}(h)$ is the unique 
maximal $\Ad(H)$-invariant cone in $\fq$ containing $h$ in its interior.
\end{prop}

In this identification the set $\{x\in G/H\: x\ge eH\}$
is mapped to $\{s^{-1} \cD_+ \: s\in  \comp({\cD_+}) \}$ 
and $\{x\in G/H \: x\le eH\}$ is mapped to $\{s \cD_+ \: s\in   \comp({\cD_+}) \}$. In particular, $gH\ge eH$ is equivalent to $\cD_+\subset g.\cD_+$
  and $eH \ge gH$ to  $g.\cD_+\subset \cD_+$.

\begin{prf} This follows from Lemma \ref{lem:Compress}. 
\end{prf}

\begin{rem} (The Riemannian case)
  Let $(\g,\theta)$ be a Riemannian symmetric Lie algebra,
  i.e., $\g = \g_r$. Then
  $H = K$, $M = G/K$ and  $h=0$. Thus $G= P^-$ and $G/P^-$ is a single point.
    Hence $\comp(\cD_+) = G$ and $\cM_{\cD_+}$ is a single point.
    Therefore Riemannian summands cannot be permitted in
    Proposition~\ref{prop:setreal}.
\end{rem}

  \begin{ex}\label{ex:torus-new} Let $G=\SL_2(\R)$
    and $h = \frac{1}{2}\pmat{1 & 0 \\ 0 & -1}$.
    Then the canonical action of $G$ on $\bP_1(\R) = \bP(\R^2)
    \cong \bS^1 =  \R\cup \{\infty\}$ is given by
    \[\begin{pmatrix} a & b\\ c & d\end{pmatrix}.x
      = \frac{a x + b}{c x + d}\]
and the stabilizer of $0$ is 
\[P^- = \Big\{ \pmat{ a & 0 \\ c & a^{-1}} \: a \not=0, c \in \R \Big\}
  = \exp(\g_{-1}(h)) \rtimes G^h, \quad
  G^h = \Big\{\pmat{a & 0 \\ 0 & a^{-1}} \Big\} \cong \R^\times. \]
  The $1$-parameter group 
  \[ a_t = \begin{pmatrix} \cosh t & \sinh t \\ \sinh t & \cosh t
    \end{pmatrix}  \]
  fixes $\pm 1$ and the orbit $\cD_+ = H.0$ of $H = \{\pm a_t \: t \in \R\}$
  is the open unit interval $\cD_+ = (-1,1)$. 
  The maximal cone in $\fq$ is generated by $\Ad (H)\R_+ h$.

  Since elements of $\bP_1(\R)$ represent one-dimensional
  linear subspaces of $\R^2$ and $\SL_2(\R)$ acts transitively on
  triples of such subspaces, it follows easily that it acts
  transitively on the set of non-dense open intervals
  $I \subeq \bS^1$, the ordered space $G/H$ can be identified
  with the ordered set of open non-dense intervals in $\bS^1$.
\end{ex}

  \begin{ex} A special case of the above construction
    is the ``complex case'' where $H$ is a connected semisimple Lie group of
hermitian type contained in a complex Lie group $G$ with Lie algebra $\fh_\C = \fh \oplus i\fh$. 
Then $G/H$ is a ncc symmetric space. 
Let $\theta_H$ be a Cartan involution on $H$. Then
$\theta_H$ extends to a Cartan involution $\theta $ on $G$.
Denote the corresponding maximal compact
subgroup of $G$ by $K$. Then $H\cap K$ is a maximal compact subgroup of $H$
and the Riemannian symmetric space $H/H\cap K$ can be realized as
complex symmetric bounded domain $\cD_+ \subeq G/P^-$. 
Let $z_0\in \fz (\fh\cap \fk)$ be 
the element determining the
complex structure on $H/H\cap K$. 
Then $h=-iz_0$ is an Euler element in $\fq = i \fh$.  
Now \eqref{eq:boundedReal} is the Harish--Chandra realization
of $H/H\cap K$ as $\cD_+$ in $G/P^-$  
  (see \cite[p.~58]{Sa80} or
\cite[Ch.~VII]{He78} for details).
 
Suppose that the complex conjugation $\tau$ of $\g$ with respect to
  $\fh$ integrates to an involution $\tau^G$ on~$G$.
  This is the case if $G$ is simply connected or if $G=\Inn{\fg}$.
We then assume that $H = G^{\tau^G}_e$.   
If $G$ is simply connected, then $H = G^{\tau^G}$ is connected
and \cite[Thm.~VI.11]{HN95} implies that $H=G^{\cD_+}$, where
$G^{\cD_+}$ is the stabilizer of the base point $\cD_+$.

But in general, if $G$ is not simply connected, then
$G^{\cD_+}$ and $G^{\tau^G}$ may differ. 

As an example, consider $H=\PSL_2(\R)\subeq G=\PSL_2(\C)\cong \Inn(\g)$
and note that $\tau^G(g) = \tau g \tau$ in this case. 
Then
\[ G^{\tau^G} =\PSL_2(\C)^{\tau^G} \cong \PGL_2(\R)  \cong \Aut(\fsl_2(\R)),\] 
which is not connected because it also contains
the image of $\pmat{i & 0 \\ 0 & -i}$. 
The domain $\cD_+ = H.i \subeq  \C_\infty$ 
(the Riemann sphere) is the upper half plane and the stabilizer
subgroup of $\cD_+$ is 
\[ G^{\cD_+} = \PSL_2(\R) \not= \PGL_2(\R) = G^{\tau^G}.\] 
The reflections in $\GL_2(\R)$ exchange the two open $H$-orbits. 
\end{ex}

\begin{rem} \mlabel{rem:5.11}
The flag manifolds $\cM = G/P^- \cong K/K^h$ appearing in this section are 
compact symmetric spaces on which the maximal compactly embedded subgroup $K \subeq G$ 
acts by automorphisms. These spaces are called {\it symmetric R-spaces}.

Defining a symmetric R-space as a compact 
symmetric space $\cM$ which is a real flag manifold, 
Loos shows in \cite[Satz~1]{Lo85} that this implies the existence of 
an Euler element $h \in \cE(\g)$ such that 
$\cM \cong G/P^-$, so that 
$\cM \cong K/K\cap P^- = K/K^h$ as a Riemannian symmetric space
(see \cite[\S 7.2]{MNO23} for more details).

If $G$ is hermitian of tube type, then $\cM \cong K/K^h$ can be identified
  with the \v Shilov boundary of the corresponding bounded
  symmetric domain $\cD_G \cong G/K$, and this leads to a
  $G$-invariant causal structure on $\cM$.
  As $\dim Z(K) = 1$, with respect to the
  $K$-action, we have a natural $1$-parameter family of $K$-invariant
  Lorentzian structures on $\cM$. They correspond to $K^h$-invariant Lorentzian
  forms on $T_{eK^h}(\cM) \cong \fq_\fk = \fz(\fk) \oplus [\fh_\fk,\fq_\fk]$
  which are positive definite on $\fz(\fk)$ and negative definite
  on its orthogonal space $[\fh_\fk,\fq_\fk]$.
  \end{rem}

\section{Observer domains associated to modular geodesics} 
\mlabel{sec:6} 

In this section we associate to any modular causal geodesic $\gamma$ 
in an ncc semisimple symmetric space $M= G/H$
an {\it observer domain}  $W(\gamma)$. 
It is an open connected subset of $M$ invariant under the 
centralizer $G^h$ of the corresponding causal Euler element~$h$. 
We then show that, for $h \in C^\circ$, the domain 
$W(\gamma)$ coincides with the connected component 
$W \subeq  W_M^+(h)$ of the base point $eH$
of the corresponding positivity domain.
In Section~\ref{sec:7} we show that $W_M^+(h)$
is connected for $G = \Inn(\g)$, which implies that~$W = W_M^+(h)$ in
this case. 

\begin{defn} \mlabel{def:4.1}
  Let $(G,\tau^G, H,C)$ be a non-compactly causal
  symmetric Lie group and  $M = G/H$ be the corresponding
  ncc symmetric space. We assume that $\g = \g_s$, i.e., that
  $(\g,\tau)$ is a direct sum of irreducible ncc symmetric Lie algebras. 

\nin (a) We write $\leq$ for the order  on $M$ defined by the closed 
Olshanski semigroup $S = H \exp(C) = \exp(C)H$
which always exists because $\fz(\g) = \{0\}$
(\cite[Thm.~3.1]{La94} or Theorem~\ref{thm:lawson}in Appendix~\ref{sec:law}) 
via 
\[ g_1H \leq g_2 H \quad \mbox{ if } \quad g_1^{-1} g_2  \in SH/H = \Exp_{eH}(C) \] 
and write order intervals as 
\[ [x,y] = \{ z \in M \: x \leq z \leq y\} 
= \up x \cap \down y, \]
where 
\[ \up x = \{ z \in M \: x \leq z\} \quad \mbox{ and } \quad 
\down y = \{ z \in M \: z \leq y\}.\] 

\nin (b) A subset $X\subeq M$ is called {\it order convex} if 
\[ [a,b] \subeq X \quad \mbox{ for } \quad  a,b \in X.\] 
As the intersection of order convex subsets is order convex, 
we can defined the {\it order convex hull} 
\[ \oconv(D) := \bigcap \{ D' \subeq M \: D \subeq D', D' \ 
  \mbox{ order convex} \}.\]
Clearly $\oconv(D)$ is the smallest order convex subset of $X$
  containing $D$.

\nin (c) For a modular geodesic 
$\gamma \: \R \to M$, we call 
\[ W(\gamma) 
:= \down \gamma(\R) \cap \up \gamma(\R) 
= \bigcup_{t < s} [\gamma(t), \gamma(s)] \] 
the {\it  observer domain associated to $\gamma$}. Note that this 
domain depends on the cone $C \subeq \fq$ specifying the order on~$M$. 
\end{defn}

\begin{lem} \mlabel{lem:wgammaprops} 
The subset $W(\gamma)$ has the following properties: 
  \begin{itemize}
  \item[\rm(a)]  $W(\gamma)\subeq M$ is  open and connected. 
  \item[\rm(b)] $W(\gamma) = \oconv(M^h_{eH})$ for
    $M^h_{eH} = \Exp_{eH}(\fq_\fp)$.  
  \item[\rm(c)] Suppose that $H_K = K^h$ and $C = C_\fq^{\rm max}$ and identify 
$M = G/H$ with $M_{\cD_+}$ {\rm(Proposition~\ref{prop:setreal})}. Then 
\begin{equation}
  \label{eq:wgammad}
 W(\gamma) = \{ g.{\cD_+} \: 0 \in g.{\cD_+}, g.{\cD_+} \ \mbox{ bounded}\}
\end{equation}
and this domain is $G^h$-invariant.   
  \end{itemize}
\end{lem}

\begin{prf} (a)  
To see that $W(\gamma)$ is open, we first observe that 
$\gamma(s) \in (\up \gamma(t))^\circ$ for 
$t < s$. For real numbers $t_j \in \R$ with 
$t_1 < t_2 < t_3 < t_4$, this implies that 
\[ [\gamma(t_2), \gamma(t_3)] \subeq [\gamma(t_1), \gamma(t_4)]^\circ.\] 
This shows that $W(\gamma)$ is open. 

To see that $W(\gamma)$ is connected, we recall that the order 
on $M$ is globally hyperbolic, in particular all order 
intervals $[x,y]$ are compact. As all elements $z \in [x,y]$ 
lie on causal curves from $x$ to $y$ (\cite[Thm.~4.29]{HN93}), the order intervals 
are pathwise connected. As an increasing union of the order intervals 
$[\gamma(-n), \gamma(n)]$, the wedge domain $W(\gamma)$ is connected. 

\nin (b)  Order intervals are convex and directed unions of 
convex sets of convex. Therefore 
\[ W(\gamma) = \bigcup_{t < s} [\gamma(t), \gamma(s)] \] 
is convex, whence $W(\gamma)  = \oconv(\gamma(\R)).$ 

From the fact that 
$h$ is central in $\fh_\fk + \fq_\fp$, it easily follows that, 
in the symmetric space $M^h_{eH} = \Exp_{eH}(\fq_\fp)$ the geodesic line  
$\gamma(\R)$ is cofinal in both directions because we have in $\fq$: 
\[ \bigcup_{s < t} (sh + C) \cap (th - C) \supeq \fq.\] 
For $x \in \fq_\fp$, we thus find 
$s,t \in \R$ with $x  \in sh + C^\circ$ and $x \in st - C^\circ$. Then 
\[ \Exp(sh) < \Exp(x) < \Exp(th) \] 
in $M^h_e$. This implies that 
\[ W(\gamma) \supeq \Exp_{eH}(\fq_\fp) = M^h_{eH} = G^h_e.eH 
\supeq \gamma(\R).\] 
This completes the proof.   

\nin (c) The modular group 
acts on $\cB \cong  N^+.eP^- \subeq G/P^-$ by $\exp(th).x = e^tx$. 
Therefore $\gamma(t) = e^t {\cD_+}$ enlarges ${\cD_+}$ for $t > 0$ and shrinks 
${\cD_+}$ for $t < 0$ (Theorem~\ref{thm:contheo}). 
As $\gamma$ is strictly increasing, this implies that 
\[ \down \gamma(\R) = \{ g.{\cD_+} \: (\exists t \in \R) g.{\cD_+} \subeq e^t {\cD_+} \} 
= \{ g.{\cD_+} \: g.{\cD_+}\ \mbox{ bounded in } \ \cB\}.\]
Further 
\[ \up \gamma(\R) = \{ g.{\cD_+} \: (\exists t \in \R) g.{\cD_+} \supeq e^t {\cD_+} \} 
= \{ g.{\cD_+} \: 0 \in g.{\cD_+}\},\]
so that \eqref{eq:wgammad} follows. 
As any $g \in G^h = P^+ \cap P^-$ acts by linear maps 
on the Bruhat cell $\cB \cong \g_1(h)$, 
\eqref{eq:wgammad} implies that $G^h$ leaves the set 
$W(\gamma)$ of all bounded domains $g.{\cD_+}$ containing 
$0$ invariant.
\end{prf}

\begin{ex} \mlabel{ex:8.15} (de Sitter space) We consider de Sitter space 
\[ M = \dS^d = \{ (x_0,\bx) \in \R^{1,d} \: \beta(x,x) = -1\}, 
\quad \mbox{ where } \quad  \beta(x,y) = x_0 y_0 - \bx\by  \] 
is  the canonical Lorentzian form on $\R^{1,d}$ 
(cf.\ Section~\ref{subsec:dS1}). 
Here 
\[ G = \SO_{1,d}(\R)^\up = \SO_{1,d}(\R)_e, \quad H = G_{\be_1} = \SO_{1,d-1}(\R)^\up \]
and 
\[ C \subeq T_{\be_1}(M) = \be_1^\bot \quad \mbox{ given by } \quad 
 C = \{  (x_0,\bx) \: x_1 = 0, x_0 \geq 0, x_0^2 \geq \bx^2\}.\]
We claim that, for the modular geodesic 
\[ \gamma(t) = \cosh(t) \be_1 + \sinh(t) \be_0  = e^{th}\be_1, \] 
we have 
\begin{equation}
  \label{eq:desitterwedge}
 W(\gamma) = W_{\dS^d}(h) = \{ x \in \dS^d \: x_1 > |x_0|\} 
= W_R \cap \dS^d,
\end{equation}
where $W_R = \{ (x_0,\bx) \: x_1 > |x_0|\}$ 
(cf. Appendix~D in \cite{NO23b}). 
As the right wedge $W_R \subeq \R^{1,d}$ is causally complete, we 
clearly have $W(\gamma) \subeq W_R \cap \dS^d = W_{\dS^d}(h)$. 
For the converse inclusion, let $x \in W_R$. We have to find a 
$t \in \R$ with $x \leq \gamma(t)$, i.e., 
\[ x_0 < \gamma(t)_0 = \sinh(t) \] 
and 
\[ 0 < \beta(\gamma(t) -x, \gamma(t) -x) 
= (\sinh(t)-x_0)^2 - (\cosh(t)-x_1)^2 - x_2^2 - \cdots- x_d^2.\]
Since $\beta(\gamma(t),\gamma(t)) = -1$, we obtain for the right hand side 
\[ \beta(\gamma(t)-x,\gamma(t)-x) 
= \beta(\gamma(t), \gamma(t)) - 2 \beta(\gamma(t),x\big) + \beta(x,x)
= -1 - 2\beta(\gamma(t),x)+ \beta(x,x).\] 
Further 
\[ -2 \beta(\gamma(t),x) = 2 x_1 \cosh(t) - 2 x_0 \sinh(t) 
\approx e^t (x_1 - x_0) \quad \mbox{ for } \quad t \> 0,\] 
and if $x_1 > |x_0|$, this expression is arbitrarily large 
for $t \to \infty$. 
This shows that $W_R \subeq \down \gamma(\R)$, and we likewise see that 
$W_R \subeq \up \gamma(\R)$. 
\end{ex}

\begin{prop} \mlabel{prop:inc1} 
  If $H_K = K^h$ and $C = C_\fq^{\rm max}$, then
  \begin{itemize}
  \item[\rm(a)] $W(\gamma) \subeq W = W_M^+(h)_{eH}$. 
  \item[\rm(b)] $h + \cD_\g \subeq\fh + C^\circ$.
  \end{itemize}
\end{prop}

\begin{prf} If $gH \in W(\gamma) \subeq G/H$, then 
the corresponding subset $g.{\cD_+} \subeq \cB$ is convex 
by Theorem~\ref{thm:contheo}, and it contains~$0$ by \eqref{eq:wgammad}. 
Therefore the curve 
\[ \eta \: \R \to M,\quad \eta(t) := \exp(th)gH \] 
is increasing because $t \mapsto e^t g.{\cD_+}$ is an increasing family 
of subsets of $\cB$. The invariance of the order thus implies that 
\[ g^{-1}.\eta'(0) = p_\fq(\Ad(g)^{-1}h) \in C_\fq^{\rm max}.\] 

We also know that $g.{\cD_+} \in P^+.{\cD_+}$ (Theorem~\ref{thm:contheo}
and Lemma~\ref{lem:wgammaprops}(c)), so that there exist 
$g_1 \in G^h$ and $y \in \g_1(h)$ with 
$g.H = g_1 \exp(y).H$. Thus 
\begin{equation}
  \label{eq:x1}
  \Ad(g)^{-1}h \in \Ad(H) e^{-\ad y}h \in \fh + C_\fq^{\rm max},
\end{equation}
and therefore 
\[ e^{-\ad y}h = h - [y,h] = h + y  \in \fh + C_\fq^{\rm max}.\] 

Recall the definition of $\cD_\g$ in \eqref{def:Dg}. The condition 
\[eP^- \in g.{\cD_+} = g_1 \exp(y).{\cD_+} = g_1.\exp(y + \cD_\g)P^-\]
is equivalent to $-y \in \cD_\g = - \cD_\g$, showing that 
\begin{equation}
  \label{eq:wgamma3}
 W(\gamma) = G^h \exp(\cD_\g).\cD_+ \subeq M_{\cD_+}
\end{equation}
(cf.\ Lemma~\ref{lem:wgammaprops}(c)). 
We therefore derive from \eqref{eq:x1}
that $h + \cD_\g \subeq \fh + C_\fq^{\rm max},$ 
and since $h \in C_\fq^{\rm  max,\circ}$ and ${\cD_+}$ is starlike with respect
to~$0$, 
we obtain 
\begin{equation}
  \label{eq:hdw}
  h + \cD_\g \subeq\fh + C_\fq^{\rm max, \circ}.
\end{equation}
We thus obtain 
$\Ad(g)^{-1}.h \in \fh + C_\fq^{\rm max,\circ}$, i.e.,  $gH \in W_M^+(h)$. 
This shows that $W(\gamma) \subeq W_M^+(h)$, and the connectedness of 
$W(\gamma)$ (Lemma~\ref{lem:wgammaprops}(a)) yields $W(\gamma) \subeq W$. 
\end{prf}

\begin{rem} From \eqref{eq:wgamma3} it follows that, as a subset of $M$, 
  \begin{equation}
    \label{eq:wgam}
    W(\gamma)
    = G^h \exp(\cD_\g).H     = G^h_e \exp(\cD_\g).H.
  \end{equation}
For the quotient map $q \: G \to G/H$, this means that 
\[ q^{-1}(W(\gamma)) = G^h \exp(\cD_\g) H \subeq G.\] 
This is a $G^h\times H$-invariant domain in $G$ specified by its intersection 
with the abelian subgroup $N^+ = \exp(\g_1(h))$; see
\cite[Rem.~6.2]{MNO23}. 
\end{rem}

Combined with Theorem~\ref{thm:connwedgedom} below, that
asserts the connectedness of $W_M^+(h)$, the following
result implies that $W_M^+(h) \subeq W(\gamma)$.
  
\begin{prop} If $H_K = K^h$ and $C = C_\fq^{\rm max}$, then 
$W \subeq W(\gamma)$. 
\end{prop}

\begin{prf} As both sides are $G^h_e$-invariant (Lemma~\ref{lem:wgammaprops}), 
  the Positivity Domain Theorem (Theorem~\ref{thm:pos-polar}) implies that
  we have to verify the inclusion 
\[ \Exp_{eH}(\Omega_{\fq_\fk}) \subeq W(\gamma).\] 
Invariance of both sides under $(H_K)_e$ and 
$\Ad((H_K)_e)\ft_\fq = \fq_{\fk}$ further 
reduce the problem to the inclusion 
\begin{equation}
  \label{eq:finalincl}
 \Exp_{eH}(\Omega_{\ft_\fq}) \subeq W(\gamma).
\end{equation}

To this end, we use the Lie subalgebra $\fl \subeq \g$ generated by 
$h$ and $\ft_\fq$ (Proposition~\ref{prop:testing}). 
Then $[\fl,\fl] \cong \fsl_2(\R)^s$
and $\ft_\fq \cong \so_2(\R)^s$. 
This reduces the verification of the inclusion \eqref{eq:finalincl} to 
the case where $\g = \fsl_2(\R)^s$, $\fh = \so_{1,1}(\R)^s$ 
and $\ft_\fq \cong \so_2(\R)^s$. 

As this is a product situation, it suffices to consider the case where 
\[ \g= \fsl_2(\R)\supeq \fh =  \so_{1,1}(\R), \quad \ft_\fq = \so_2(\R)
\quad \mbox{ and } \quad 
h = \frac{1}{2}\pmat{1 & 0 \\ 0 & -1}.\] 
By \eqref{eq:wgam}, we  have to show that 
\begin{equation}
  \label{eq:sl2tilde-fact} 
 \exp(tx) \in G^h_e \exp(\cD_\g) H
\quad \mbox{ for } \quad 
|t| < \pi/2 \quad \mbox{ and } \quad 
x = \frac{1}{2} \pmat{0 & -1 \\ 1 & 0}.
\end{equation}
We identify $\fsl_2(\R)$ with $3$-dimensional Minkowski space $\R^{1,2}$, via 
\[ \be_0 := \frac{1}{2} \pmat{0 & -1 \\ 1 & 0}, \quad
\be_1 := \frac{1}{2} \pmat{0 & 1\\ 1 & 0}, \quad 
\be_2 := h =\frac{1}{2} \pmat{1 & 0\\ 0 & -1}.\] 
In the centerfree group $G := \Inn(\g) \cong \SO_{1,2}(\R)_e$, we have 
\[ K := G_{\be_0} \cong \SO_2(\R) \quad \mbox{ and } \quad K^h = \{e\} = H_K,\] 
so that $H := G_{\be_1} = \exp(\fh_\fp) = \SO_{1,1}(\R)_e$ is connected. 
Therefore 
$G/H \cong G.\be_1 = \dS^2$ (de Sitter space) 
and $\exp(tx)H$ corresponds to 
\[ \exp(tx).\be_1= \cos(t)\be_1 + \sin(t) \be_2.\] 
Now $|t| < \pi/2$ implies $\cos(t) > 0$, hence that 
\begin{equation}
  \label{eq:exptx}
 \exp(tx).\be_1\in W_{\dS^2}(\gamma) \quad \mbox{ for } \quad 
\gamma(t) = \cos(t) \be_1 + \sin(t) \be_0
\end{equation}
(Example~\ref{ex:8.15}). 
We write elements of ${\cD_+}$ as $y = s e$, $|s| < 1$ (see~Example~\ref{ex:sl2}). 
Then $\exp(y).\be_1$ corresponds to 
\begin{align*}
 e^{\ad y}.\frac{1}{2}(e+f)
& =  \frac{1}{2}e^{\ad se}(e+f)
 =  \frac{1}{2}(e+ e^{\ad se} f)
 =  \frac{1}{2}\big(e+ f + s[e,f] + \frac{s^2}{2}[e,[e,f]]\big) \\
& =  \frac{1}{2}(e+ f) + sh + \frac{s^2}{2}[e,h] 
 =  \frac{1}{2}(e+ f) + sh - \frac{s^2}{2}e, 
\end{align*}
so that 
\[ \exp(y).\be_1 
 = \be_1 + s \be_2  - \frac{s^2}{2}(\be_1 - \be_0)
 =\Big(\frac{s^2}{2}, 1 - \frac{s^2}{2}, s\Big).\] 
This elements lies in the wedge domain $W_{\dS^2}(h)$ if and only if 
$1-s^2/2 > s^2/2$ (Example~\ref{ex:8.15}), 
which is equivalent to $|s| < 1$. Then its $G^h_e$-orbit contains 
the element $(0, \sqrt{1-s^2}, s)$. 
For $|t| < \pi/2$, the element $\exp(tx).\be_1$ is of this form, showing that 
$\exp(tx) \in G^h \exp(y) H$. This completes the proof.
\end{prf}

Combining the preceding two propositions, we get the main result 
of this section. It shows that the observer domain
$W(\gamma)$ coincides with a connected component of
the positivity domain~$W_M^+(h)$. This result provides two complementary
perspectives on this domain. 

\begin{thm} \mlabel{thm:equality}
  {\rm(Observer Domain Theorem)} 
  Let $(\g,\tau,C)$ be a non-compactly causal semisimple 
  symmetric Lie algebra with causal Euler element
  $h \in C^\circ \cap \fq_\fp$ with $\tau = \tau_h \theta$ 
  and let  $G$ be a connected Lie group with Lie algebra $\g$
  and $H := K^h \exp(\fh_\fp)$.
  If $C = C_\fq^{\rm max}$, then $W =  W(\gamma)$. 
\end{thm}

We can even extend this result to coverings: 
\begin{cor} \mlabel{cor:equality2}
If $H' \subeq H = K^h\exp(\fh_\fp)$  is an open subgroup and $C = C_\fq^{\rm max}$, then 
$W' :=  W_{M'}^+(h)_{eH'} = W(\tilde\gamma)$ holds in $M'= G/H'$ for 
$\tilde\gamma(t) := \Exp_{eH'}(th)$. 
\end{cor}

\begin{prf}  Let $q \: M' = G/H' \to G/H \cong M_{\cD_+}$ be the canonical 
  equivariant covering from \cite[Lemma~7.11]{MNO23}.

First we show that $W' \subeq M'$ is order convex. So let 
$x\leq y \leq z$ in $M'$ with $x,z \in W'$ 
and let $\eta \: [0,2] \to M'$ be a causal curve with 
\[ \eta(0) =x, \quad \eta(1) = y, \quad \eta(2) = z.\] 
Then $q(\eta(t)) \in [q(x),q(z)] \subeq W$ for $t \in [0,2]$ holds 
because $W = W(\gamma)$ is order convex in $M$. 

As $W$ is contractible by
Theorem\ref{thm:pos-polar}(b), 
it is in particular simply connected. 
Therefore $q^{-1}(W)$ is a disjoint union of open subsets 
$(W_j')_{j \in J}$ mapped by $q$ diffeomorphically onto $W$. 
By definition, $W'$ is one such connected component, so that 
\[ q_W := q\res_{W'} \: W' \to W \] 
is a diffeomorphism. 
Therefore $\eta$ is the unique continuous lift of $q \circ \eta$ in $M'$, 
hence contained in~$W'$. This implies that $y \in W'$, so that 
$W'$ is order convex. 

As $q_W \: W' \to W$ 
is an isomorphism of causal manifolds, it also is an order isomorphism. 
Finally $W(\gamma) = \oconv(\gamma(\R)) = W$ implies that 
$W'(\gamma) = \oconv(\tilde\gamma(\R)) = W'$. 
\end{prf}

  \begin{rem} It is not clear to
    which extent $W(\gamma)$ depends on the specific cone $C$.
    In particular it would be interesting to see if
    the minimal and maximal cones lead to the same domain 
    $W(\gamma)$.
    We have already seen that the positivity domain 
    $W_M^+(h)$ depends non-trivially on the cone $C$ (\cite[Ex.~6.8]{MNO23})
    so one may expect that this is also the case  for $W(\gamma)$.
  \end{rem}

\begin{lem} \mlabel{lem:tauWinv}
The involution $\tau_M$ on $M$ defined by
    $\tau_M(gH) = \tau^G(g)H$ satisfies
    \begin{equation}
      \label{eq:tauWM}
      \tau_M(W_M^+(h)) = W_M^+(h)
      \quad \mbox{ and } \quad \tau_M(W(\gamma)) = W(\gamma).
                \end{equation}
\end{lem}

\begin{prf}  (a) The condition $gH \in W_M^+(h)$ is equivalent to
    $\Ad(g)^{-1}h \in \cT_C$ by \eqref{eq:tauWM}, and this implies that
    \[ \Ad(\tau(g))^{-1}(-h) = \tau(\Ad(g)^{-1}h) \in \tau(\cT_C)
    = -\cT_C,\]
    so that $\Ad(\tau(g))^{-1}h \in \cT_C$, i.e.,
    $\tau_M(gH) \in W_M^+(h)$. As $\tau_M$ is an involution, it follows that
    $\tau_M(W_M^+(h)) = W_M^+(h)$.

    \nin (b) As $\tau(C) = - C$, the involution $\tau_M$ reverses the causal
    structure on $M$. Moreover,
    $\tau_M(\gamma(t)) = \gamma(-t)$, so that
    \[ \tau_M(W(\gamma))
    = \bigcup_{t < s} [\tau_M(\gamma(s)), \tau_M(\gamma(t))]
    = \bigcup_{t < s} [\gamma(-s), \gamma(-t)] = W(\gamma).   \qedhere   \] 
\end{prf}

We have seen above that, for the  modular geodesic
$\gamma(t) = \Exp_{eH}(th)$ in $M$, we have $W(\gamma) = W$.
The modular geodesic $\gamma$ is a specific orbit of the modular
flow inside $W$. Now we show that all other $\alpha$-orbits in $W$
lead to the same ``observer domain''. 

\begin{prop}
  Let $m \in W$ and consider the
  curve
  \[ \beta \:  \R \to W, \quad \beta(t) = \alpha_t(m) = \exp(th).m.\]
  Then
  \begin{equation}
    \label{eq:W}
    W = W(\beta) = \bigcup_{s < t} [\beta(s), \beta(t)].
  \end{equation}
\end{prop}

\begin{prf}  Using the subset realization of
  $M = G/H$ as
$M_{\cD_+} = \{ g.\cD_+ \: g \in G \}$ 
from Proposition~\ref{prop:setreal}, we have
  \[  W(\gamma) = \{ g \cD_+ \: 0 \in g.\cD_+, g.\cD_+\
    \mbox{ bounded in } \exp(\g_1(h)).P^-\}\]
    (Lemma~\ref{lem:wgammaprops}(c))
    and $W = W(\gamma)$ by Theorem~\ref{thm:equality}.  
  So we can write
  \[ \beta(t) = e^t.\cD' \quad \mbox{ for some } \quad \cD' \in W(\gamma).\]
  As $\beta(\R) \subeq W(\gamma)$, the order convex hull $W(\beta)$ of
  $\beta(\R)$ is contained in $W(\gamma) = W$. 
  To verify the converse inclusion, let $\cD'' \in W$. Then
  $0 \in \cD''$, and since $\cD'$ is bounded, there exists a $t \in \R$ with
  $\beta(t) \subeq \cD''$. Likewise the boundedness of $\cD''$ implies the
  existence of some $s \in \R$ with $\cD'' \subeq \beta(s)$.
  Hence $\cD'' \in [\beta(t),\beta(s)] \subeq W(\beta)$. This shows that
  $W \subeq W(\beta)$, and hence equality in \eqref{eq:W}.   
\end{prf}

\begin{rem} A similar result also holds in Minkowski space.
  If $\bx \in W_R = \{ \by \in \R^{1,d} \: y_1 > |y_0| \}$ and
  \[ \beta(t) = e^{th} \bx
    = (\cosh(t) x_0 + \sinh(t) x_1, \cosh(t) x_1 + \sinh(t) x_0, x_2,
    \ldots, x_d),\]
  then any other element $\by \in W_R$ satisfies
$\by \in [\beta(t), \beta(s)]$ 
  for suitable $t < s$, i.e., 
  $\by - \beta(t) \in V_+$ and $\beta(s) - \by \in V_+$.
  In fact, $\beta(t)_0 \sim e^t(x_0 + x_1) \to \infty$ for $t \in \infty$
  and $\beta(t_0) \sim e^{-t}(x_0 - x_1) \to -\infty$ for $t \to -\infty$.
  Moreover, for $s \to \infty$
  \begin{align*}
&    (\cosh(s) x_0 + \sinh(s) x_1 - y_0)^2 - (\cosh(s) x_1 + \sinh(s) x_0 - y_1)^2 \\
    &    \sim \Big(e^s \frac{x_0 + x_1}{2} -y_0\Big)^2
    - \Big(e^s \frac{x_0 + x_1}{2} -y_1\Big)^2
      \sim  e^s(x_0 + x_1)(y_1-y_0) \to \infty
  \end{align*}
  and, for $t \to -\infty$, 
  \begin{align*}
 &  (\cosh(t) x_0 + \sinh(t) x_1 - y_0)^2 - (\cosh(t) x_1 + \sinh(t) x_0 - y_1)^2 \\
    &   \sim \Big(e^{-t} \frac{x_0 - x_1}{2} -y_0\Big)^2
    - \Big(e^{-t} \frac{x_1 -x_0}{2} -y_1\Big)^2
    \sim e^{-t} (x_0 - x_1) (-y_0)  - e^{-t}(x_1 - x_0) (-y_1) \\
&    = e^{-t} (x_1 - x_0) (y_0 + y_1) \to\infty.
  \end{align*}
  This shows that $W(\beta) = W_R$ for all integral curves of the modular flow
  in $W_R$. 
  \end{rem}

\begin{rem} On the de Sitter space $M = \dS^d \subeq \R^{1.d}$, the involution
    $\tau_h$ can be implemented naturally by
    \[ \tau_{h,M}(x) = (-x_0, -x_1, x_2, \ldots, x_d). \]
    This involution does not fix the base point $\be_1$,
    it reverses the causal structure and it commutes with modular flow.
    Accordingly, we have the relation
    \[ \tau_{h,M}(W^+(h)) = W^+(-h).\] 
    As we shall see in the next section, such a relation can only 
    be realized because $-h \in \Ad(G)h$, i.e., the direction
    of the boost can be reversed by an element of $G$.
    If $-h \not\in \Ad(G)$ ($h$ is not symmetric),
    then we shall see in Corollary~\ref{cor:7.3}
    below that $W^+(-h) = \eset$, so that there is no
    involution on $M$ mapping $W^+(h)$ to $W^+(-h)$.

    However, as $\tau_h = \tau \theta$ (as involutions on $\g$),
    and there are natural implementations $\tau_M$ and $\theta_M$ on
    $M = G/H$, both fixing the base points,
    the involution $\tau_M\theta_M$ implements
    the involution $\tau_h$ on $M$ and fixes the base point,
    but it also fixes the wedge region 
    \[ \tau_M\theta_M(W^+(h)) = W^+(h) \]
    because it preserves $h$ and the causal structure.
    This is not desirable because we would prefer that $\tau_h$
    maps $W^+(h)$ to some ``opposite'' wedge region (cf.~\cite{MN21}).
    Possible ways to resolve this problem and ideas how to implement
    locality conditions on non-compactly causal symmetric spaces
    are briefly discussed in \cite[\S 4.3]{MNO23}.
\end{rem}

\section{Existence of positivity domains for Euler elements} 
\mlabel{sec:7}

In this section we show that, for the maximal cone  
$C= C_\fq^{\rm  max}$ and a simple Lie algebra $\g$,  the
  real tube domain $\cT_C = \fh + C^\circ$ intersects the set
  $\cE(\g)$ of Euler elements in a connected subset  (Theorem~\ref{thm:x.1}). 
This implies that, for an Euler element $h' \in \g$, the positivity
  domain $W_M^+(h')$ is non-empty if and only if $h'$ and $h$
  are conjugate   (Corollary~\ref{cor:7.3}). 

\begin{thm} \mlabel{thm:x.1}
 Suppose that $(\g,\tau,C)$ is an irreducible simple ncc
 symmetric Lie algebra with $C = C^{\rm max}_\fq$, $\cT_C := \fh + C^\circ$, 
 $G = \Inn(\g)$, $H = K^h \exp(\fh_\fp)$ and $M = G/H$.
 Then $\cE(\g) \cap \cT_C$ is connected and a subset of $\cO_h$.
 More precisely, 
\begin{equation}
  \label{eq:ident1a}
  \cE(\g) \cap \cT_C  = \cO_h \cap \cT_C   = \Ad(H_e)(h + \cD_\g),
\end{equation}
where $\cD_\g = \{ u \in \g_1(h) \:  \|u\| < 1\}$
is the  open unit ball for which $\exp(\cD_\g) P^- = H.P^- \subeq G/P^-.$ 
\end{thm}

\begin{prf} We recall from Proposition~\ref{prop:setreal}
the open subsets $\cD_\pm := H.eP^\mp\subeq G/P^\mp$ which are the open
orbits of the base point under $H = K^h \exp(\fh_\fp)$.
Then
\[ \comp(\cD_\pm) = H \exp(\mp C) \]
follows from Proposition~\ref{prop:setreal}, applied to
the causal Euler element $h$ and its negative.
These semigroups have the Lie wedges 
\[ \L(\comp(\cD_\pm)) = \fh  \mp C. \] 

Let $x \in \cE(\g) \cap \cT_C$ for
$\cT_C = \fh + C^\circ = \L(\comp(\cD_-))^\circ$.
We then have $s_t := \exp(tx) \in \comp(\cD_-)^\circ$ for $t > 0$.
We conclude that $s_t(\oline{\cD_-}) \subeq \cD_-$
and that there exists a complete metric on $\cD_-$ for which
each $s_t$ is a strict contraction
(cf.\ \cite[Thm.~II.4]{Ne01}),\begin{footnote}
  {This reference deals with bounded symmetric domains in complex
    spaces, but $\cD$ can be embedded into such a domain $\cD_\C$ by
    embedding $\g \into \g_\C \cong \g^c_\C$. If $C_{\g^c} \subeq \g^c$
    is an invariant cone with $C = \g \cap i C_{\g^c}$, then
    $(\g,\tau, C) \into (\g_\C, \tau_{\g_\C}, i C_{\g^c})$ is a causal
    embedding and $\cD_+ = H.eP^- \subeq G^c.eP^-_\C = \cD_+^\C$ is a real form
    of a complex bounded symmetric domain $\cD_+^\C$;
    see \cite[Lem. 1.4]{Ol91} or \cite[Lem 5.1.11]{HO97} for more details.
    }
\end{footnote}
so that the Banach Fixed Point
Theorem implies the existence of a unique attracting fixed point
$m_- \in \cD_-$ for the vector field $X^{G/P^+}_x \in \cV(G/P^+)$
defined by~$x$. We now have
\[ m_-  \in \cD_- = H.eP^+ = H_e.eP^+.\]
Hence there exists
$g_1 \in H_e$ with $g_1.m_- = eP^+$, and thus
\begin{equation}
  \label{eq:ab1x}
  y := \Ad(g_1) x \in \fp^+ = \g_1(h) \rtimes \g_0(h).
\end{equation}
Then $y \in \cT_C \cap \fp^+$ is an Euler element, and a similar
argument shows that the vector field
$X^{G/P^-}_y$ has a unique repelling fixed point $m_+ \in \cD_+$.
So $m_+ = \exp(-z) P^-$ for some $z \in \g_1(h)$,
and $\exp(z).m_+ = eP^-$. 
Hence the base point $eP^-\in G/P^-$
is a repelling fixed point of the Euler element
$y' := e^{\ad z} y \in \g_0(h)$,
and $eP^+$ is an attracting fixed point in $G/P^+$.
The attracting and repelling properties of the fixed points imply that
\[ \g_1(h) \subeq \g_1(y') \quad \mbox{ and } \quad 
\g_{-1}(h) \subeq \g_{-1}(y'), \]
so that we also have
\[ \g_0(h) = [\g_1(h), \g_{-1}(h)] \subeq \g_0(y').\]
As $h$ and $y'$ are Euler elements, this entails that
$\g_\lambda(h) = \g_\lambda(y')$ for
$\lambda = -1,0,1$. 
This shows that $\ad h = \ad y'$ and hence that~$y' = h$ because
$\fz(\g) = \{0\}$.

We conclude that 
\[ x = \Ad(g_1)^{-1} y = \Ad(g_1)^{-1} e^{-\ad z} h \quad \mbox{
    with }\quad g_1 \in H_e, z \in \cD_\g.\]
Conversely, we have seen in Proposition~\ref{prop:inc1} that
\begin{equation}
  \label{eq:inc4}
  e^{\ad \cD_\g} h = h + \cD_\g \subeq \cT_C.
\end{equation}
We finally obtain~\eqref{eq:ident1a}.
\end{prf}

\begin{rem} Note that the preceding proof is based on the natural
  embedding
  \[ \cO_h \cong G/G_h \to G/P^- \times G/P^+ \]
  which maps the Euler element $\Ad(g)H$ to
  $(m_+,m_-)$, where $m_+$ is the unique repelling fixed point of the
  flow defined by $h$ in $G/P^-$ and $m_- \in G/P^+$ is the unique
  attracting fixed point.   
\end{rem}

\begin{cor} \mlabel{cor:7.3} {\rm(The set of positivity domains in $M$)} 
  If $h_1 \in \cE(\g)$ is an Euler element for which
  the positivity domain
  \[ W_M^+(h_1) = \{ m \in M = G/H \: X^M_{h_1}(m) \in C_m^\circ \} \]
  is non-empty, then there exists a $g \in G$ with $h_1 = \Ad(g)h$ and
 \[ W_M^+(h_1) = g.W_M^+(h).\]   
\end{cor}

\begin{prf} As $X^M_{h_1}(g_1H) \in C_{g_1H}^\circ$ is equivalent to
  $\Ad(g_1)^{-1} h_1 \in \fh + C^\circ$ by
  (see Lemma~\ref{lem:wm+}), 
  Theorem~\ref{thm:x.1} implies that $h_1 = \Ad(g) h\in \cO_h$ for
  some $g \in G$. The relation $W_M^+(h_1) = g.W_M^+(h)$ now follows
  directly from the definitions.   
\end{prf}

The preceding corollary shows that any wedge domain of the type
$W_M^+(h_1)\subeq M$, $h_1 \in \cE(\g)$, is a $G$-translate of the wedge domain
$W_M^+(h)$, where $h \in C^\circ \cap \fq_\fp$ is a causal Euler element.
So the action of $G$ on the ``wedge space'' $\cW(M)$ of $M$ is transitive. 

\begin{cor} If the causal Euler element $h$ is not symmetric,
  then $W_M^+(-h) = \eset$.   
\end{cor}

\begin{rem} (Extensions to the non-simple case) 
If $(\g,\tau)$ is a direct sum of irreducible
ncc symmetric Lie algebra $(\g_j,\tau_j)$
and $h = \sum_j h_j$ accordingly, 
then
\[ C_\fq^{\rm max}(h) = \prod_j C_{\fq_j}^{\rm max}(h_j) \]
(cf.\ \eqref{eq:cmax-gen}). 
Projecting to the ideals $\g_j$, we obtain with
Theorem~\ref{thm:x.1}
for $C = C_\fq^{\rm max}(h)$ and 
$C_j = C_{\fq_j}^{\rm max}(h_j)$ the relation 
\begin{equation}
  \label{eq:eulprod}\cE(\g) \cap \cT_C
\subeq \prod_j \cE(\g_j) \cap \cT_{C_j}
\subeq \prod_j \cO_{h_j} = \cO_h.
\end{equation}
Further,
\[ \cO_h \cap  \cT_C = \prod_j \cO_{h_j} \cap \cT_{C_j} \]
and $\cD_\g = \prod_j \cD_{\g_j}$ imply
\eqref{eq:ident1a} for this case.

Note that the situation corresponds to
$\g = \g_s$ (see \eqref{eq:qs}). 
In the general situation, where we assume only that all ideals of
$\g$ contained in $\h$ are compact, we have
\[ \g = \g_k \oplus \g_r \oplus \g_s,\]
where $\g_k \subeq \fh$ is compact, $\g_r$ is a direct sum of
Riemannian symmetric Lie algebras and $\g_s$ is a direct sum of irreducible
ncc symmetric Lie algebras. All Euler elements are contained in~$\g_r + \g_s$. 
If $\g$ is only reductive, we assume
$\fz(\g) \subeq \g^{-\theta}$, so that $\fz(\g) \subeq \g_r$.
Then $h = h_r + h_s$ and 
\[ C_\fq^{\rm max}(h)= \fq_r \oplus C_{\fq_s}^{\rm max}(h_s).\] 
We conclude that
\[\cE(\g) \cap \cT_C \subeq \cE(\g) \cap \cT_{C_\fq^{\rm max}(h)}
= \big( (\cE(\g_r)\cup\{0\}) \times \cE(\g_s)\big)
\cap \cT_{C_{\fq_s}^{\rm max}(h_s)}\subeq  (\cE(\g_r)\cup\{0\}) \times \cO_{h_s}.\] 
This shows that, for any Euler element
$k \in \g$ with $W_M^+(k) \not=\eset$ we must have $k_s \in \cO_{h_s}$,
but there is no restriction on the Riemannian component
$k_r \in \cE(\g_r)$. 
\end{rem}

\section{Connectedness of the positivity domain} 
\mlabel{subsec:7.3}

In this section we show that, if $G \cong \Inn(\g)$ is the adjoint group,
then the positivity domain $W_M^+(h)$ is connected.
This contrasts the situation for compactly causal symmetric spaces,
where wedge regions are in general not connected.
A typical example is anti-de Sitter spacetime (cf.~\cite[Lemma~11.2]{NO23b}).

\begin{thm} {\rm(Connectedness of positivity domains)} 
\mlabel{thm:connwedgedom}  
Suppose that $(\g,\tau,C)$ is an irreducible simple ncc
symmetric Lie algebra with $C = C^{\rm max}_\fq$ and
the causal Euler element $h \in C^\circ \cap \fq_\fp$.
Let $M = G/H$ for $G = \Inn(\g)$ and $H = K^h \exp(\fh_\fp)$. 
Then the positivity domain $W_M^+(h)$ is connected. 
\end{thm}

\begin{prf} From Theorem~\ref{thm:x.1}
  we derive that 
\[ G^+(h) := \{ g \in G \: \Ad(g)^{-1}h \in \cT_C \}
= G^h \exp(\cD_\g) H_e, \]
and this leads with Lemma~\ref{lem:wm+} to
\[ W_M^+(h) =  G^+(h).eH
  = G^h \exp(\cD_\g).eH.\]
Since $G^h$ has at most two connected components, this set is
either connected or has two connected components
(\cite[Thm.~7.8]{MNO23}). 
As $G^h = K^h \exp(\fq_\fp)$, we have
$G^h = K^h G^h_e$, and
$\Ad(K^h)$ preserves the open unit ball in $\g_1(h)$. We thus derive 
from $K^h  = H_K$: 
\begin{align*}
 W_M^+(h)
 &  = G^h \exp(\cD_\g).eH   = G^h_e K^h\exp(\cD_\g).eH
 = G^h_e \exp(\cD_\g) K^h.eH = G^h_e \exp(\cD_\g).eH,
\end{align*}
which is connected. 
\end{prf}

\begin{cor} \mlabel{cor:equal} $W(\gamma) = W_M^+(h)$.   
\end{cor}

\begin{prop} \mlabel{prop:GW} {\rm(The stabilizer group of the observer
  domain)} 
  If $\g = \g_s$, then $G^h$
  coincides with the
  stabilizer group
  \[ G_{W(\gamma)}:= \{ g \in G \: g W(\gamma) = W(\gamma)\} \]
  of the observer domain $W(\gamma)  \subeq M = G/H$.
\end{prop}

\begin{prf} We work with the subset realization of 
  $M = G/H$ as
$M_{\cD_+} = \{ g.\cD_+ \: g \in G \}$ 
  from Proposition~\ref{prop:setreal}. 
  Then 
  \[ W(\gamma) = \{ g \cD_+ \: 0 \in g.\cD_+, g.\cD_+\
    \mbox{ bounded in } \exp(\g_1(h)).P^-\}\]
  (Lemma~\ref{lem:wgammaprops}(c)).  
  Since $\exp(\R h)$ acts on $\exp(\g_1(h))$ by dilations, it follows that
    \begin{equation}
    \label{eq:intersect}
    \bigcap_{g\cD_+\in W(\gamma)} g\cD_+ = \bigcap_{t \in \R} e^t \cD_+ = \{eP^-\}.
    \end{equation}
    Therefore $g W(\gamma) = W(\gamma)$ for the action of $g$ on
    $G/H \subeq \bP(G/P^-)$ implies that $g$ preserves
    the intersection $\{eP^-\}$ of all subsets contained in $W(\gamma)$.
    This shows that $g$ fixes $eP^-$, so that $g \in P^-$.

    Next we recall that the involution $\tau_M$ on $M$ defined by 
    $\tau_M(gH) = \tau(g)H$ leaves $W(\gamma)$ invariant
    (Lemma~\ref{lem:tauWinv}), and this leads to
    \[ G_{W(\gamma)} = \tau(G_{W(\gamma)}) \subeq P^- \cap \tau(P^-) = P^- \cap P^+ = G^h.
    \qedhere\] 
\end{prf}

The preceding proposition shows that the set
$\cW = \cW(M)$ of wedge domains in $M = G/H$ coincides with
\begin{equation}
  \label{eq:wedgespace}
  \cW = G.W(\gamma) \cong G/G^h \cong \cO_h.
\end{equation}
In particular, it is a symmetric space.
Recall that, by Corollary~\ref{cor:equal},
the observer domain coincides with the positivity domain $W_M^+(h)$. 

\section{KMS wedge regions}

With the structural results obtained so far,
we have good control over the positivity
domains $W_M^+(h)$ in ncc symmetric spaces $M = G/H$.
So one may wonder if they also have an interpretation in terms of a KMS
like condition. In \cite{NO23b}, this has been shown for modular
flows with fixed points, using such a fixed point as a base point.
In this section we extend the characterization of the wedge
domain $W$ in terms of a geometric KMS condition to general
ncc spaces. 

To simplify references, we list our assumptions and the relevant
notation below: 
\begin{itemize}
\item $\g$ is simple,
\item  $G= \Inn(\g) \subeq G_\C = \Inn(\g_\C)_e$
  (by (GP) and (Eff), \cite[Lemma~2.12]{NO23b})
\item $\sigma \: G_\C \to G_\C$ denotes the complex conjugation with
  respect to $G$. 
\item $H = G^c \cap G$, where $G^c = (G_\C^{\oline\tau})_e$
  and $K_\C \subeq G_\C^\theta$ is an open subgroup. Note that
  $H \subeq G^\tau$.  
\item $\Xi = G.\Exp_{eK}(i\Omega_\fp) \subeq G_\C/K_\C$.
\item $H_\C \subeq G^\tau_\C$ is open with $G \cap H_\C = H$ (see \S 5),
  so that $M = G/H \into G_\C/H_\C$.
\item $\tau_h^G(H_\C) = H_\C$ for the holomorphic involution
  of $G_\C$ integrating the complex linear extension of~$\tau$.
\item $\sigma(H_\C) = H_\C$ for the conjugation of $G_\C$ with respect to $G$.
\item  $\kappa_h = e^{-\frac{\pi i}{2} \ad h}$ integrates to the
    automorphism $\kappa_h^G(g) = \exp\big(-\frac{\pi i}{2}h\big)g \exp\big(\frac{\pi i}{2}h\big)$ of $G_\C$. 
\end{itemize}

 Note that
\begin{equation}
  \label{eq:kappa2}
  \tau_h^G := (\kappa_h^G)^2
\end{equation}
is a holomorphic involutive automorphism of $G_\C$ inducing $\tau_h$ on
the Lie algebra~$\g$.

Let
\[ \Xi  := G.\Exp_{eK}(i\Omega_\fp) \subeq G_\C/K_\C \]
be the crown of $G/K$. The involution $\tau_h$ on $G$ preserves $K$,
hence induces an involution on $G/K$, and we extend it to an antiholomorphic
involution $\oline\tau_h$ on $G_\C/K_\C$.
The canonical map $G \times_K i \Omega_\fp \to \Xi$ is a diffeomorphism
(\cite[Prop.~4.7]{NO23b}) 
and
\[ \oline\tau_h(g.\Exp(ix)) = \tau_h(g).\Exp(-i \tau_h(x)) \]
implies that
\begin{equation}
  \label{eq:realcrown}
 \Xi^{\oline\tau_h}
= G^{\tau_h}.\Exp(i \Omega_\fp^{-\tau_h})
= \exp(\fq_\fp).\Exp(i \Omega_{\fh_\fp}) \
\cong G^h_e \times_{K^h_e} i\Omega_{\fh_\fp}.
\end{equation}
(see the proof of \cite[Thm.~6.1]{NO23b} for details). 
This describes the fixed point as a ``real crown domain'' of the
Riemannian symmetric space $(G/K)^{\tau_h} = \Exp(\fq_\fp)$. 

For an open subgroup $H_\C \subeq G^\tau_\C$
(where $\tau$ denotes the holomorphic involution)
with $G \cap H_\C = H$, we obtain an embedding
$M = G/H \into G_\C/H_\C$. Then the stabilizer
of
\[ m_K := \Exp_{eH}\Big(\frac{\pi i}{2} h\Big)
  = \exp\Big(\frac{\pi i}{2} h\Big) H_\C \in G_\C/H_\C \]
coincides with $K$, so that $G.m_K \cong G/K$
(\cite[Thm.~5.4]{NO23b}). Accordingly, 
\[ K_\C := (\kappa_h^G)^{-1}(H_\C) \]
is an open subgroup of $G_\C^\theta$ that coincides with the stabilizer~$G_\C^{m_K}$.
In this sense $G_\C/H_\C \cong G_\C/K_\C$,
but with different base points $m_H := eH_\C$ and~$m_K$.
Recall that
$\tau = e^{\pi i \ad h}\theta = e^{\frac{\pi i}{2} \ad h}\theta e^{-\frac{\pi i}{2} \ad h}$
implies $\theta = (\kappa_h^G)^{-1} \tau \kappa_h^G$.
The invariance of $H_\C$ under $\tau_h^G$ implies that
\[ H_\C = (\kappa_h^G)^{-1}(K_\C),\]
so that $K_\C$ and $H_\C$ are exchanged by the order-$4$ automorphism
$\kappa_h^G$ and invariant under $\tau_h^G$.

As $\tau_h^G$ commutes with $\kappa_h^G$, it also leaves $K_\C$ invariant.
Moreover, $\sigma \kappa_h^G \sigma = (\kappa_h^G)^{-1}$ entails
\[ \sigma(K_\C) = \kappa_h^G(\sigma(H_\C)) = \kappa^G_h(H_\C) = K_\C.\]
Therefore the antiholomorphic extension
$\oline\tau_h^G$ also preserves $K_\C$ and
induces on $G_\C/K_\C \cong G_\C/H_\C$ an 
antiholomorphic involution $\oline\tau_h$ fixing the base point~$m_K$
with stabilizer~$K_\C$.
Then
\[ m_H' := \oline\tau_h(m_H)
  =\oline\tau_h\Big(\exp(-\frac{\pi i}{2}h\Big).m_K\Big)
  =\exp\Big(\frac{\pi i}{2}h\Big).m_K 
  =\exp(\pi i h).m_H,\]
may be different from $m_H$.

\begin{rem} The condition $m_H= m_H'$ is equivalent to
  $\exp(\pi i h) \in H_\C$. Note that $e^{\pi i \ad h} =  \tau_h \in \Aut(\g_\C)$
  is an involution that commutes with $\tau$, so that the choice of
  $H_\C$ determines whether $\exp(\pi i h)$ is contained in $H_\C$ or not.

For $\g = \fsl_2(\R)$, $G = \Inn(\g)$,
    and $h = \frac{1}{2}\pmat{1 & 0 \\ 0 & -1}$,
        we obtain  on $\SL_2(\R)\subeq \SL_2(\C)$ the involution 
    \[ \tau\pmat{a & b \\ c & d} = \pmat{d & c \\ b & a}. \]
    For $g \in \SL_2(\C)$, the condition
    $\tau\Ad(g)\tau = \Ad(g)$ is equivalent to
    $\tau(g)g^{-1} \in \ker(\Ad) = \{\pm \1\}$.
    As
    \[ e^{\pi i \ad h} = \Ad(\exp \pi ih)
      = \Ad\pmat{ i & 0 \\ 0 & -i} \quad \mbox{ and } \quad
      \tau\pmat{ i & 0 \\ 0 & -i} \pmat{ i & 0 \\ 0 & -i}^{-1}  = - \1,\]
    it follows that
    \[ \tau_h =  e^{\pi i \ad h} \in G_\C^\tau \setminus (G_\C^\tau)_e.\]
    In particular, $K_\C$ and $H_\C$ have two connected
    components in $G_\C \cong \PSL_2(\C)$.

    In $G \cong \PSL_2(\R)$, a similar argument shows that
    $\theta = \Ad\pmat{ 0 & -1 \\ 1 & 0} \in G^\tau \setminus G^\tau_e$.
    So $G^\tau$ also has two connected components, but only
    its identity component $G^\tau_e$ acts causally on $\fq$.
    Therefore $H = G^\tau_e$, but for $H_\C$ we have two choices,
    $G_\C^\tau$, or its identity component. 
\end{rem}

Comparing with the arguments in \cite[Lemma~6.3]{NO23b}, where
$\alpha_{\pi i}= \tau_h$ on $M$, 
we have to be more careful in the present context.
Here $\oline\tau_h$ restricts to a map
\[ M = G.m_H \to M' := G.m_H' = \exp(\pi i h).M,\]
and these two copies of $G/H$ may not be identical.
However, the antiholomorphic map
\[ \sigma_M :=  \alpha_{\pi i} \circ \oline\tau_h \]
maps $M$ to itself, fixes the base point $m_H$ and commutes with the
$G$-action. Hence it fixes $M$ pointwise and describes a ``complex
conjugation'' with respect to~$M$. In particular, the two maps
\[ \oline\tau_h \: M \to M' \quad \mbox{ and } \quad 
\alpha_{\pi i} \:   M \to M' \] 
coincide on $M$. 

We define the {\it KMS wedge domain} 
\begin{equation}
  \label{eq:wkms}
  W^{\rm KMS} := \{  m \in M \:  \alpha_{it}(m) \in \Xi \mbox{ for }
  0 < t < \pi\}.
\end{equation}

\begin{thm} \mlabel{thm:kms} $W^{\rm KMS} = W_M^+(h)_{eH}$.     
\end{thm}

\begin{prf} \nin ``$\subeq$'': 
For $z \in \C$ and $p \in M$, we first observe that 
\[ \oline\tau_h(\alpha_z(p))
= \alpha_{\oline z}(\oline\tau_h(p)) 
= \alpha_{\oline z}\alpha_{\pi i}(p) 
= \alpha_{\pi i + \oline z}(p).\]
For $z = \frac{\pi i}{2}$, we thus obtain
$\alpha_{\pi i/2}(p) \in M_\C^{\oline\tau_h}$.
We conclude  that
\begin{equation} \label{eq:8.1}
  \alpha_{\frac{\pi i}{2}}(W^{\rm KMS}) \subeq  \Xi^{\oline\tau_h}
      = \exp(\fq_\fp).\Exp(i \Omega_{\fh_\fp}).
\end{equation}
Hence 
\[ W^{\rm KMS} 
\subeq  \kappa_h\big(\Xi^{\oline\tau_h}\big)
= G^h_e.\Exp_{eH}(\kappa_h(i \Omega_{\fh_\fp})),\]
where
\[ \kappa_h(i\Omega_{\fh_\fp}) = \Big\{ x \in \fq_\fk \:
\|\ad x\| < \frac{\pi}{2}\Big\} =: \Omega_{\fq_\fk}.\]
This suggest to define a ``polar wedge domain'' as
\[ W_M^{\rm pol}(h) := G^h_e.\Exp(\Omega_{\fq_\fk})\subeq M.\] 
We actually know from Theorem~\ref{thm:pos-polar}, that this is
the connected component $W = W_M^+(h)_{eH} \subeq W_M^+(h)$
containing the base point. 
We thus obtain
\begin{equation}
  \label{eq:8.2}
  W^{\rm KMS}  \subeq W = W_M^+(h)_{eH}.
\end{equation}

\nin ``$\supeq$'': To see that $W_M^+(h)_{eH}\subeq W^{\rm KMS}$,
we first recall from the first part of the proof that 
\[ W_M^+(h)_{eH} = \kappa_h\big(\Xi^{\oline\tau_h}\big)
  = G^h_e.\Exp_{eH}(\kappa_h(i \Omega_{\fh_\fp}))
  = G^h_e.\alpha_{-\pi i/2}(\Exp_{eK}(i\Omega_{\fh_\fp})).\]
To see that this domain is contained in the $G^h_e$-invariant domain
$W^{\rm KMS}\subeq M$, we thus have to show
that, for $x \in \Omega_{\fh_\fp}$, we have 
\[ \alpha_{it}.\Exp_{eK}(ix) \in \Xi \quad \mbox{ for } \quad  |t| < \pi/2.\] 
Let $\ft_\fq \subeq \fq_\fk$ is a maximal abelian subspace
(they are all conjugate under $(H_K)_e$). Then
$\fa_\fh := i\kappa_h(\ft_\fq) \subeq \fh_\fp$ is also maximal abelian and
$\Omega_{\fh_\fp} = e^{\ad \fh_\fk}.\Omega_{\fa_\fh}$. 
So it suffices to show that,
for $x \in \Omega_{\fa_\fh}$ and $|t| < \pi/2$, we have
$\alpha_{it}.\Exp_{eK}(ix) \in \Xi.$ 
  By Proposition~\ref{prop:testing}, $\ft_\fq$ is contained in a 
  $\tau$-invariant subalgebra 
  $\fs \cong \fsl_2(\R)^s$, where $\R h + \fs$ is generated by $h$ and $\ft_\fq$
  and $h = h_0 + h_1 + \cdots + h_s$, where
  $h_j$, $j = 1,\ldots, s$, is an Euler element in a  simple ideal
  $\fs_j \cong \fsl_2(\R)$ of $\fs$.
  Then $\fa_\fh = i \kappa_h(\ft_\fq) \subeq \fa$ is spanned
  by $s$ Euler elements $x_1, \ldots, x_s$ and
  \[ \Omega_{\fa_\fh} = \Big\{ \sum_{j = 1}^s t_j x_j \:
    (\forall j)\ |t_j| < \pi/2\Big\}.\] 
  Let $S := \la \exp \fs \ra$ and
  $\Xi_S := S.\Exp(i(\Omega_\fp \cap \fs)) \subeq \Xi$.
  Then the discussion in Remark~\ref{rem:app-deSitter}
  implies that, for $|t| < \pi/2$ and $x = \sum_j t_j x_j \in \Omega_{\fa_\fh}$, we
  have $\alpha_{it}(\Exp_{eK}(ix)) \in \Xi_S \subeq \Xi.$ 
  \end{prf}

  The preceding proof implies in particular the following interesting
  observation: 
\begin{cor} For every $m \in \Xi^{\oline\tau_h}$, we have
  $\alpha_{it}(m) \in \Xi$ for $|t| < \pi/2$, so that the orbit map
  $\alpha^m$ extends to a holomorphic map $\cS_{\pm \pi/2} \to \Xi$.
\end{cor}

\begin{cor} 
  $\alpha_{\frac{\pi i}{2}} \: W^{\rm KMS} \to  \Xi^{\oline\tau_h}$ is a diffeomorphism
  that induces an equivalence of fiber bundles
      \[ W^{\rm KMS}     \cong G^h_e \times_{K^h_e} \Omega_{\fq_\fk} 
    \to G^h_e \times_{K^h_e} i\Omega_{\fh_\fp} \cong \Xi^{\oline\tau_h}.\]
\end{cor}

\begin{prf} Theorem~\ref{thm:kms} implies in particular that 
  $\alpha_{\frac{\pi i}{2}} \: W^{\rm KMS} \to  \Xi^{\oline\tau_h}$ is bijective. Since $W^{\rm KMS} = W_M^+(h)_{eH}$ is an open subset of
  $M$ and $\Xi^{\oline\tau_h}$ an open subset of
  $M_\C^{\oline\tau_h}$, it actually is a diffeomorphism. 
  The second assertion follows from the fact that it
  commutes with the action of the subgroup $G^h_e$.
\end{prf}

\appendix

\section{Irreducible ncc symmetric Lie algebras}
\mlabel{app:classif}

The following table lists all irreducible non-compactly 
causal symmetric Lie algebras $(\g,\tau)$ according to the following types: 
\begin{itemize}
\item Complex type: $\g= \fh_\C$ and $\tau$ is complex conjugation with 
  respect to $\fh$. In this case $\g^c \cong \fh^{\oplus 2}$, so that
  $\rk_\R(\g^c) = 2 \rk_\R(\fh)$. 
\item Cayley type (CT): $\tau = \tau_{h_1}$ for an Euler element $h_1 \in \fh$.
  Then $\rk_\R(\g^c) = \rk_\R(\g) = \rk_\R(\fh)$.
\item Split type (ST):  $\tau \not=\tau_{h_1}$ for all $h_1 \in \fh \cap \cE(\g)$ 
and $\rk_\R \fh = \rk_\R \g^c$:
\item Non-split type (NST): 
$\tau \not=\tau_{h_1}$ for all $h_1 \in \fh \cap \cE(\g)$ 
and $\rk_\R \fh =\frac{ \rk_\R \g^c}{2}$:
\end{itemize}

In the table we write $r = \rk_\R(\g^c)$ and $s = \rk_\R(\fh)$.
Further $\fa \subeq \fp$ is maximal abelian of dimension~$r$.
For root systems $\Sigma(\g,\fa)$ of type $A_{n-1}$,
there are $n-1$ Euler elements $h_1, \ldots, h_{n-1}$,
but for the other root systems there are less;
see \cite[Thm.~3.10]{MN21} for the concrete list.
For $1 \leq j <n$ we write $j' := \min(j,n-j)$. \\[4mm]

\hspace{-17mm}
\begin{tabular}{||l|l|l|l|l|l|l|l||}\hline
$\g$ &  $\g^c = \fh + i \fq$ & $r$ & $\fh = \g^{\tau_h\theta}$ 
& $s$ & $\Sigma(\g,\fa)$  & $h$ &  $\g_1(h)$  \\ 
\hline
\hline 
Complex type \phantom{\Big (} &&& &&  &&\\
\hline 
$\fsl_n(\C)$ & $\su_{j,n-j}(\C)^{\oplus 2}$  \phantom{\Big(}& $2j$
& $\su_{j,n-j}(\C)$ & $j'$ & $A_{n-1}$ & $h_j$ & $M_{j,n-j}(\C)$  \\
 $\sp_{2n}(\C)$  & $\sp_{2n}(\R)^{\oplus 2}$ & $2n$ &  $\sp_{2n}(\R)$ 
& $n$ & $C_n$ & $h_n$ & $\Sym_n(\C)$   \\
$\so_n(\C), n > 4$ &  $\so_{2,n-2}(\R)^{\oplus 2}$ & 4 & $\so_{2,n-2}(\R)$ 
& 2 & $D_{[\frac{n}{2}]}, B_{[\frac{n}{2}]}$  &$h_1$ & $\C^{n-2}$   \\
$\so_{2n}(\C)$ & $\so^*(2n)^{\oplus 2}$ & $2[\frac{n}{2}]$ & $\so^*(2n)$&
$[\frac{n}{2}]$ &$D_{n}$ & $h_{n-1}, h_{n}$ & $\Skew_{n}(\C)$      \\
 $\fe_6(\C)$ &  $(\fe_{6(-14)})^{\oplus 2}$ & $4$ &$\fe_{6(-14)}$ & $2$ &
$E_6$ & $h_1, h_6$ & $M_{1,2}(\bO)_\C$    \\
 $\fe_7(\C)$ &  $(\fe_{7(-25)})^{\oplus 2}$ & $6$ & $\fe_{7(-25)}$ & $3$ & 
$E_7$ & $h_7$ & $\Herm_3(\bO)_\C$     \\
\hline
Cayley type  \phantom{\Big (}&& &&&  &&\\
\hline 
$\su_{r,r}(\C)$ &  $\su_{r,r}(\C)$ & $r$ & $\R \oplus \fsl_r(\C)$& $r$ &
$C_r$ & $h_r$ & $\Herm_r(\C)$ \\
$\sp_{2r}(\R)$  & $\sp_{2r}(\R)$ & $r$ & $\R \oplus \fsl_r(\R)$ & $r$ & 
$C_r$ & $h_r$ & $\Sym_r(\R)$   \\
$\so_{2,d}(\R), d> 2$ &  $\so_{2,d}(\R)$ & $2$ &  
$\R \oplus \so_{1,d-1}(\R)$  & $2$ &   $C_2$ &$ h_2$ & $\R^{1,d-1}$  \\
$\so^*(4r)$ &  $\so^*(4r)$ & $r$ & 
$\R \oplus \fsl_r(\H)$& $r$ & $C_r$ & $h_r$ & $\Herm_r(\H)$ \\
 $\fe_{7(-25)}$ & $\fe_{7(-25)}$ & $3$ & 
$\R \oplus \fe_{6(-26)}$ & $3$  & $C_3$ & $h_3$ & $\Herm_3(\bO)$  \\
\hline
Split type \phantom{\Big (} &&&&&  &&\\
\hline 
$\fsl_{n}(\R)$ &$\su_{j,n-j}(\C)$ & $j'$ &  $\so_{j,n-j}(\R)$ &
$j'$ & 
$A_{n-1}$ & $h_j$ & $M_{j,n-j}(\R)$ \\
$\so_{n,n}(\R)$ & $\so^*(2n)$ & $[\frac{n}{2}]$ & $\so_{n}(\C)$ &
$[\frac{n}{2}]$ &$D_{n}$ & $h_{n-1}, h_{n}$ & $\Skew_n(\R)$  \\
  $\so_{p+1,q+1}(\R)$ &$\so_{2,p+q}(\R)$& $2$ & 
$\so_{1,p}(\R) \oplus \so_{1,q}(\R)$ & $2$
& $ B_{p+1}\, (p<q)$  & $h_1$ &$\R^{p,q}$    \\
$p,q > 1$ && && & $D_{p+1}\, (p = q)$  & &    \\
$\fe_6(\R)$ & $\fe_{6(-14)}$& $2$ & $\fu_{2,2}(\H)$
& $2$ &  $E_6$ & $h_1, h_6$ &
$M_{1,2}(\bO_{\rm split})$   \\
$\fe_7(\R)$ & $\fe_{7(-25)}$ & $3$ &$\fsl_4(\H) = \su^*(8)$ & $3$ & $E_7$ & $h_7$ &
$\Herm_3(\bO_{\rm split})$   \\
\hline 
Non-split type \phantom{\Big (} && &&&  &&\\
\hline 
$\fsl_n(\H)$ &$\su_{2j,2n-2j}(\C)$ & $2j'$ & $\fu_{j,n-j}(\H)$
& $j'$ & $A_{n-1}$ & $h_j$ & $M_{j,n-j}(\H)$   \\
 $\fu_{n,n}(\H)$ & $\sp_{4n}(\R)$ & $2n$ & $\sp_{2n}(\C)$ & $n$ 
& $C_{n}$ & $h_n$ & $\Aherm_n(\H)$  \\
$\so_{1,d+1}(\R)$  & $\so_{2,d}(\R)$ & $2$ & $\so_{1,d}(\R)$ & $1$ 
& $A_1$ & $h_1$ & $\R^d$ \\
$\fe_{6(-26)}$ & $\fe_{6(-14)}$  &$2$  & $\ff_{4(-20)}$  & $1$ &
$A_2$ & $h_1, h_2$ &
$M_{1,2}(\bO)$  \\
\hline
\end{tabular} \\[2mm] 
{\rm Table 1: Irreducible ncc symmetric Lie algebras 
with corresponding causal Euler elements~$h \in \fa$} 

\section{Geodesics in symmetric spaces}

This appendix contains some elementary observations  concerning
  geodesics in symmetric spaces.

\begin{lem} \mlabel{lem:geodesic} 
  Let $M = G/H$ be a symmetric space with symmetric Lie algebra
  $(\g,\tau)$,
  $x \in \g$ and $y \in \fq$. Then 
\begin{equation}
  \label{eq:curveq}
 \exp(tx) H = \exp(ty) H \quad \mbox{ for all } \quad t \in \R 
\end{equation}
holds if and only if $p_\fq(x) = y$ and $[x,y] = 0$. 

In particular, $\gamma(t) := \exp(tx)H$ is a geodesic in $M$ if and only if 
$[x, \tau(x)] = 0$. 
\end{lem}

\begin{prf} The relation \eqref{eq:curveq} is equivalent to 
\[ \exp(-ty) \exp(tx) \in H \subeq G^{\tau^G} \quad \mbox{ for all } \quad t \in \R.\]  
Applying $\tau^G$, we obtain 
\[ \exp(ty) \exp(t \tau(x)) = \exp(-ty) \exp(tx), \] 
which leads to 
$ \exp(2ty)  =  \exp(tx) \exp(-t \tau(x)).$ 
Evaluating the derivative of this curve in the right trivialization of 
$T(G)$, we get 
\[ 2y = x + e^{t\ ad x}(-\tau(x))= x - e^{t\ ad x}(\tau(x))\quad \mbox{ for all } \quad t \in \R.\] 
For $t = 0$ we get $p_\fq(x) = y$, and taking derivatives in $0$ shows that 
$[x,\tau(x)] = 0.$ 

If, conversely, this condition is satisfied, then 
$x = x_\fh + x_\fq$ with $x_\fh \in \fh$ and $x_\fq \in \fq$, 
where 
\[ 0 = [x,\tau(x)] = 2 [x_\fh, x_\fq].\] 
Therefore 
\[ \exp(tx)H 
  = \exp(t x_\fq) \exp(t x_\fh) H 
= \exp(t x_\fq) H  = \Exp_{eH}(t x_\fq) \] 
is a geodesic in $M$. 
\end{prf}

The following lemma
provides important information on the subset~$M^x$.

\begin{lem} \mlabel{lem:mx}
Let $x \in \g$ and write 
\[ M^x := \{ gH \in M \: \Ad(g)^{-1}x \in \fq \}.\] 
Then $M^x$ is a submanifold of $M$ which is 
invariant under the action of $G^x$, and the orbits or 
$G^x_e$ are the connected components of $M^x$. 
\end{lem}

\begin{prf} Let $m_0 = g_0 H \in M^x$ and $x_c := \Ad(g_0)^{-1} x$. 
For $y \in \fq$ we have 
\[ \Exp_{m_0}(g_0.y) = g_0.\Exp_{eH}(y) = g_0 (\exp y) H 
= \exp(\Ad(g_0)y).m_0\] 
and 
\[ \Ad(g_0\exp(y))^{-1} x =  e^{-\ad y} x_c 
= \cosh(\ad y)x_c - \sinh(\ad y) x_c
= \underbrace{\cosh(\ad y)x_c}_{\in \fq}
- \underbrace{\frac{\sinh(\ad y)}{\ad y} [y,x_c]}_{\in \fh}. \] 
Let $U \subeq \fq$ be a $0$-neighborhood for which 
$\Exp_{eH}\res_U$ is a diffeomorphism onto an open subset of $M$ 
and the spectral radius of $\ad y$ is smaller than $\pi$ for $y \in U$. 
Then $\frac{\sinh(\ad y)}{\ad y} \: \fh \to \fh$ is invertible. 
With the above formula, we thus conclude for $y \in U$ that 
$\Exp_{m_0}(g_0.y) \in M^x$ is equivalent to $[y,x_c] = 0$, 
which is equivalent to $\Ad(g_0)y \in \g^x$. 
This shows that $M^x$ is a submanifold of $M$. 

As $\Exp_{m_0}(g_0.y) = \exp(\Ad(g_0)y).m_0 \in \exp(\g^x).m_0$, 
it further follows that the the orbit 
of $m_0$ under the connected group $G^x_e$ contains a neighborhood 
of $m_0$. This shows that 
the orbits of $G^x_e$ in $M^x$ are connected open subsets, hence 
coincide with its connected components. 
\end{prf}

\begin{rem}  For $x \in \fq$ the centralizer 
$\g^x$ is $\tau$-invariant, so that 
$\g^x = \fh^x \oplus \fq^x$ and the dimension 
of the $G^x_e$-orbit through $eH$ is $\dim \fq^x$. 
We have 
\[  M^x_{eH} \cong  G^x_e/ (H \cap G^x_e), \] 
and Lemma~\ref{lem:mx} 
shows that 
the geodesic $\Exp_{eH}(\R x)$ is central in the symmetric space 
$M^x_{eH}$ in the sense that its tangent space 
$\R x$ is central in the Lie algebra $\g^x$ (cf.\ \cite{Lo69}). 
\end{rem}

\begin{lem} \mlabel{lem:3.11}
For $y \in \fq$, the equality $M^y = G^y.eH$ is equivalent to 
\begin{equation}
  \label{eq:cond2}
 \cO_y \cap \fq = \Ad(H) y.
\end{equation}
\end{lem}

\begin{prf} As $y \in \fq$, the base point $eH$ is contained in $M^y$, 
and thus $G^y.eH \subeq M^y$. 
So the equality $M^y = G^y.eH$ means that $M^y \subeq G^y.eH$, 
i.e., $\Ad(g)^{-1}y \in \fq$ implies $gH \in G^y.eH$, resp., 
$g \in G^y H$. This  in turn is equivalent to 
$\Ad(g)^{-1}y \in \Ad(H)y$. 
\end{prf}

\section{Lawson's Theorem}
\mlabel{sec:law}

Let $(G,\tau^G,H,C)$ be a causal symmetric Lie group,
i.e., $\tau^G$ is an involutive automorphism of $G$,
$H \subeq G^{\tau^G}$ an open subgroup and
$C \subeq \fq$ a hyperbolic
pointed generating $\Ad(H)$-invariant closed convex cone.
We write $\g = \fh \oplus \fq$ for the corresponding decomposition of $\g = \L(G)$.

According to \cite[Lemma~2.3]{La94}, $\exp\res_{C}$ is injective if and only if
$\Gamma_Z := \fz(\g) \cap \exp_G^{-1}(e)$ satisfies
\[ \Gamma_Z \cap \fq = \Gamma_Z \cap (C-C) = \{0\}.\]
If $\fz(\g) \subeq \fq$, this condition is satisfied if and only if
$\exp\res_{\fz(\g)}$ is injective. This condition is always satisfied
if $\g$ is semisimple because $\fz(\g) = \{0\}$ in this case.

Suppose that $\Gamma_Z = \{0\}$. By \cite[Lemma~2.4]{La94}, 
$\exp\res_C$ is a homeomorphism onto a closed subset
of $G$ if and only if, 
for no non-zero $x \in C \cap \fz(\g)$, the subgroup $\oline{\exp(\R x)}$
is compact. By \cite[Thm.~3.1]{La94}, this in turn is equivalent to the
polar map
\[ \Phi \: C \times G^{\tau^G} \to \exp(C) G^{\tau^G}\]
being a homeomorphism onto a closed subset of~$G$.
\cite[Thm.~3.1]{La94} further shows that
$\exp(C) G^{\tau^G}$ is a subsemigroup of $G$.
If $G$ is $1$-connected, then the subgroup $G^{\tau^G}$
is connected and $Z(G)$ is simply connected, so that
all requirements from above are satisfied
(\cite[Cor.~3.2]{La94}).

\begin{thm} {\rm(Lawson's Theorem)} \mlabel{thm:lawson}
  Let $(G, \tau^G, H, C)$ be a non-compactly causal reductive
  symmetric Lie group. Suppose that $\fz(\g) \subeq \fq$
  and that $\Gamma_Z = \{0\}$.
  Then $S := \exp(C)H$ is a closed subsemigroup of $G$
  with Lie wedge $\L(S) = \fh + C$.   
\end{thm}

\begin{prf} Our assumption implies that
  $\exp\res_{\fz(\g)} \: \fz(\g) \to Z(G)_e$ is bijective,
hence a diffeomorphism onto the closed subgroup $Z(G)_e$.
It follows in particular that $\exp(\R x) \cong \R$ is non-compact
for each non-zero
$x \in \fz(\g)$. Therefore the polar map $\Phi$ is a homeomorphism
onto a closed subset and the remaining assertions
follow from \cite[Thm.~3.1]{La94}.
  
\end{prf}

\begin{rem}
(a) If $G$ is reductive, then $G = (G,G)_e Z(G)_e$
and if $x \in \fz(\g)$ satisfies $\exp z \in (G,G)_e$, then
$\exp z \in Z((G,G)_e)^{-\tau^G}$, which is a discrete group.
We shall see below that this group may be infinite, even  if
$\Gamma_Z = \{0\}$. 

\nin (b) If $M = G_\C/G$ is of complex type 
  and $G$ is hermitian, then $Z(G_\C)$ is finite.

\nin (c) If $M$ is of non-complex type and irreducible, then $\g^c$
  is simple hermitian with $\fz(\fk^c) = \R i h$, where
  $h \in \cE(\g)$ is a causal Euler element. 
  If $Z(G)$ is infinite, then $\fg$ is also hermitian, hence of tube
  type because it contains an Euler element   (\cite{MN21}).
 Then all Euler elements
  in $\g$ are conjugate and this implies that $(\g,\tau)$ is of Cayley type. 
  So $\fz(\fk) \subeq \fq_\fk$ and thus $Z(G)^{-\tau}$ is infinite if
  $G$ is simply connected. This shows that it is possible that
  $(G,G)_e \cap Z(G)_e$ is infinite. 

  A concrete example is the group
  \[ G := (\tilde \SL_2(\R) \times \R)/D, \]
  where $D \subeq Z(\tilde\SL_2(\R)) \times \R\cong \Z \times \R$
  is the graph   of a non-zero homomorphism
  $\gamma \: Z(\tilde\SL_2(\R)) \to \R$.
  Then $Z(G) \cong \R$ and $Z(G) \cap (G,G) \cong Z(\tilde\SL_2(\R)) \cong \Z$. 
\end{rem}

\begin{rem} Suppose that $h_0 \in \fh \cap \cE(\g)$ is such that 
$-\tau_{h_0}(C) = C$, 
    then $C \cap \fz(\g)$ is contained in $C \cap -C = \{0\}$.
    Therefore the condition on
    $C \cap  \fz(\g)$ in Lawson's Theorem (cf.\ Appendix~\ref{sec:law})
    is satisfied.  
\end{rem}

\section{de Sitter space} 
\mlabel{subsec:dS1} 

In this appendix we collect some concrete observations
concerning de Sitter space $\dS^d$, which is an important
example of a non-compactly causal symmetric space. Some facts on
$2$-dimensional de Sitter space are used in particular in some
of our proofs to verify the corresponding assertions for
$\g = \fsl_2(\R)$. 

In $(d+1)$-dimensional Minkowski space $\R^{1,d}$, we write the 
Lorentzian form as 
\[ \beta(x,y) = x_0y_0 - \bx \by 
\quad \mbox{ for } \quad x = (x_0,\bx), y = (y_0, \by).\] 
We consider $d$-dimensional de Sitter space 
  \[ M := \dS^d := \{ x = (x_0, \bx) \in \R^{1,d} 
\: x_0^2 - \bx^2 = -1\},\] 
$G = \SO_{1,d}(\R)_e$ and the Euler element
$h \in \so_{1,d}(\R)$, defined by
\[ h.(x_0, \ldots, x_d) = (x_1, x_0, 0,\ldots, 0).\]
It generates the Lorentz boost in the $x_0$-$x_1$-plane. 
The fixed point set of the modular flow in $M = \dS^d$ is 
\[  M^\alpha = M \cap \Spann\{\be_2,\ldots, \be_d\}
  = \{(0,0,x_2,\ldots, x_d) \: x_2^2 + \cdots + x_d^2 = 1\} \cong \bS^{d-2}.\]
This submanifold is connected for $d > 2$ and consists of
two points for $d = 2$. 
The corresponding wedge domain is the connected subset 
\[ W_M^+(h) = M \cap W_R = \{ x \in \dS^d \: x_1 > |x_0|\}.\]

By \cite[Prop.~D.3]{NO23b}, the timelike geodesics of~$M$ of velocity $1$ 
take the form 
\[ \gamma(t) = \Exp_x(tv)=
 \cosh(t) x + \sinh(t) v, \qquad \beta(v,v) = 1, \beta(x,v) = 0, 
\beta(x,x) = 1 \] 
whereas the trajectories of the modular 
flow are 
\[ \alpha_t(x) = e^{th}x 
= (\cosh(t) x_0 + \sinh(t) x_1, \cosh(t) x_1 + \sinh(t) x_0, x_2, 
\ldots, x_d).\]
Comparing  both expressions leads for $h$-modular geodesics to the 
conditions 
\[ x_2 = \cdots = x_d = 0 \quad \mbox{ and }\quad 
v = h.x = (x_1, x_0,0,\ldots, 0).\] 
Therefore exactly two orbits of the modular flow are 
timelike geodesics. If we also ask for the geodesic to be 
positive with respect to the causal structure, then 
$x_1 > 0$ determines the geodesic uniquely. 

We infer from \cite[Prop.~D.3]{NO23b} that 
\[ \Exp_{\be_2}(t \be_1) = \cos(t) \be_2 + \sin(t) \be_1 \] 
is a closed space-like geodesics. For $0 < t < \pi$, its values are contained 
in $W_M^+(h)$, and this geodesic arc connects the two fixed points 
$\be_2$ to $-\be_2$ of the modular flow.

\begin{rem} \mlabel{rem:app-deSitter}
  In addition to
  $h$, we also consider the Euler elements defined by
\[   h_d(x_0,\ldots, x_d) = (x_d,0,\ldots, 0, x_0).\]
The involution corresponding to $h$ acts on $\R^{1+d}$ by 
\[ \tau_h(x_0, x_1, \ldots, x_d) = (-x_0, -x_1,x_2, \ldots, x_d),\]
and its antilinear extension acts on $\C^{1+d}$ by
\[ \oline\tau_h(z_0, z_1, \ldots, z_d) = (-\oline{z_0},
  -\oline{z_1},\oline{z_2}, \ldots, \oline{z_d}).\]
Note that, in $\g$, we have $\tau_h(h) = h$ and $\tau_h(h_d) = - h_d$,
so that $h \in \fq_\fp$ and $h_d \in \fh_\fp$.

In $\C^{1+d}$, we consider the domain
\[ \Xi := \{ z = x + i y \in\C^{1+d} \: y_0 > 0,
  y_0^2 > y_1^2 + \cdots + y_d^2 \}.\]
On $\Xi$ the antiholomorphic 
involution $\oline\tau_h$
has the fixed point set
\begin{align*}
  \Xi^{\oline\tau_h}
&   = \Xi \cap (i\R^2 \oplus \R^{d-1}) \\ 
  &  = \{ (i x_0, i x_1, x_2, \ldots, x_d) \:
    x_0 > |x_1|, - x_0^2 + x_1^2 - x_2^2 - \ldots - x_d^2= -1\}\\
  &  = \{ (i x_0, i x_1, x_2,\ldots, x_d) \:
    x_0 > |x_1|,  \underbrace{x_0^2 -x_1^2}_{>0}  + x_2^2 + \cdots
    + x_d^2= 1\}.
\end{align*}
It follows in particular that
\[ x_0^2 - x_1^2 \in (0,1].\] 
The analytic extension of  the modular flow $(\alpha_t)_{t \in \R}$
acts on $\Xi$ by 
\[ \alpha_{it}(z_0, \ldots, z_d)
  = (\cos t \cdot z_0 + i \sin t \cdot z_1,
  i \sin t \cdot z_0 + \cos t \cdot z_1, z_2, \ldots, z_d).\]
Starting with a $\oline\tau_h$-fixed element
$z = (ix_0, ix_1, x_2, \ldots, x_d)$
 in $\Xi$, this leads to
 \[ \alpha_{it}(ix_0, ix_1, x_2,\ldots, x_d)
  = (\cos t \cdot ix_0 -  \sin t \cdot x_1,
  - \sin t \cdot x_0 + \cos t \cdot i x_1, x_2, \ldots, x_d)\]
 with imaginary part 
 \[ (x_0 \cos t, x_1 \cos t, 0,\ldots, 0), \]
 so that we obtain for $|t| < \pi/2$ that
 \[ 
   |x_0 \cos t| = x_0 \cos t > |x_1| \cos t, \]
 which implies that
 \begin{equation}
   \label{eq:curve-in-xi}
 \alpha_{it}(z) \in \Xi \quad \mbox{ for } \quad
   z \in \Xi^{\oline\tau_h} \quad \mbox{ and }  \quad
   |t| < \pi/2.
 \end{equation}
\end{rem}

\begin{ex}
For the special case $d = 2$, we have $\fsl_2(\R) \cong \R^{1,2}$ and
the Euler element
\[ h^0= \frac{1}{2} \pmat{1 & 0 \\ 0 & -1}
    \in \fsl_2(\R)\cong \R^{1,2} \] 
corresponds to the base point $\be_2$  (see \cite{NO23b}), so that 
\[ \fsl_2(\R) \supeq  \cO_h \cong \dS^2 \subeq \R^{1,2}.\]
Accordingly, 
\[  C = \cone(e^0, f^0), \quad 
C^c = \cone(e^0, -f^0), \quad 
 \quad \mbox{ and } \quad 
x_0 := \frac{1}{2}(e^0 - f^0) = \frac{1}{2}\pmat{ 0 & 1 \\ -1 & 0} \in C^c.\] 
For $g_t := \exp(t x_0)$ we then have 
\[ \Ad(g_{\pi/2}) h^0 = -h^1, \quad \Ad(g_{\pi/2}) h^1 = h^0 \quad \mbox{ and } 
\quad \Ad(g_\pi) h^0 = -h^0.\] 
We also note that, for $0 < t < \pi$, the Lie algebra 
element $\Ad(g_t) h^0$ corresponds to 
$\Exp_{\be_2}(t \be_1) \in W_M^+(h^0)$. Note that 
\[ g_\pi \in K^{\tau_h^G} = K^{\tau^G}.\] 
\end{ex}

\bigskip

\noindent 
Vincenzo Morinelli\\ Dipartimento di Matematica \\
   Universit\`a di Roma ``Tor Vergata'' , Italy\\
   morinell@mat.uniroma2.it
\bigskip

\noindent Karl-Hermann Neeb\\ 
Department Mathematik \\ 
Friedrich-Alexander-Universit\"at \\ 
Erlangen-N\"urnberg \\ 
Cauerstrasse 11 \\ 
91058 Erlangen, Germany\\ 
neeb@math.fau.de
\bigskip

\noindent Gestur \'Olafsson\\
Department of mathematics \\ 
Louisiana State University \\ 
Baton Rouge, LA 70803\\ 
olafsson@math.lsu.edu
\newpage

\end{document}